\documentclass[11pt]{amsart}

\usepackage{amssymb,amsmath} 
\usepackage[latin1]{inputenc}
\usepackage[francais]{babel}
\usepackage[matrix,arrow,curve]{xy}
\usepackage[mathscr]{eucal}
\usepackage{epic,eepic}
\usepackage{bbm}
\usepackage{yhmath}

\newtheorem{theoreme}[subsection]{Th\'eor\`eme}
\newtheorem{lemme}[subsection]{Lemme}
\newtheorem{remarque}[subsection]{Remarque}
\newtheorem{corollaire}[subsection]{Corollaire}
\newtheorem{prop}[subsection]{Proposition}
\newtheorem{exemple}[subsection]{Exemple}
\newtheorem{hypothese}[subsection]{Hypoth\`ese}
\newtheorem{conj}[subsection]{Conjecture}
\newtheorem{notation}[subsection]{Notation}

\newenvironment{pf}
{\medskip\noindent {\it Preuve --- \ }}
{\hfill\nobreak $\Box$ \par\bigbreak}
\def\hypo{hypoth\`ese\,}

\newcommand{\Qp}{{\mathbb Q}_p }
\newcommand{\Q}{{ \mathbb Q } }
\newcommand{\R}{{\mathbb R}}
\newcommand{\C}{{\mathbb C}}
\newcommand{\Z}{{ \mathbb Z  }}
\newcommand{\N}{{ \mathbb N  }}
\newcommand{\Qb}{\overline{\Q}}
\newcommand{\Qpb}{{\overline{\Q}_p}}
\newcommand{\AAA}{\mathbb{A}}
\newcommand{\A}{{\mathbb A}}

\newcommand{\GG}{\mathbb G}

\newcommand{\OO}{{\cal O}}

\newcommand{\cal}{\mathcal}
\newcommand{\lra}{\longrightarrow}
\newcommand{\got}{\mathfrak}
\newcommand{\ps}{\par \smallskip}
\newcommand{\ba}{\backslash}
\newcommand{\lgr}{\longrightarrow}

\newcommand{\Hom}{\text{Hom}}
\newcommand{\Ker}{{\text{Ker}\,}}

\renewcommand{\ker}{{\text{Ker}\,}}
\newcommand{\isomo}{\overset{\sim}{\rightarrow}}

\newcommand{\trace}{{\rm trace}}

\newcommand{\Gal}{{\rm Gal}}

\newcommand{\W}{{\rm W}}
\newcommand{\WD}{{\rm WD}}
\newcommand{\Supp}{{\rm Supp}}

\newcommand{\Lef}{{\rm Lef}}
\newcommand{\St}{{\rm St}}

\newcommand{\GL}{{\rm GL}}
\newcommand{\SO}{{\rm SO}}
\newcommand{\Sp}{{\rm Sp}}
\newcommand{\SL}{{\rm SL}}
\newcommand{\G}{{\rm G}}

\newcommand{\Tpt}{\Theta_{\pi,\theta}}

\newcommand{\TO}{{\rm TO}}
\newcommand{\STO}{{\rm STO}}
\newcommand{\ON}{{\rm O}}

\newcommand{\Norm}{{\cal N}}
\renewcommand{\No}{\Norm}
\newcommand{\Wh}{{\rm W}}
\newcommand{\Cc}{{\cal C}_c}

\title{Corps de nombres peu ramifi\'es et formes automorphes autoduales}

\author{G. Chenevier et L. Clozel}

\begin{document}

\maketitle

\section*{Introduction}

Soient $S$ un ensemble fini de nombre premiers, $\Q_S \subset \overline{\Q}$
l'extension alg\'ebrique maximale de $\Q$ non ramifi\'ee hors de $S$ (et
l'infini) et $\G_S=\Gal(\Q_S/\Q)$. Un r\'esultat bien connu de Minkowski
affirme que si $S=\emptyset$ alors $\G_S=\{1\}$. En revanche, si $S$ est non
vide, la structure de ces groupes $\G_S$ est tr\`es mal connue, et ce
malgr\'e leur omnipr\'esence en g\'eom\'etrie arithm\'etique. 
Par exemple, un r\'esultat d'Hermite assure que $\G_S$ n'a qu'un nombre fini de sous-groupes ferm\'es 
d'indice donn\'e, mais on ne sait pour aucun $S\neq \emptyset$ si $\G_S$ est topologiquement engendr\'e par un nombre fini d'\'el\'ements ! 
Un autre probl\`eme du folklore consiste \`a d\'eterminer les sous-groupes de d\'ecompositions de $\G_S$. Malgr\'e 
le peu d'indices dont nous disposons pour appr\'ehender cette question, 
il semble commun\'ement esp\'er\'e que ces groupes sont aussi gros que la restriction 
impos\'ee sur la ramification le permet : les r\'esultats de cet article
vont dans cette direction.

Plus pr\'ecis\'ement, supposons $S\neq \emptyset$ et fixons $p \in S$ un nombre premier. La donn\'ee d'un plongement $\Q_S \longrightarrow \Qpb$ definit
un morphisme continu
\begin{equation}\label{decomap} \Gal(\Qpb/\Qp) \longrightarrow \G_S
\end{equation}
dont la classe de conjugaison est ind\'ependante du plongement choisi. Nous 
nous int\'eressons dans cet article \`a la question, soulev\'ee notamment par R. Greenberg, 
de l'injectivit\'e de ce morphisme, o\`u ce qui revient au m\^eme \`a la densit\'e de $\Q_S$ dans $\Qpb$.

Ce probl\`eme a \'et\'e r\'ecemment reconsid\'er\'e dans \cite{Ch}, auquel
nous renvoyons le lecteur pour une discussion plus compl\`ete. Par exemple, il est d\'emontr\'e
{\it loc. cit.} que (\ref{decomap}) est injective d\`es que $S$ contient un
nombre premier $\ell \neq p$ tel que $\ell \equiv 3 \bmod 4$ et que $-\ell$ est un carr\'e modulo $p$. 
Apr\`es quelques r\'eductions \'el\'ementaires, la preuve donn\'ee {\it loc. cit.} consiste \`a d\'emontrer l'existence de suffisament de
repr\'esentations automorphes (sur certains groupes unitaires) ayant des composantes locales inertielles partout
prescrites, et de leur appliquer les travaux de Harris-Taylor
\cite{HT} concernant les repr\'esentations galoisiennes associ\'ees (\'etendant des r\'esultats
an\-t\'e\-ri\-eurs de Clozel et
Kottwitz). Les corps de nombres obtenus sont alors ultimement extraits de
l'action galoisienne sur la cohomologie $\ell$-adique de certains quotients arithm\'etiques
\og explicites\fg\,  des espaces sym\'etriques attach\'es aux groupes 
unitaires r\'eels ${\rm U}(n,1)(\R)$, et ce pour tous les $n\geq 1$.

Notre objectif principal dans ce texte est de supprimer ces hypoth\`eses parasites sur
$\ell$, {\it i.e.} de d\'emontrer le r\'esultat suivant (Th\'eor\`eme
\ref{applicationgalois}).
\ps
\ps
{\bf Th\'eor\`eme A:} {\it Si $|S|\geq 2$, alors (\ref{decomap}) est injective.}
\ps
\ps
En particulier, {\it pour tout entier $m\geq 1$, il existe un corps de nombres de degr\'e multiple de $m$ et non ramifi\'e hors de $S$}. \ps

Ce r\'esultat avait \'et\'e conjectur\'e dans \cite{Ch}, et 
ramen\'e \`a des propri\'et\'es encore largement conjecturales de certaines
formes modulaires de Siegel. La m\'ethode employ\'ee ici est similaire, \`a ceci pr\`es que nous nous passons des groupes symplectiques et
raisonons directement sur le groupe lin\'eaire $\GL_{2n}$. En contrepartie, comme le verrons, les propri\'et\'es
d'autodualit\'e requises rendent les questions d'existence de repr\'esentations automorphes avec
propri\'et\'es locales partout prescrites nettement plus subtiles et ardues. 

Plus pr\'ecis\'ement, fixons $n\geq 1$ arbitraire, choisissons un $\ell \in S-\{p\}$ et fixons une composante de Bernstein 
$\got{c}_p$ du groupe $\GL_{2n}(\Q_p)$. Nous voulons d\'emontrer l'existence d'une repr\'esentation automorphe cuspidale 
$\Pi$ de $\GL_{2n}(\AAA)$ ayant les propri\'et\'es suivantes\footnote{Dans cet article, $\AAA=\AAA_{\Q}$ d\'esignera les ad\`eles de $\Q$.}: 
\begin{itemize}
\item[(P1)] $\Pi$ est autoduale, $\Pi_\infty$ est alg\'ebrique r\'eguli\`ere,
et $\Pi_\ell$ est essentiellement de carr\'e int\'egrable, (de sorte que $\Pi$ rentre dans le
cadre des travaux de Harris et Taylor)
\item[(P2)] $\Pi$ est non ramifi\'ee hors de $\{\infty,\ell,p\}$ et $\Pi_p$ est dans
la composante $\got{c}_p$ fix\'ee.
\end{itemize}
\ps
Le choix de la repr\'esentation de carr\'e int\'egrable $\Pi_\ell$ n'a pas
d'importance pour notre application, et il nous sera en fait commode
d'imposer que \ps
\begin{itemize}
\item[(P3)] $\Pi_\ell$ est la repr\'esentation de Steinberg. 
\end{itemize} 
\ps 
Bien s\^ur, la composante $\got{c}_p$ ne peut \^etre
quelconque puisque $\Pi_p$ est n\'ecessairement autoduale. Ce n'est en fait pas la seule obstruction, nous y reviendrons un peu plus loin. 
En ce qui concerne notre application au th\'eor\`eme plus haut, il nous suffira 
de consid\'erer les composantes $\got{c}_p(\omega)$ induites, \`a partir du sous-groupe de Levi $\GL_n(\Q_p) \times
\GL_n(\Q_p)$, d'une supercuspidale de la forme $\omega \times
\check{\omega}$, o\`u de plus $\omega$ n'a aucun twist non ramifi\'e isomorphe \`a sa contragr\'ediente
$\check{\omega}$ ; mais nos r\'esultats sont en fait de port\'ee plus vaste.
L'essentiel de nos efforts sera dirig\'e vers la preuve du r\'esultat
suivant (Th\'eor\`eme \ref{mainthm}).\ps\ps
{\bf Th\'eor\`eme B}: Pour tout $\omega$ comme ci-dessus et
$\got{c}_p=\got{c}_p(\omega)$, il existe une
re\-pr\'e\-sen\-ta\-tion automorphe cuspidale $\Pi$ de $\GL_{2n}(\AAA)$
satisfaisant les propri\'et\'es (P1), (P2) et (P3).
\ps\ps

Notre m\'ethode pour construire la repr\'esentation $\Pi$ repose sur la
formule des traces d'Arthur-Selberg pour le groupe $\GL_{2n}$ tordu par
l'automorphisme\footnote{Pour des raisons techniques nous serons en fait
amen\'es \`a consid\'erer des variantes b\'enignes de cet automorphisme, 
nous n\'egligeons cet aspect dans cette introduction.} $$\theta(g)={}^t\!g^{-1}.$$
La formule des traces est une identit\'e de distributions $I_{\rm spec}(\cdot)=I_{\rm
geom}(\cdot)$ dont l'utilisation pour ce type de probl\`emes est bien connue: appliqu\'ee \`a des fonctions tests 
$$f = \otimes_v' f_v$$ ne tra\c{c}ant que dans les donn\'ees spectrales qui
nous int\'eressent (ou presque), il s'agit de montrer la non nullit\'e de son c\^ot\'e
g\'eom\'etrique $I_{\rm geom}(f)$.

Hors de $\infty,\, \ell,\, p$, \, $f_v$
sera pour nous simplement la fonction 
caract\'eristique de $\GL_{2n}(\Z_v)$. En les places $\infty,\, \ell$ et
$p$, les fonctions tests n\'ecessaires seront des pseudocoefficients ou des
fonctions de Bernstein bien choisies dont nous devrons ma\^{i}triser les 
int\'egrales orbitales $\theta$-tordues, et une
partie du travail sera de les d\'efinir et d'\'etablir certaines de leurs
propri\'et\'es. Notons que bien que $\GL_{2n}(\R)$ n'ait pas de s\'erie
discr\`ete d\`es que $n>1$, il admet des s\'eries $\theta$-discr\`etes, et
ce sont des repr\'esentations de ce type qui nous int\'eressent \`a l'infini (en
revanche, les composantes $\got{c}_p(\omega)$ d\'ecrites plus haut 
ne sont pas essentiellement discr\`etes, ni m\^eme parmi les composantes
autoduales). Pour ces fonctions tests, une version simplifi\'ee de la formule des 
traces due \`a Arthur \cite{ITF} s'applique, dont le c\^ot\'e g\'eom\'etrique 
se r\'eduit aux termes port\'es par les \'elements $\theta$-semisimples et $\Q$-elliptiques.

Compte tenu de la rigidit\'e de notre probl\`eme, l'unique libert\'e dont nous disposons est de faire varier 
$f_\infty$ parmi les pseudocoefficients de s\'eries $\theta$-discr\`etes
cohomologiques.
Ces repr\'esentations sont en fait naturellement param\'etr\'ees par le poids
extr\'emal $\lambda$ d'une repr\'esentation irr\'eductible $V_\lambda$ du groupe compact
$\SO_{2n+1}(\R)$. Nous d\'emontrerons alors ultimement que lorsque $\lambda$ tend vers
l'infini en s'\'eloignant des murs, le c\^ot\'e g\'eom\'etrique de la formule des
traces devient asymptotiquement \'equivalent \`a un terme unique que nous
appelons le {\it terme principal}. Ce terme est (\`a un scalaire $>0$ pr\`es) 
l'int\'egrale orbitale tordue $\TO_{\gamma_0}(f)$ de $f$ en un certain \'el\'ement $\theta$-semisimple 
$\Q$-elliptique $$\gamma_0 \in
\GL_{2n}(\Q)$$ dont le centralisateur tordu est le groupe symplectique
$\Sp_{2n}$ sur $\Q$. Ainsi,
\begin{equation} \label{geomasymptotique} I_{\rm spec}(f)=I_{\rm geom}(f)
\underset{ \lambda \rightarrow \infty }{\sim} \TO_{\gamma_0}(f) = C \cdot
\dim(V_\lambda),
\end{equation}
\noindent o\`u $C$ est une
constante explicite non nulle ne d\'ependant que de $f^{\infty}$, et
$f_\infty=f_{\infty,\lambda}$. Compte tenu de notre choix de
$f^{\infty}$, le Th\'eor\`eme B est bien s\^ur cons\'equence de la formule
(\ref{geomasymptotique}).\footnote{En fait, une mani\`ere d'interpr\'eter
(\ref{geomasymptotique}) est de la
voir comme une formule de type Riemann-Roch donnant asymptotiquement la dimension
d'un certain espace de formes automorphes de poids variables et de niveau fix\'e
(et engendrant des $\Pi$ satisfaisant (P1), (P2) et (P3)).}\ps\ps

Pour d\'emontrer ces r\'esultats, nous devons surmonter un
certain nombre de difficult\'es dont le traitement est r\'eparti 
de la mani\`ere suivante. Le premier chapitre contient des
pr\'eliminaires sur les caract\`eres des repr\'esentations des groupes
compacts connexes ainsi qu'une illustration de notre m\'ethode dans un contexte
non tordu. Les chapitres $2$, $3$ et $4$ sont consacr\'es au groupe
$\GL_{2n}$ tordu par l'automorphisme $\theta$ et \`a la preuve du Th\'eor\`eme
$B$. Enfin, dans un dernier et court chapitre $5$, nous montrons comment le Th\'eor\`eme $B$ entra\^ine le Th\'eor\`eme
$A$ ; le lecteur peut commencer par celui-l\`a en guise de motivation.
D\'ecrivons maintenant lin\'eairement et plus en d\'etail le contenu des chapitres
$1$ \`a $4$.
\ps\ps

Une id\'ee sous-jacente \`a la m\'ethode expos\'ee ci-dessus (formule
(\ref{geomasymptotique})) est que la formule des traces se simplifie asymptotiquement \og en faisant tendre le poids vers
l'infini\fg . Cette id\'ee se trouve d\'ej\`a dans un argument de
Serre \cite{SerreJams} dans le contexte du groupe $\GL_2$ non tordu. Dans le premier chapitre, nous l'\'etendons \`a son cadre g\'en\'eral
naturel. En guise d'application, nous
\'etablissons le r\'esultat suivant (Th\'eor\`eme \ref{casseriediscrete}). Soit $G$ un $\Q$-groupe r\'eductif connexe tel que $G(\R)$ admet
des s\'eries discr\`etes, et dont les composantes d\'eploy\'ees sur $\R$ et $\Q$ de son
centre co\"incident. \ps
\ps 
{\bf Th\'eor\`eme C}: {\it Soient $\pi$ est une repr\'esentation supercuspidale de
$G(\Q_p)$ et $K$ un sous-groupe compact ouvert de $G(\AAA_f^{\{\infty
,p\}})$. Il existe une repr\'esentation automorphe cuspidale $\Pi$ de $G$ telle
que $\Pi^{\{\infty,p\}}$ admet des $K$-invariants non nuls, $\Pi_p$ est
isomorphe \`a une torsion non ramifi\'ee de $\pi$, et $\Pi_\infty$ est dans
la s\'erie discr\`ete. }
\ps
\ps
La m\'ethode employ\'ee est similaire \`a celle d\'ecrite plus haut \`a ceci
pr\`es qu'elle utilise la formule des traces dans un cas non tordu, ce qui
introduit un certain nombre de simplifications. Par exemple, 
les propri\'et\'es n\'ecessaires des fonctions tests $f_\infty$
sont d\'ej\`a connues. Nous avons besoin toutefois de d\'emontrer le 
r\'esultat suivant sur les caract\`eres des groupes de Lie
compacts connexes (Proposition \ref{asymptotlambda}). Soient $H$ un tel groupe compact et $T$ un tore maximal de $H$, 
notons $V_\lambda$ la representation
irr\'eductible de $H$ de poids extremal $\lambda \in X^*(T)$. Si $\gamma \in H$ est non central, alors 
$$\trace( \gamma ,\, V_\lambda)/\dim(V_\lambda) \longrightarrow 0$$
lorsque $\lambda$ tend vers l'infini dans une direction convenable de
$X^*(T) \otimes \R$.

Le Th\'eor\`eme $C$ admet diverses
variantes: nous pourrions par exemple demander de prescrire de plus les
composantes de Bernstein $\got{c}_v$ de $\Pi_v$ en un
ensemble fini fix\'e de places $v$ diff\'erentes de $\{p,\infty\}$, ce qui
imposerait cependant en g\'en\'eral quelques restrictions sur ces composantes
$\got{c}_v$.\footnote{Un probl\`eme est que nous ne savons pas en g\'en\'eral
montrer qu'une si repr\'esentation automorphe cuspidale $\Pi$ de $G$ est telle que $\Pi_\infty$
est essentiellement discr\`ete de param\`ete suffisament r\'egulier, alors
$\Pi$ est temp\'er\'ee \`a toutes les places.} Nous n'irons pas dans cette direction,
notamment car nous le ferons plus loin dans le cadre de $\GL_{2n}$ tordu. Le r\'esultat
ci-dessus et ses variantes n'entra\^inent pas le Th\'eor\`eme $B$ car $\GL_{2n}$ 
n'a pas les propri\'et\'es requises pour $n>1$. Cependant, nous pourrions l'appliquer au groupe orthogonal d\'eploy\'e
$G=\SO^*_{2n+1}$, de sorte que l'existence des repr\'esentations que nous
recherchons d\'ecoulerait en fait de la conjecture de transfert de $\SO^*_{2n+1}$ vers
$\GL_{2n}$. Malheureusement, les cas actuellement connus de ce transfert n\'ecessitent tous
notamment une hypoth\`ese de g\'en\'ericit\'e de la repr\'esentation \`a
transf\'erer, et nous ne voyons pas comment assurer que nous construisons de telles
repr\'esentations par la formule des traces. C'est pourquoi nous
raisonnons par la suite directement sur le groupe $\GL_{2n}$
tordu\footnote{Une autre m\'ethode aurait consist\'e \`a construire des $\Pi$ comme plus haut
qui sont g\'en\'eriques en produisant directement des s\'eries de Poincar\'e
(cf. \cite[\S 5]{Shahidi2}). Cependant, il semble alors plus d\'elicat de prescrire
$\Pi$ en toutes les places (plut\^ot que toutes sauf une), et nous n'avons pas
poursuivi cette voie. En contrepartie, nous n'utilisons pas les r\'esultats difficiles de transfert suscit\'es.} par $\theta$.  \ps \ps

Le chapitre $2$ contient le travail n\'ecessaire \`a la place archim\'edienne. 
On y d\'efinit et \'etudie en d\'etail les propri\'et\'es des fonctions $f_\infty$ dont nous avons besoin
pour d\'emontrer le Th\'eor\`eme $B$ et dont
nous avons d\'ej\`a parl\'e plus haut. Si $\lambda$ est un poids extr\'emal
d'une repr\'esentation irr\'eductible $V_\lambda$ de $\SO_{2n+1}(\R)$, il
lui correspond d'apr\`es Langlands une unique repr\'esentation $\theta$-discr\`ete
cohomologique $\pi_\lambda$ de $\GL_{2n}(\R)$ (et r\'eciproquement). D'apr\`es
une version du th\'eor\`eme de Paley-Wiener due \`a Mezo \cite{Mezo}, cette repr\'esentation
admet un pseudocoefficient $f_{\infty,\lambda}$ dont la trace sur
$\GL_{2n}(\R)\theta$ isole $\pi_\lambda$ dans le spectre temp\'er\'e autodual de
$\GL_{2n}(\R)$. Le r\'esultat principal de ce chapitre est le suivant
(Th\'eor\`eme \ref{staborbreel}). \ps \ps

{\bf Th\'eor\`eme D:} Soit $\gamma \in \GL_{2n}(\R)$ un \'el\'ement
$\theta$-semisimple. Si $\gamma$ n'est pas $\theta$-elliptique, alors
$\TO_{\gamma}(f_{\infty,\lambda})=0$. Sinon, pour un choix convenable
de mesure positive invariante sur la classe de $\theta$-conjugaison de $\gamma$, on a
$$\TO_\gamma(f_{\infty,\lambda}) = e(\gamma,\lambda)
\trace(\Norm\gamma,V_\lambda),$$  
o\`u $e(\gamma,\lambda)$ est un signe ne d\'ependant que de $\gamma$ et $\lambda$, et $\Norm\gamma
\in \SO_{2n+1}(\R)$ est la \og norme\fg\,  de $\gamma$. En particulier, ces
int\'egrales orbitales sont stables.
\ps\ps
Pour l'\'etude de l'application norme dans ce contexte, nous renvoyons \`a
un article de Waldspurger \cite{Wald}. Les r\'esultats de ce chapitre sont en fait plus complets.
En utilisant des r\'esultats de Bouaziz \cite{bouaziz}, nous commen\c{c}ons par v\'erifier que 
le caract\`ere tordu de $\pi_\lambda$ sur les \'el\'ements elliptiques 
fortement $\theta$-r\'eguliers de $\GL_{2n}(\R)$ co\"incide avec le caract\`ere de $V_\lambda$ via
l'application norme (pour une normalisation convenable de l'op\'erateur
d'entrelacement, cf. Th\'eor\`eme \ref{stabilitetheta}). En particulier il est stable, ce qui est l'analogue d'un r\'esultat de Waldspurger dans le 
cas $p$-adique \cite{Wald2}. Ceci implique le Th\'eor\`eme $D$ pour les int\'egrales orbitales tordues {\it stables}. Un argument 
simple ram\`ene la forme pr\'ecise du th\'eor\`eme au cas o\`u $\pi_{\lambda}$ est \`a cohomologie non triviale pour le syst\`eme 
de coefficients constant, cas o\`u il a \'et\'e d\'emontr\'e par Labesse\footnote{Labesse nous a assur\'e qu'une r\'edaction ult\'erieure 
pr\'eciserait cette d\'emonstration.} \cite{labesse}.  \ps\ps

	Le chapitre $3$ donne la preuve esquiss\'ee plus haut
de la construction de la repr\'esentation $\Pi$. Dans un paragraphe \S
\ref{preliminaireLef}, nous
d\'efinissons et \'etudions la fonction $f_\ell$ dont nous avons
besoin. C'est un pseudocoefficient tordu de la repr\'esentation de Steinberg
de $\GL_{2n}(\Q_\ell)$ dont il nous faut calculer les int\'egrales orbitales
tordues. Il n'est
pas plus long ici de de mener cette \'etude dans le cadre d'une groupe r\'eductif connexe
$G$ g\'en\'eral, et d'un automorphisme $\theta$ d'ordre fini quelconque, et
c'est le choix que nous adoptons. Nous imitons pour cela une m\'ethode de 
Kottwitz \cite{Ko} (voir aussi \cite[\S 9]{BLS})) consistant \`a r\'ealiser g\'eom\'etriquement la fonction $f_\ell$ comme une fonction
d'Euler-Poincar\'e pour l'automorphisme $\theta$ de $G$, et reposant
ultimement sur les propri\'et\'es de l'immeuble de Bruhat-Tits de
$G(\Q_\ell)$, des
travaux de Serre sur les mesures d'Euler-Poincar\'e, et des r\'esultats de
Casselman et Borel-Wallach sur la cohomologie lisse de $G(\Q_\ell)$. \`A
l'aide des r\'esultats des chapitres $1$ et $2$, nous d\'emontrons
la formule (\ref{geomasymptotique}) conditionnellement au r\'esultat
suivant (Th\'eor\`eme \ref{thmcle}) qui fera l'objet du chapitre $4$, et qui permet de montrer la
non-annluation (cruciale) de la constante $C$. \ps\ps

{\bf Th\'eor\`eme E:} Il existe une fonction lisse \`a support compact
$f_p$ sur $\GL_{2n}(\Q_p)$ dont les traces tordues sont nulles hors de la composante de Bernstein
$\got{c}_p(\omega)$ et telle que $$\TO_{\gamma_0}(f_p) \neq 0.$$\ps\ps

Le r\^ole particulier jou\'e par $\gamma_0$ r\'esulte de ce que c'est
l'unique classe de conjugaison $\Q$-elliptique $\theta$-semisimple de
$\GL_{2n}(\Q)$ dont la norme est centrale ({\it i.e.} triviale) dans $\SO_{2n+1}$.
En particulier, cette classe co\"incide avec sa classe de conjugaison stable. En la place
$\ell$, nous d\'emontrons aussi un r\'esultat analogue au Th\'eor\`eme $E$ pour la fonction
d'Euler-Poincar\'e $f_\ell$ (Proposition \ref{nbleftheta}).\ps\ps

	Le chapitre $4$ est vou\'e \`a la preuve du Th\'eor\`eme $E$. Il s'agit du coeur technique de
cet article. Avant d'en dire plus sur sa d\'emonstration, il convient d'en
discuter les tenants et les aboutissants. Oublions temporairement que
$\got{c}_p$ a la forme n\'ecessaire \`a notre application et supposons que
c'est une composante autoduale quelconque. Comme nous l'avions sous-entendu
plus haut, il y a une obstruction \`a ce que l'on puisse construire $\Pi$
avec $\Pi_p$ appartenant \`a $\got{c}_p$. En effet, il est n\'ecessaire que 
le $L$-param\`etre de $\Pi_p$ soit {\it symplectique}. Cela peut se d\'eduire simplement des
th\'eor\`emes de Harris-Taylor et Taylor-Yoshida. De m\^eme, le formalisme du groupe de Langlands sugg\`ere qu'une obstruction \`a inclure un ensemble
fini quelconque de supercuspidales comme composantes locales d'une
repr\'esentation automorphe cuspidale autoduale de $\GL_{2n}(\AAA)$ est
qu'elles soient soit toutes symplectiques, soit toutes orthogonales (ceci
avait d\'ej\`a \'et\'e observ\'e par Prasad et Ramakrishnan dans \cite[\S 3]{PR}, auquel les r\'esultats de cet 
article apportent certaines r\'eponses). 
D'une mani\`ere ou d'une autre, ce type d'hypoth\`ese devra donc appara\^itre
dans notre construction de $\Pi$. Du point de vue de notre m\'ethode, cette obstruction
se traduit exactement par la nullit\'e ou non des int\'egrales orbitales
tordues en l'\'el\'ement $\gamma_0$ des fonctions de Bernstein tordues
$f_p$ de $\got{c}_p$. Dans le cas o\`u $\got{c}_p$ est la composante d'une supercuspidale
$\pi$ autoduale, ceci est en parfait accord avec un r\'esultat de Shahidi
\cite[Prop. 5.1]{Shahidi2}: il montre que cette int\'egrale orbitale est
non nulle pour un coefficient de $\pi$ bien choisi si et seulement si le
$L$-param\`etre de $\Pi$ est symplectique\footnote{Ainsi formul\'e, et ainsi
qu'il l'est expliqu\'e {\it loc. cit.}, le
r\'esultat de Shahidi est conditionel au r\'esultat suivant d\'emontr\'e
ult\'erieurement par Henniart \cite[Thm. 1.3]{He}: $L(\pi,\Lambda^2,s)=L(\Lambda^2 {\rm rec}(\pi),s)$
o\`u ${\rm rec}(\pi): \W_{\Q_p} \longrightarrow \GL_{2n}(\C)$ est le $L$-param\`etre de $\pi$.}. Ainsi, nous avons
montr\'e que l'analogue du Th\'eor\`eme $B$ vaut si l'on prend pour
$\got{c}_p$ une telle composante avec $\pi$ symplectique. Le cas des composantes $\got{c}_p(\omega)$ n\'ecessaires \`a notre
application est \'etudi\'e en d\'etail dans ce chapitre et semble nouveau. Notons que les $L$-param\`etres des
repr\'esentations autoduales de $\got{c}_p(\omega)$ \'etant \`a la fois
symplectiques et orthogonaux, nous nous attendons en fait \`a ce qu'il n'y ait pas d'obstruction
dans ce cas. 

Bien que le th\'eor\`eme $E$ soit de nature locale, notre d\'emonstration utilise des arguments globaux. 
De fa\c{c}on naturelle d'apr\`es la th\'eorie d'Arthur \cite{arthurlivre}, l'alternative symplectique/orthogonale 
pour les repr\'esentations est \'etroitement li\'ee \`a la stabilisation de la formule des traces, la partie \og stable\fg\,  provenant de 
$\SO(2n+1)$ \'etant donn\'ee par les rep\'esentations symplectiques. Pour certaines fonctions $f_{p}$ particuli\`eres 
d\'etermin\'ees par $\got{c}_{p}(\omega)$ (les \og pseudo-coefficients
positifs\fg\,), nous voulons d\'emontrer la non-annulation de 
$\TO_{\gamma_{0}}$. Gr\^ace \`a l'\'etude du caract\`ere tordu des repr\'esentations de $\got{c}_p(\omega)$, 
on v\'erifie que les int\'egrales orbitales stables de $f_{p}$ ne sont pas identiquement nulles. Une version simplifi\'ee et stabilis\'ee de la formule 
des traces d'Arthur produit des $\Pi$ v\'erifiant les conditions pr\'ec\'edentes mais pr\'esentant peut-\^etre de la ramification parasite. Un argument 
nouveau de positivit\'e \S\ref{argumentpositivite}, utilisant un lemme simple sur les caract\`eres des groupes 
compacts (Prop.~\ref{ancien49}), nous permet alors de montrer la 
non-nullit\'e du terme principal de $I_{\rm geom}(f)$ et donc de $\TO_{\gamma_{0}}(f_{p})$. Il est pour ceci crucial de normaliser partout l'op\'erateur 
d'entrelacement, associ\'e \`a l'automorphisme $\theta$, de fa\c{c}on \`a fixer le vecteur de Whittaker. Les sorites n\'ecessaires sont regroup\'es au paragraphe \S\ref{waldwitt}.

Cet argument de positivit\'e semble nouveau et susceptible d'autres
applications ; nous en esquissons quelques unes
dans un dernier paragraphe \S\ref{flocauxcusp} sur les propri\'et\'es
locales-globales des repr\'esentations automorphes cuspidales autoduales de $\GL(2n)$. Nous montrons
tout d'abord que les repr\'esentations automorphes 
\'etudi\'ees par Clozel et Harris-Taylor sont toujours symplectiques
(Th\'eor\`eme~\ref{thmsymplectique}): \ps\ps
{\bf Th\'eor\`eme F:} {\it Soient $F$ un corps totalement r\'eel et $\pi$ une
rep\'esentation automorphe cuspidale de $\GL_{2n}(\AAA_F)$. On suppose que  
$\pi$ est autoduale, essentiellement de carr\'e int\'egrable en au moins une
place finie, et cohomologique \`a toutes les places archim\'ediennes.
Alors pour toute place $v$ de $F$, le $L$-param\`etre de $\pi_v$ pr\'eserve
une forme bilin\'eaire symplectique non d\'eg\'en\'er\'ee.}\ps\ps
En particulier, les repr\'esentations galoisiennes $\ell$-adiques associ\'ees sont aussi
symplectiques. En fait, nous donnons des caract\'erisations des composantes
locales essentiellement discr\`etes des repr\'esentations automorphes ci-dessus
(Th\'eor\`eme \ref{propsymplstein}), \'eclairant notamment certains
aspects du spectre temp\'er\'e de $\GL(2n)$ tordu, sur un corps
$p$-adique. Par exemple, nous obtenons le r\'esultat
suivant, pr\'ecisant ceux de Shahidi discut\'es plus haut. Ici, $K$ est une
extension finie de $\Q_p$. \ps\ps
{\bf Th\'eor\`eme G:} {\it Soit $\pi$ une repr\'esentation supercuspidale autoduale
de $\GL_{2n}(K)$ dont le $L$-param\`etre est orthogonal. Alors les
int\'egrales orbitales stables des pseudocoefficients tordus de $\pi$ sont
toutes nulles.}\ps\ps

Nous terminons en \'enon\c{c}ant une conjecture sur 
la distribution de Plancherel sur le spectre temp\'er\'e autodual de
$\GL_{2n}(K)$ : elle est concentr\'ee sur la vari\'et\'e des repr\'esentations
symplectiques et c'est une mesure sur cette derni\`ere (Conjecture
\ref{conjectureplancherel}). \ps\ps

Pour finir, notons qu'il est sans doute possible de raffiner notre m\'ethode pour d\'emontrer des versions plus fortes du Th\'eor\`eme $B$ (par exemple : remplacer la composante $\got{c}_p(\omega)$ 
par une composante autoduale \og symplectique\fg\,  quelconque,  demander que $\Pi_\ell$ soit non ramifi\'ee, etc...). 
Cependant, m\^eme en admettant ces r\'esultats, ainsi que des g\'en\'eralisations convenables des travaux de Harris-Taylor, 
nous ne voyons pas comment am\'eliorer le Th\'eor\`eme $A$ ({\it i.e.}\textit{} autoriser $S=\{p\}$, cf. 
\cite[\S 4.2]{Ch}).\ps\ps

Nos d\'emonstrations reposent \'evidemment sur l'aride formule des traces tordue d'Arthur. De plus, comme on l'a dit, la forme pr\'ecise des th\'eor\`emes est 
\'etroitement li\'ee aux travaux annonc\'es par lui sur la fonctorialit\'e entre groupes classiques et $\GL(n)$. Certains de nos r\'esultats, 
sans aucun doute, feront partie de l'expos\'e final de sa th\'eorie.

\bigskip
\bigskip

\section{Existence de repr\'esentations en niveau minimal : cas o\`u la caract\'eristique d'Euler-Poincar\'e est non nulle}

\subsection{\'Enonc\'e du r\'esultat}\label{enoncecassimple} Dans ce chapitre, et \`a titre de galop d'essai, 
nous d\'emontrons le r\'esultat naturel d'existence de repr\'esentations automorphes pour 
un groupe r\'eductif $G$ sur $\mathbb{Q}$, v\'erifiant des conditions prescrites en deux places $\{\infty,p\}$, 
lorsque le groupe adjoint a une mesure d'Euler--Poincar\'e non nulle au sens de Serre~\cite{Serre}.
Soit donc $G$ un groupe r\'eductif connexe d\'efini
sur $\Q$. Notons $Z$ le centre de $G$, $S$ le 
sous--tore maximal de $Z$ d\'eploy\'e sur $\Q$ (``composante d\'eploy\'ee de $Z$'').

\begin{hypothese}\label{hyoptore} Les composantes d\'eploy\'ees de $Z$ sur $\Q$ et $\R$ co\"incident. 
\end{hypothese}

On a ainsi une suite exacte
\begin{equation}\label{eqtore}
1 \lgr S \lgr Z \lgr C \lgr 1 
\end{equation}
de $\Q$--groupes diagonalisables, la composante neutre de $C$ \'etant aniso\-tro\-pe sur $\Q$ et $\R$. On note $$A=S(\R)^+$$ 
(on d\'esigne par ${}^+$ les composantes neutres topologiques). 

On suppose enfin que la mesure d'Euler--Poincar\'e sur les quotients arithm\'etiques de $G(\R)/Z(\R)$ -- 
ou, ce qui revient au m\^eme, de $G(\R)/S(\R)$ -- est non nulle ; ceci revient \`a dire que $G(\R)/S(\R)$ a 
une s\'erie discr\`ete, ou qu'il admet une forme int\'erieure compacte.

Fixons un nombre premier $p$. Soit $D$ le plus grand tore d\'eploy\'e quotient de $G$ sur $\Q_p$ 
(noter que $D$ d\'epend de $p$). Soit $\pi_p$ une repr\'esentation supercuspidale de $G(\Q_p)$. 
L'orbite inertielle de $\pi_p$ est l'ensemble des repr\'esentations $\{\pi_p\otimes \chi\}$ o\`u $\chi$ 
parcourt les caract\`eres non ramifi\'es de~$D(\Q_p)$.

On fixe une mesure (finie) $G(\AAA)$-invariante sur $A\,  G(\Q) \backslash G(\AAA)$ et on consid\`ere l'espace 
\begin{equation*}
\mathcal{A}_G = L^2(A\,  G(\Q) \backslash G(\AAA))
\end{equation*}
muni de la repr\'esentation naturelle de $G(\AAA)$.  Une {\it repr\'esentation cuspidale} de $G(\AAA)$ 
sera, par d\'efinition, une repr\'esentation irr\'eductible apparaissant dans le sous-espace des fonctions cuspidales 
(au sens usuel) de $\mathcal{A}_G$. Si $K\subset G(\AAA_f^p)$ est un sous--groupe compact ouvert, elle est de niveau $K$ si son espace des $K$--invariants est non~nul.

\begin{theoreme}\label{casseriediscrete}
Soit $\pi_p$ une repr\'esentation cuspidale de $G(\Q_p)$, et $K\subset G(\AAA_f^p)$. Il existe une repr\'esentation cuspidale $\Pi=\otimes_v \Pi_v$ de $G(\AAA)$, de niveau $K$, telle que \begin{itemize}
\item[(i)] $\Pi_\infty$ est une repr\'esentation de la s\'erie discr\`ete de $G(\R)/A$,
\item[(ii)] $\Pi_p$ appartient \`a l'orbite inertielle de~$\pi_p$.
\end{itemize}
\end{theoreme}

\begin{remarque}\label{remcasseriediscrete}{\rm
La d\'emonstration se simplifie quand $G$ est semisimple. Nous n'avons pas voulu faire cette \hypo\,  car les groupes 
apparaissant dans les applications naturelles des formes automorphes \`a l'arithm\'etique -- cf.~\cite{HT}~-- sont 
rarement semisimples. Ces questions de passage d'un groupe \`a un groupe isog\`ene (au centre pr\`es) rec\`elent 
des ph\'enom\`enes non~triviaux.}
\end{remarque}

\subsection{Fonctions locales : cas r\'eel}\label{fonccasreel}

Nous d\'ecrivons des fonctions particuli\`eres sur $G(\R)$ et $G(\Q_p)$ adapt\'ees \`a notre probl\`eme.

Consid\'erons d'abord la place archim\'edienne. Puisque $G(\R)/A$ a une s\'erie discr\`ete, 
$G$ a une forme int\'erieure r\'eelle $G^*$ anisotrope modulo le centre. On dira qu'une repr\'e\-sen\-ta\-tion 
de $G(\R)$ ou $G^*(\R)$ est dans la s\'erie discr\`ete si elle est unitaire et de carr\'e int\'egrable modulo le centre. 
D'apr\`es Langlands et Shelstad \cite{Shelstadlemme}, il y a une bijection entre repr\'esentations unitaires irr\'eductibles de $G^*(\R)$ 
(dont on notera $\widehat{G^*(\R)}$ l'ensemble) et $L$--paquets de s\'eries discr\`etes de $G(\R)$. 
Pour $\delta \in \widehat{G^*(\R)}$, soit $\Pi(\delta )$ le $L$--paquet associ\'e. 
Le caract\`ere central $\omega$ de toute repr\'esentation $\pi\in\Pi(\delta )$ co\"incide avec celui de $\delta$. 
Nous nous int\'eressons aux repr\'esentations telles que $\omega|_A=1$.

Notons $\overline G$ le groupe $G(\R)/A$. D'apr\`es Clozel--Delorme \cite{CloDel} et Labesse \cite{Labch1}, il existe, 
pour tout $\pi\in\Pi(\delta )$, une fonction $f_\pi\in C_c^\infty(\overline G)$ -- d'ailleurs $K_\infty$--finie pour un 
sous--groupe compact maximal $K_\infty$ de $G(\R)$ -- telle que
\begin{equation}\label{pseudocoeffinfty}
\langle\trace \,\pi,f_\pi\rangle \, = \, 1\,,
\end{equation}
la trace de $f_\pi$ dans toute autre repr\'esentation temp\'er\'ee irr\'eductible de $\overline G$ \'etant nulle. Il en r\'esulte que
$$
f_\pi(zg) = \omega(z)^{-1} f_\pi(g)
$$
si $g\in\overline G$ et $z$ appartient au centre $\overline Z$ de $\overline G$.
La formule (\ref{pseudocoeffinfty}) suppose choisie une mesure de Haar $d\overline g$ sur $\overline G$.

Soit $h$ une fonction $C^\infty$ \`a support compact sur $A$ ; $G(\R)$ est isomorphe \`a $A\times \overline G$ et 
les fonctions $h\otimes f_\pi$ sont donc des fonctions sur $G(\R)$; si $\overline f_\delta ^*$ est un coefficient 
de $\delta$  (avec $\trace \,\delta (\overline f_\delta ^*)\,\neq\, 0)$, $$h \otimes \overline f_\delta ^*\,=:\,f_\delta ^*$$ est 
une fonction sur $G^*(\R)$. Les fonctions $f_\delta =\sum\limits_{\pi\in\Pi} h\otimes f_\pi$ et $f_\delta ^*$, sur $G(\R)$ et $G^*(\R)$, 
sont associ\'ees au sens de Shelstad et de l'appendice de \cite{CloDel}. Il en r\'esulte que l'on a les propri\'et\'es suivantes. Soit $\gamma\in G(\R)$ un \'el\'ement semisimple, $I\subset G$ son centralisateur ; pour $f\in \Cc^\infty (G(\R))$ consid\'erons
$$
\ON_\gamma(f) = \int_{I(\R)\backslash G(\R)} f(g^{-1} \gamma g) \frac{dg}{di}
 $$
 o\`u $dg$, $di$ sont des mesures de Haar. Alors :
 
 \begin{lemme}\label{lemmepseudocorb}\begin{itemize}\item[(i)] Si $\gamma$ n'est pas $\R$--elliptique, $O_\gamma(f_\delta )=0$.
\item[(ii)] Si $\gamma$ est $\R$--elliptique, associ\'e \`a un \'el\'ement $\gamma^*$ de $G^*(\R)$,
$$
O_\gamma (f_\delta)  = e(\gamma) O_{\gamma^*} (f_\delta ^*)
$$
pour une normalisation convenable des mesures sur les centralisateurs $I(\R)$ et $I^*(\R)$, et o\`u $e(\gamma)=\pm 1$ ne d\'epend que de $\gamma$.
\item[(iii)] En particulier, pour $\gamma$ $\R$-elliptique,
$$
O_\gamma(f_\delta)  = e(\gamma) h(\alpha) \theta_\delta ({\overline\gamma^{*}}^{-1})\,,
$$
o\`u $\gamma=\alpha\overline \gamma$ selon l'isomorphisme choisi entre $G(\R)$ et $A\times \overline G$, et o\`u $\overline\gamma^{*}$ est l'image de $\gamma^{*}$ dans $G^{*}(\R)/A$, 
$\theta_\delta$ \'etant le caract\`ere de la repr\'esentation~$\delta $.
\end{itemize}
\end{lemme}

\begin{pf} La partie (i) est bien connue et r\'esulte des propri\'et\'es des fonctions $f_\pi$. La partie (ii) est due \`a Shelstad 
\cite{Shelstadlemme}. 
(Nous ne d\'ecrirons pas les normalisations des mesures, pour nous inessentielles). Le signe $e(\gamma)$ est d\'ecrit par Shelstad. 
Si $\gamma$ est central, $e(\gamma)=1$. Enfin (iii) r\'esulte simplement des relations d'orthogonalit\'e de Schur sur~$\overline G$.
\end{pf}

Noter que quand $\delta $ varie, le support des fonctions $f_\pi$ (et donc $f_\delta$) peut \^etre choisi contenu dans un compact {\it fixe}. 
Ceci r\'esulte de~\cite{CloDel}, ou d'ailleurs de l'argument de Labesse~\cite{Labch1}. Nous devrons enfin contr\^oler la variation 
avec $\delta $ de $\theta_\delta (\overline\gamma)$. Noter que l'on peut \'evidemment faire l'\hypo   suivante :

\begin{hypothese}\label{hypothtore} 
$S$ ne contient aucun sous--tore $S'\not=1$ qui soit facteur direct dans~$G$.
\end{hypothese}

Il en fait \'equivalent de demander qu'un sous-tore (automatiquement
rationnel) $S' \subset S$ soit facteur direct
dans $G$ sur $\R$ et sur $\Q$. En effet, le cocentre $G/G_{\rm der}$ \'etant
$\Q$-isog\`ene \`a $Z$, ses composantes d\'eploy\'ees sur $\R$ et sur $\Q$
co\"incident, de sorte que tout caract\`ere r\'eel $G \rightarrow \GG_m$ est
d\'efini sur $\Q$. Enfin, comme le centre et le cocentre de deux formes
int\'erieures sont canoniquement isomorphes, l'hypoth\`ese \ref{hypothtore}
est encore \'equivalente \`a demander que $G^*$ n'ait aucun tore central
d\'eploy\'e qui soit facteur direct sur $\R$.

\begin{lemme}\label{lemmeconnexe}
Sous l'\hypo \ref{hypothtore}, $G^*(\R)$ est connexe, et $G^*(\R)/A$ est donc compact connexe.
\end{lemme}

\begin{pf} Il existe un
sous-groupe alg\'ebrique $H \subset G^*$ d\'efini sur $\R$ tel que $H(\R)$
soit un sous-groupe compact connexe maximal de $G^*(\R)$. Les alg\`ebres de
Lie de $S$ et de $H$ sont alors en somme directe dans ${\rm Lie}(G^*)$, de
sorte que le $\R$-morphisme naturel $S \times H \longrightarrow G^*$
est surjectif (sur les $\C$-points), de noyau $\mu := S \cap H$ fini, {\it
i.e.} 
\begin{equation}\label{devissagelemme} 
1 \longrightarrow \mu \longrightarrow (S \times H) \longrightarrow G^*
\longrightarrow 1.
\end{equation}

En fait, la suite ci-dessus reste exacte apr\`es passage aux $\R$-points,
{\it i.e.} $G^*(\R)=H(\R)\cdot S(\R)$. En effet, $G^*/S$ 
est r\'eductif connexe et par construction $(G^*/S)(\R)=G^*(\R)/S(\R)$ est
compact:  ce dernier est donc aussi connexe. Cela conclut car $H(\R)$ est connexe et ${\rm Lie}(H(\R))
\rightarrow {\rm Lie}((G^*/S)(\R))$ est un isomorphisme. 

Le tore $S$ \'etant d\'eploy\'e, $S(\R)= A \times \{\pm 1\}^{{\rm dim}(S)}$
et $\mu \subset \{\pm 1\}^{{\rm dim}(S)}$. L'exactitude de (\ref{devissagelemme})
sur les $\R$-points assure que 
$$ G^*(\R) \simeq A \times \{\pm 1\}^{{\rm dim}(S)}/\mu \times H(\R).$$
Mais si $\mu \subsetneq \{\pm 1\}^{{\rm dim}(S)}$, alors $S$ admet un sous-tore strict - nec\'essairement facteur direct sur $\R$
- contenant $\mu$, de sorte que la suite exacte (\ref{devissagelemme})
contredit l'hypoth\`ese \ref{hypothtore}, ce qui conclut.
\end{pf}

Soit alors $\overline T$ un tore maximal de $\overline G^* = G^*(\R)/A$. On peut param\'etrer 
les re\-pr\'e\-sen\-ta\-tions $\delta $ par 
leur plus haut poids $\lambda\in X = X^*(\overline T)$. Si $\overline \gamma$ appartient au 
centre $\overline Z$ de $\overline G^*$, on a
$$
\theta_\delta (\overline\gamma)=\deg(\delta )\omega(\overline \gamma)\,.
$$
Notons $\delta (\lambda)$ la repr\'esentation associ\'ee \`a $\lambda\in X$. Alors $\deg \delta (\lambda)$ est donn\'e par le polyn\^ome de Weyl 
$$
P(\lambda) = {\prod_{\alpha}\frac{\langle \alpha,\lambda+\rho\rangle}{\langle\alpha,\rho\rangle}}
$$
les produits portant sur un ensemble de racines positives dont $\rho$ est la demi--somme. \ps\ps

La proposition qui suit sera appliqu\'ee \`a $\overline G^*$ mais est vraie pour tout groupe de Lie compact connexe. 
Jusqu'\`a la fin du \S\ref{fonccasreel} $G$ d\'esignera un tel groupe, $T$ un tore maximal de $G$ et $X$ 
le groupe des caract\`eres de $T$. On fixe un ensemble de racines positives pour $(G,T)$. 
Pour $\lambda\in X$ dominant, soit $\theta_\lambda$ le caract\`ere de la repr\'esentation de plus 
haut poids $\lambda$ de $G$.

\begin{prop}\label{asymptotlambda} Soit $\gamma\in G$. Pour $\lambda$ dominant
$$
\theta_\lambda(\gamma) = \sum_i E_i(\gamma,\lambda) \ P_i(\lambda)\,.
$$

La somme est finie ; les $E_i(\gamma,\lambda)$ sont des polyn\^omes (dont les degr\'es d\'ependent de $\lambda$) 
en les $\gamma^\chi$, \`a $\chi$ parcourt une base de $X$ ; $E_i(\gamma,\lambda)$ est {\rm 
uniform\'ement born\'e} quand $\lambda$ varie. De plus, $P_i(\lambda)$ est un polyn\^ome de degr\'e strictement inf\'erieur \`a celui de $P(\lambda)$ si $\gamma$ n'est pas central.
\end{prop}

\begin{pf} Si $\gamma$ est r\'egulier ou central, la proposition n'est autre que la formule du caract\`ere ou du degr\'e de Weyl. Si le centralisateur de $\gamma$ est un sous--groupe de Levi, elle se d\'eduit de la formule de Kostant. En g\'en\'eral, nous imitons l'une des d\'emonstrations de celle--ci.

Soit $G_\textrm{der}$ le groupe d\'eriv\'e de $G$. On v\'erifie aussit\^ot qu'il existe un rev\^etement connexe fini $\widetilde G$ 
de $G$ dont le groupe d\'eriv\'e est simplement connexe, et que la proposition pour $G$
r\'esulte\footnote{Remarquer notamment que l'image inverse du centre de $G$ est le centre de $\widetilde
G$.} du r\'esultat pour $\widetilde G$. On suppose donc $G_\textrm{der}$ simplement connexe ; le centralisateur de tout \'el\'ement $\gamma$ est alors connexe.

Fixons $\gamma$ un tel \'el\'ement, on peut supposer quitte \`a le conjuguer que $\gamma \in T$. Notons $M$ le centralisateur de $\gamma$ dans $G$, 
soit $R(G,T)$ et $R(M,T)$ les ensembles de racines associ\'es, $R_+ \subset R$ les racines positives, 
$\Delta \subset R_+$ les bases, pour $G$ et $M$, $\rho$, $\rho_M$ les demi--sommes de racines associ\'ees. On suppose 
que $R_+(M,T)=R_+(G,T) \cap R(M,T)$ ; quitte \`a prendre un nouveau rev\^{e}tement on peut supposer que 
$\rho$ et $\rho_M$ appartiennent \`a $X$. Soit $W$ le groupe de Weyl de $(G,T)$, $W_M$ celui de $(M,T)$. Soit
$$
W^M = \{ w \in W \,\,:\,\, w^{-1} \alpha\in R_+ (G,T)\, \,  \forall \alpha\in\Delta_M\}\,.
$$
Comme dans le cas parabolique :

\begin{lemme}\label{lemmepara}
Si $w\in W$, $w$ s'\'ecrit de mani\`ere unique $w=w_s w_u$ avec $w_s \in W_M$ et $w_u\in W^M$.
\end{lemme}

\begin{pf} Soient en effet $C_G^+$, $C_M^+$ les int\'erieurs  des chambres de Weyl positives, dans $X\otimes\R$, relatives \`a $G$ et $M$. Fixons $H\in C_G^+$. Alors $w\in W$ appartient \`a $W^M$ si et seulement~si
$$
\langle\alpha, w(H)\rangle  >0 \qquad (\alpha\in \Delta_M)\,,
$$
c'est--\`a--dire  si  $w(H) \in C_M^+$.

Soit $w \in W$ : il existe alors $w_s\in W_M$ tel que $w_s^{-1} w (H) \in C_M^+$ : 
donc $w_s^{-1}w \in W^M$, d'o\`u la d\'ecomposition cherch\'ee. Si $w=w_sw_u=w_s'w_u'$,  $w_u(H)$ et $w'_u(H)$ 
sont dans $C_M^+$ et $w_s C_M^+$ rencontre $w_s' C_M^+$, donc~$w_s=w_s'$. 
\end{pf}

\begin{lemme}\label{lemmedom}
Si $\lambda\in X$ est dominant pour $G$ et $w_u\in W^M$, $w_u(\lambda+\rho)-\rho_M$ est dominant pour~$M$.
\end{lemme}

\begin{pf}Il faut en effet v\'erifier que
$$
2 \frac{\langle w_u(\lambda+\rho),\alpha\rangle}{\langle \alpha,\alpha \rangle}\ - 
2 \frac{\langle\rho_M,\alpha\rangle}{\langle \alpha,\alpha \rangle}\ \geq 0\qquad (\alpha\in \Delta_M)\,.
$$

Le premier terme est entier, positif car $w_u^{-1}\alpha\in R^+(G,T)$, strictement positif car $\lambda+\rho$ est r\'egulier. 
Le second est \'egal~\`a~$1$.\end{pf}

D\'emontrons alors la proposition. Ecrivons d'abord, formellement (= dans le groupe de Grothendieck  de~$T$) :
$$
\theta_\lambda = {\frac{N}{D}}
$$
o\`u 
\[
\begin{array}{lll}
  N& =\displaystyle\sum_{w\in W} \varepsilon(w) e^{w(\lambda+\rho)}\,,    \\
  D&=\displaystyle\sum_{w\in W} \varepsilon(w) e^{w\rho} \,,     \\
\end{array}
\]
$\varepsilon$ \'etant le signe sur $W$. Alors
\[
\begin{array}{lll}
  &N   &= \displaystyle\sum_{w_u \in W^{M}} \varepsilon(w_u) \displaystyle\sum_{w_s\in W_{M}} \varepsilon(w_s) e^{w_s w_u(\lambda+\rho)}   \\
  &   &= \displaystyle\sum_{w_u \in W^{M}} \varepsilon(w_u) \displaystyle\sum_{w_s \in W_{M}} \varepsilon(w_s) e^{w_s(\lambda_u+\rho_M)}   \\
\end{array}
\]
o\`u l'on a \'ecrit $w_u(\lambda+\rho)=\lambda_u+\rho_M$ d'apr\`es le Lemme~\ref{lemmedom}, avec $\lambda_u$ dominant pour $M$. 
Chacun des termes de la somme ci-dessus index\'e par $w_u$ est donc un num\'erateur de Weyl pour $M$. Par~ailleurs, s'il on pose 
$S(M)=R^+(G,T)\backslash R^+(M,T)$, 
\[
\begin{array}{lll}
   D  &= e^\rho \underset{\alpha \in R^{+}(G,T)}{\prod} (1-e^{-\alpha})                         \\
        &= D_{M} \, e^{\rho-\rho_M} \, \underset{\alpha \in  S(M)}{\prod}(1-e^{-\alpha})\,.
\end{array}
\]
Ainsi
$$
\frac{N}{D} = e^{\rho_M-\rho} \underset{\alpha \in S(M)}{\prod} (1-e^{-\alpha})^{-1} \left( \sum_{w_u} \varepsilon(w_u)\,\theta_{M_,\lambda_u}\right)
$$
o\`u $\theta_{M,\lambda_u}$ est le caract\`ere de $M$ associ\'e \`a $\lambda_u$. Cette expression donne la valeur de 
$\theta_\lambda$ pour un \'el\'ement  $G$--r\'egulier $t\in T$. Si $t\lgr \gamma$, les termes $(1-t^{-\alpha})$ pour 
$\alpha\notin R(M,T)$ restent non nuls ; $\theta_{M,\lambda_u}$ a pour limite $\omega_{\lambda_u}(\gamma)$ 
$P_M(\lambda_u)$, o\`u $\omega_{\lambda_u}(\gamma)=\gamma^{w_u(\lambda+\rho)-\rho_M}$ et $P_M$ est le polyn\^ome de Weyl pour $M$. 
Ainsi
$$
\theta_\lambda(\gamma)= \frac{\gamma^{\rho_M-\rho}}{\Pi_{\alpha \in S(M)}(1-\gamma^{-\alpha})} \sum_{w_u \in W^{M}} \varepsilon(w_{u}) \gamma^{w_u(\lambda+\rho)-\rho_M} P_M(\lambda_u)
$$
ce qui d\'emontre la proposition. L'assertion sur le degr\'e des $P_{i}$ r\'esulte des formules du degr\'e de Weyl pour $G$ et $M$ 
et de ce que $S(M) \neq \emptyset$ si et seulement si $\gamma$ n'est pas central.
\end{pf}

\begin{corollaire}\label{corlimitecar} Si $\gamma \in G$ n'est pas central, alors $\frac{\theta_{\lambda}(\gamma)}{P(\lambda)}$ 
tends vers $0$ lorsque $\lambda \in X^{*}(T)\otimes \R$ tend vers l'infini en s'\'eloignant des murs des chambres de 
Weyl\,\footnote{Nous entendons par l\`a que pour tout $\alpha \in R^{+}(G,T)$, $\langle \lambda, \alpha \rangle$ tend vers l'infini.}.
\end{corollaire}

\begin{pf} Dans les notations pr\'ec\'edentes, nous avons explicitement
$$\frac{P_{M}(\lambda_{u})}{P(\lambda)} = \frac{\prod_{\alpha \in R^{+}(G,T)} \langle \rho, \alpha \rangle}{\prod_{\alpha \in R^{+}(M,T)}\langle
 \rho_{M},\alpha \rangle} \,\left(  \underset{\alpha \in R^{+}(G,T)\backslash w_{u}^{-1}  R^{+}(M,T)}{\prod} 
 \langle \lambda + \rho, \alpha \rangle \right)^{-1}.$$
\end{pf}

\subsection{Fonctions locales en $p$}${}^{}$\ps\ps

L'ensemble des caract\`eres non ramifi\'es $\chi$ de $D(\Q_p)$ forme un tore complexe ; soit $F_p(\chi)$ une fonction polynomiale sur celui--ci et qui de plus est invariante par le sous-groupe (fini) des $\chi$ tels que 
$\pi_p \otimes \chi \simeq \pi_p$. Par la th\'eorie de Bernstein, on sait alors qu'il existe une fonction $f_p$ sur $G(\Q_p)$ ne tra\c{c}ant de mani\`ere non nulle que dans l'orbite inertielle de $\pi_p$, et telle que
\begin{equation}\label{pseudoex19}
\langle \trace (\pi_p\otimes \chi), \ f_p \rangle = F_p(\chi). 
\end{equation}
La formule de Plancherel montre que $f_p(1)=1$ pour un choix convenable de $F_p$.

Soit par ailleurs $S(\Q_p)\subset G(\Q_p)$ (\S\ref{enoncecassimple}) et soit $\varepsilon$ un \'el\'ement de 
$\{\pm1\}^{\dim S}$, plong\'e dans $S(\Q_p)$. Puisque $\chi(\varepsilon)=1$, la formule de Plancherel implique~que
$$
f_p(\varepsilon)=\omega_p(\varepsilon)^{-1}
$$
o\`u $\omega_p$ est le caract\`ere central de~$\pi_p$.

\subsection{D\'emonstration du Th\'eor\`eme \ref{casseriediscrete}}

Notons $R$ la repr\'esentation de $G(\AAA)$ dans la partie cuspidale de $\mathcal{A}_G$. Soit $f$ la fonction 
sur $G(\AAA)$ donn\'ee par  $$f=f_\infty\otimes f_p\otimes f^{p,\infty}$$ o\`u $f_\infty=f_\delta $ et $f_p$ viennent 
d'\^etre d\'efinies, et $f^{p,\infty}$ est la fonction caract\'eristique de $K$. Puisque $f_p$ est cuspidale, et que les int\'egrales orbitales de 
$f_\infty$ en les \'el\'ements non elliptiques s'annulent, la formule des traces simplifi\'ee de Deligne et Kazhdan (cf. Henniart~\cite{henniartftr}) s'applique :
$$
\trace\,  R(f)\, =\, \sum_\gamma v(\gamma) \, O_\gamma(f_\infty) \,O_\gamma(f_p) \,O_\gamma(f^{p,\infty})
$$
o\`u $\gamma$ parcourt les classes de conjugaison elliptiques de $G(\Q)$ -- en fait $\R$--elliptiques vu 
les propri\'et\'es de $f_\infty$ -- et $v(\gamma)>0$ est un volume. Puisque les supports de nos fonctions sont 
contenus dans un compact fixe, cette somme est finie, uniform\'ement quand $f_\infty$ varie (on pourrait aussi, bien s\^ur, d\'eduire (1.10) des r\'esultats d'Arthur~[ITF]).

Pour $f_\infty=f_\delta $ et $\gamma$ central, $O_\gamma(f_\infty)=f_\infty(\gamma)$ est essentiellement \'egale au polyn\^ome de Weyl en le param\`etre 
$\lambda$ de $\delta $. Consid\'erons une suite de $\lambda\in X^*(\overline T)$ tendant vers l'infini dans une direction r\'eguli\`ere. 
Si $\gamma$ n'est pas central, $O_\gamma(f_\delta )$ est n\'egligeable par rapport \`a $P(\lambda)$ d'apr\`es le
Lemme~\ref{lemmepseudocorb} et
le Corollaire~\ref{corlimitecar}. On en d\'eduit~que
\begin{equation} \label{traceasymp}
P(\lambda)^{-1}\  \trace \,R(f) \,= \,v\,\sum_{\gamma\in Z(\Q)}h(\alpha)\omega_\lambda^{-1}(\overline \gamma) f_p(\gamma) f^{p,\infty}(\gamma) \,+\, \, o (1)\,,
\end{equation}
o\`u l'on a utilis\'e les notations du Lemme~\ref{lemmepseudocorb} : $\gamma=\alpha\overline \gamma$, $\alpha\in A$, $\overline \gamma \in G(\R)/A$. 
La fonction $h \in {\rm C}_c^\infty(A)$ est pour l'instant arbitraire.

Utilisons la suite exacte (\ref{eqtore}), qui reste exacte apr\`es passage
aux points rationnels. Si $F$ est la fonction de $\gamma$ figurant dans la somme
(\ref{traceasymp}), on consid\`ere donc
\begin{equation} \label{traceasymp2}
\sum_{c \in C(\Q)} \sum_{s\in S(\Q)} F(s \gamma),
\end{equation}
o\`u l'on a choisi un repr\'esentant $\gamma\in Z(\Q)$ de $c$. 
\'Ecrivons $s =s ^+\varepsilon$, o\`u $\varepsilon\in \mathcal{E}:=
\{\pm1\}^{{\rm dim}(S)} \subset S(\Q)$ et $s^+\in S(\Q)\cap A$. 
Si~$\gamma=\alpha\overline \gamma$,
$$
f_\infty(s \gamma) = h (s^+\alpha) \omega_\lambda^{-1} (\overline \varepsilon\overline \gamma)\,.
$$

La somme \'etant finie, on peut supposer que les seuls termes pr\'esents v\'erifient 
$s^+\alpha=1$ en prenant le support de $h$ suffisament proche de 1 ; 
quitte \`a changer le repr\'esentant $\gamma$ on peut donc supposer $\alpha=1$ et donc
$s^+=1$ pour tous les termes de
(\ref{traceasymp2}), qui s'\'ecrit alors 
\begin{equation}\label{traceasymp3}
\sum_{c\in C(\Q)}\  \sum_{\varepsilon\in\mathcal{E}} \omega_\lambda^{-1}(\varepsilon\overline \gamma) f_p(\varepsilon\gamma) f^{p,\infty}(\varepsilon\gamma)\,. 
\end{equation}
Si $c=1$, le terme correspondant est
$$\sum_\varepsilon \omega_\lambda^{-1} (\varepsilon) f_p(\varepsilon) f^{p,\infty}(\varepsilon)\,.
$$
Puisque $f_p(\varepsilon)=\omega_p^{-1}(\varepsilon)$, il est \'egal \`a
$$
\sum_\varepsilon f^{p,\infty}(\varepsilon) \, \, \geq  \, \, f^{p,\infty}(1) \,\,>\,\,0
$$
si $\omega_\lambda|_\mathcal{E}= \omega_p^{-1}|_\mathcal{E}$. Fixons $\lambda_0$ 
v\'erifiant cette condition. Soit $\lambda=\lambda_0+\mu$, avec $\mu\in X$ trivial 
sur $\mathcal{E}\subset \overline G^*$. Soit $\overline Z$ le centre de $\overline G^*$, 
donc $\overline Z=Z(\R)/A$. La suite exacte
$$
1 \lgr S(\R) \lgr Z(\R) \lgr C(\R) \lgr 1
$$
implique que $\overline Z/\mathcal{E}=C(\R)$. Pour de tels $\lambda$,
(\ref{traceasymp3}) se r\'e\'ecrit
\begin{equation}\label{traceasymp4}
\sum_{c \in C(\Q)} a(c) \mu^{-1}(c) 
\end{equation}
o\`u l'on a identifi\'e dans l'\'ecriture $\mu(c)$ l'\'el\'ement $c\in C(\Q)$ \`a l'image de $\overline \gamma\in\overline Z$ 
dans $C(\R)$. De plus $a(1)\not= 0$ et $a(c)$ est \`a support finie. 
Puisque, tout en faisant tendre $\lambda$ vers l'infini dans les
directions semisimples de $\overline G^*$, on peut prendre pour
$\mu_{|C(\R)}$ des caract\`eres arbitraires, la th\'eorie de Fourier sur ce tore compact 
implique que (\ref{traceasymp4}) prend des valeurs non~nulles.

D'apr\`es (\ref{traceasymp}), on en d\'eduit que pour $\lambda$ assez grand
et soumis aux conditions sp\'ecifi\'ees ci-dessus,
$$
P(\lambda)^{-1}\ \trace \, R(f) \not=0\,.
$$

Mais 
$$
\trace \, R(f)\, =\, \sum_{\Pi\, \,  {\rm cuspidale}} \trace \, \Pi_\infty(f_\infty)\, \, \trace
\, \Pi_p (f_p) \, \dim (\pi^{p,\infty})^{K_{p,\infty}}\,,
$$
avec $\trace \, \Pi_\infty(f_\infty)\, = \, (\int_A h(a)da)\, \trace\, \overline
\Pi_\infty(f_\delta )$, o\`u $\overline \Pi_\infty$ est la restriction de
$\Pi_\infty$ \`a~$\overline G$.
Si le param\`etre $\lambda$ est assez r\'egulier, toute repr\'esentation de 
$\overline G$ v\'erifiant $\trace\,  \overline \Pi_\infty(f_\delta )\not=0$ 
appartient au $L$--paquet de s\'eries discr\`etes $\Pi(\delta )$ (\cite{V}). 
Ceci termine la d\'emonstration.

\vskip2mm

\begin{remarque}\label{remarthur}{\rm Comme on l'a indiqu\'e, le recours 
\`a la formule des traces simplifi\'ee n'est pas n\'ecessaire. En particulier, on aurait pu consid\'erer une repr\'esentation $\pi_p$ 
appartenant \`a la s\'erie discr\`ete, pour laquelle on peut construire une fonction $f_p$ jouissant de propri\'et\'e analogues. 
Mais il n'est pas vrai alors, en g\'en\'eral, que la trace de $f_p$ n'est non--nulle que pour les repr\'esentations obtenues 
par torsion de $\pi_p$ (prendre $G=GL_2$, $\pi_p =$ repr\'esentation de Steinberg). Si les repr\'esentations $\pi_p\otimes \chi_p$ 
sont les seules repr\'esentations {\it unitaires} ayant cette propri\'et\'e, la d\'emonstration reste valide. }
\end{remarque}

\bigskip
\bigskip

\section{Le cas de $\GL(2n)$ : \\ repr\'esentations et fonctions \`a la place
archim\'edienne}
\setcounter{equation}{0}

Comme on l'a dit dans l'introduction le reste de l'article  est consacr\'e au groupe $\GL(2n)$ sur $\Q$. 
Dans ce chapitre nous voulons imiter dans ce cas les r\'esultats r\'eels du \S\ref{fonccasreel}. 
Puisque $\GL(2n,\R)$, pour $n>1$, n'a pas de s\'erie discr\`ete, les m\'ethodes du Ch.~1 ne s'appliquent pas. 
On doit donc consid\'erer l'automorphisme ext\'erieur (essentiellement : $g\mapsto {}^{t}\!g^{-1}$) de $\GL(2n)$ ; 
tous les objets (repr\'esentations, int\'egrales orbitales, formule des traces\dots) seront les objets tordus pour cet 
automorphisme.

\subsection{Repr\'esentations $\theta$--stables}\label{prelimglnreel} Soit $J$ la matrice
antidiagonale
\[
\left(
\begin{array}{ccc}
  &   &  1 \\
  &\adots   &   \\
 1 &   &    
\end{array}
\right) \in \GL(2n)\,.
\]

On note $\G$ le groupe $\GL(2n)$ ; on le munit de l'automorphisme $g \mapsto \theta(g) = J {}^{t}\!g^{-1}J$. 
Noter que $\theta$ pr\'eserve le sous--groupe de Borel usuel de $\G$ ; 
il ne pr\'eserve pas un \'epinglage (cf.~\cite{KS}) (on revient la--dessus au chapitre 4).

Nous suivons Waldspurger \cite{Wald}, ce qui nous permettra d'utiliser ses r\'esultats sans changement.

Soit $\pi$ une repr\'esentation (admissible) irr\'eductible de $\G(\R)$. Alors $\pi$ est $\theta$--invariante 
$(\pi\cong \pi^\theta$ o\`u $\pi^\theta= \pi\circ \theta)$ si et seulement si $\pi$ est isomorphe \`a sa duale $\widetilde \pi$. Notons $\widehat\G
=\GL(2n,\C)$ le groupe dual. La classification de Langlands associe \`a $\pi$ une repr\'esentation semisimple
$$
r(\pi) : \W_\R \lgr \widehat \G\,.
$$

D'apr\`es des faits bien connus, $\pi\cong \widetilde \pi$ si, et seulement si, $r(\pi)$ est isomorphe 
\`a~$\widetilde{r(\pi)}$.

Soit $p_1,\ldots,p_n\in \frac{1}{2}+\Z$. Si $p\in \frac{1}{2}\Z$ on lui associe 
un caract\`ere $\chi$ de $\C^\times$, not\'e d'apr\`es Langlands 
$z \mapsto z^p(\overline z)^{-p}$. C'est un caract\`ere unitaire. 
Soit $r(\chi)$ la repr\'esentation de $\W_\R$ induite de $\chi$ : 
Langlands lui associe une repr\'esentation de la s\'erie discr\`ete (unitaire) de $\GL(2,\R)$, 
que l'on notera $\delta (p)$. Si $p\in \frac{1}{2}+\Z$, $\chi|_{\R^\times}$ est \'egal au signe. 
La formule bien connue
$$
\det \mathrm{ind}(\chi)=(\chi|_{\R^\times}) \varepsilon_{\C/\R}
$$
montre que le caract\`ere central de $\delta (p)$ est trivial.

Si $r:\W_\R \lgr \GL(2n,\C)$ est la somme des $r(\chi_i)$ pour $\chi_i$ associ\'e \`a
$p_i$ (dans $\frac{1}{2} +\Z$), la repr\'esentation $\pi$ associ\'ee est l'induite
unitaire $$
\mathrm{ind}_P^{\G} (\delta _1 \otimes \cdots\otimes \delta _n) \qquad (\delta _i = \delta (p_i))
$$
o\`u $P$ est le parabolique de type $(2,\ldots,2)$. Elle est temp\'er\'ee. Puisque le caract\`ere central de $\delta (p)$ est trivial, 
celle--ci est autoduale ; il en est de m\^eme~de~$\pi$.

La repr\'esentation $\pi$ est {\it alg\'ebrique} au sens de \cite{AA}, \cite{HT}. 
Supposons de plus les $p_i$ distincts. Alors $\pi$ est cohomologique. Plus pr\'ecis\'ement, soit 
$$
p(\pi) = (p_1, p_2,\ldots, p_n, -p_n,-p_2,\ldots,-p_1)
$$
o\`u l'on a suppos\'e $p_1>p_2\cdots > p_n >0$, et
\[
\begin{array}{lll}
  &m(\pi)   & = \Big(p_1 - \frac{2n-1}{2},\, p_2 - \frac{2n-3}{2},\ldots, p_n, - \frac{1}{2},\, -p_n
+\frac{1}{2},\ldots, -p_1 +\frac{2n-1}{2}\Big)   \\
  &   &= (m_1,m_2,\ldots,m_{2n}) \,.   \\ 
\end{array}
\]
Alors $m(\pi)$ est le plus haut poids d'une repr\'esentation rationnelle $V$ de $\GL(2n)$, et
$$
H^{\bullet }(\mathfrak{g}, K_\infty ; \pi\otimes V)\not= 0
$$
(cf. \cite[Lemme 3.14]{AA} ; noter que $V$ est autoduale).

Pour expliquer les calculs qui suivent, d\'ecrivons la factorisation de $r=r(\pi)$ donn\'ee 
par la fonctorialit\'e d'Arthur. Comme on l'a remarqu\'e, chaque repr\'esentation $r(\chi)$ est de 
d\'eterminant trivial, donc
$$
r(\chi) : \W_\R \lgr \SL(2,\C).
$$

Ainsi $r(\chi)$ laisse invariante une forme symplectique sur $\C^2$ ; on en d\'eduit que dans une base convenable,
$$
r : \W_\R \lgr \Sp(2n,\C) \subset \GL(2n,\C)\,.
$$
Le groupe complexe $\widehat H = \Sp(2n,\C)$ est le groupe dual du groupe sp\'ecial orthogonal
d\'eploy\'e de type $B_n$, que l'on notera
$H=\SO^*(2n+1)$ (sur $\R$ ou $\Q)$. Par cons\'equent $r$ d\'efinit un param\`etre de Langlands (r\'eel) pour $H$. On v\'erifie aussit\^ot que $r$ ne se factorise par aucun sous--groupe de Levi de $\widehat H$ ; d'apr\`es Langlands $r$ d\'efinit donc un $L$--paquet de s\'eries 
discr\`etes de $H(\R)$, ou d'ailleurs de toute forme int\'erieure de celui--ci ; en particulier $r$  d\'efinit une repr\'esentation de dimension finie~$\pi_H$~de
$$
H_c(\R) =\SO(2n+1,\R)\,.
$$

D\'ecrivons celle--ci. Le tore maximal $T_H$ de $H_c$ s'identifie naturellement \`a $(T_a)^n$ o\`u $T_a$ est le tore
r\'eel anisotrope de dimension $1$ ; la demi--somme des racines, pour un ordre convenable, est alors
$$
\rho_H = \Big(\frac{2n-1}{2}\,,\, \frac{2n-3}{2},\ldots,\frac{1}{2}\Big) \in X^* (T_H) \otimes \frac{1}{2}\Z = \Big(\frac{1}{2}\Z\Big)^n\,.
$$

La repr\'esentation $\pi_H$ a pour plus haut poids $m_H=p_H -\rho_H \in X^*(T_H)$, o\`u
$$
p_H = (p_1,\ldots,p_n) \, ;
$$
$p_H$ est le param\`etre du caract\`ere infinit\'esimal de~$\pi_H$.

\subsection{Application norme}\label{parnorme}

Dans ce paragraphe nous d\'ecrivons rapidement suivant
Waldspurger~\cite[\S III.2]{Wald} -- auquel nous renvoyons le lecteur pour 
les d\'etails -- l'application (ou plut\^ot la correspondance), g\'en\'eralement appel\'ee
{\it norme}, entre classes de conjugaison tordue dans $\GL(2n,\R)$ et classes de conjugaison dans
$\SO^*(2n+1,\R)$.\ps

Notons $\G^+$ le produit semi--direct de $\{1,\theta\} \simeq \Z/2\Z$ par
$G$, $\theta$ op\'erant par
 $x\mapsto J {}^{t}\!x^{-1}J$. On a $\G^+ = \G \coprod \theta \G = 
\G \coprod \tilde\G$. Si $g$, $h\in \G$ on dit que $g$ et $h$ sont $\theta$--conjugu\'es si $g=x^{-1}h x^\theta$ pour un $x\in \G$ ; il revient au m\^eme de dire que $\theta g$, $\theta h$ sont conjugu\'es par $\G\subset \G^+$. 
Ces notions ont un sens pour les points \`a valeurs dans un corps quelconque.

$\G^+$ \'etant r\'eductif (non connexe) il y a une notion naturelle d'\'el\'ement semisimple dans $\G^+$ ;  
$\widetilde g = \theta g$ $(g\in\G)$ est semisimple si et seulement si $\widetilde g^2 =(g^\theta)g$ est 
semisimple. On dira que $g$ est $\theta$--semisimple. Si $g$ est semisimple, son centralisateur $Z_{\G}(g)$ 
dans $\G$ est r\'eductif.
On dit que $g$ est fortement r\'egulier si c'est un tore\footnote{Cette 
d\'efinition diff\`ere en apparence de celle de Waldspurger \cite[\S
I.1]{Wald}, mais sa classification 
des centralisateurs $Z_{\G}(g)$ (\cite[{\rm fin du} \S I.3]{Wald}) montre qu'elle est \'equivalente pour~$\GL(2n)$.}.  
Si $g$, $h$ sont fortement r\'eguliers, on dit qu'ils sont stablement conjugu\'es s'il existe $x\in \G(\C)$ tel que 
$xg x^{-1}=h$. On fait des d\'efinitions analogues sur $\Q$ et sur un corps $p$--adique,
cf.~\cite{Wald}. Par transport de structure on a donc d\'efini des \'el\'ements ``stablement $\theta$--conjugu\'es'', ``$\theta$--fortement r\'eguliers'' dans~$\G(\R)$.  

On sait d\'efinir de m\^eme des \'el\'ements fortement r\'eguliers, ainsi que la conjugaison stable, dans $H(\R)$
(Kottwitz~\cite{kott}).

Soit $g\in \G(\R)$ un \'el\'ement fortement $\theta$--r\'egulier, et soit $\Lambda(g)$ l'ensemble des valeurs propres (complexes) de $g^\theta\cdot g$ : celles--ci sont distinctes. Alors il existe $h\in H(\R)$ tel que l'ensemble des valeurs propres $\Lambda(h)$ de $h$ soit \'egal \`a $-\Lambda (g) \cup \{1\}$. On dit  que $h$ est une norme de $g$. 

\vskip2mm

\begin{notation}\label{notationnorme} $h=\mathcal{N}g$.
\end{notation}

\begin{prop}[Waldspurger]\label{normewalds}
La correspondance $\mathcal{N}$ d\'efinit une bijection entre classes de $\theta$--conjugaison 
stable d'\'el\'ements fortement $\theta$--r\'eguliers (dans $\G$) et classes de conjugaison stable d'\'el\'ements fortement r\'eguliers (dans~$H$).

\end{prop}

On dit qu'un \'el\'ement (fortement $\theta$--r\'egulier) de $\G$ est elliptique si sa norme est une classe elliptique de $H(\R)$. 
On peut alors la consid\'erer comme une classe de conjugaison dans $H_c(\R)$, conjugaison et conjugaison stable \'etant identiques 
dans un groupe compact.

\subsection{Stabilit\'e des caract\`eres tordus
$\Tpt$}\label{parstabilitetheta} Soit $\pi$ une repr\'esentation alg\'ebrique, r\'eguli\`ere, temp\'er\'ee et (donc) autoduale (\S\ref{prelimglnreel}) de $\G(\R)$. Le choix d'un op\'erateur d'entrelacement involutif $A:\mathfrak{H}_\pi \lgr\mathfrak{H}_\pi$ $(A^2 =1)$ entrela\c cant $\pi$ en $\pi\circ \theta$ permet d'\'etendre $\pi$ en une repr\'esentaton $\pi^+$ de $\G^+(\R)$. On note $\Tpt$ le caract\`ere de $\pi^+$ sur $\widetilde \G(\R)$ :
$$
\Tpt (g) = \Theta_{\pi^+}(\theta g)\qquad (g\in \G(\R))\,.
$$
A priori, $\Tpt$ est une distribution en $g$, mais un th\'eor\`eme d'Harish--Chandra (cf.~Bouaziz~\cite{bouaziz}) 
montre que c'est une fonction analytique sur les \'el\'ements (fortement) $\theta$--r\'eguliers.

\begin{theoreme}\label{stabilitetheta}
Pour un choix convenable de $A$, on a pour tout $g\in \G(\R)$ de norme fortement r\'eguli\`ere et elliptique :
$$
\Tpt (g) = \Theta_{\pi_H}(\mathcal{N} g)\,.
$$

En particulier, $\Tpt$ est invariant par conjugaison (tordue) stable sur les \'el\'ements de norme elliptique.
\end{theoreme}

Le th\'eor\`eme r\'esulte du travail de Bouaziz~\cite{bouaziz}, qu'il nous suffit d'interpr\'eter. Comme l'article de Bouaziz est 
\'ecrit dans un langage tr\`es diff\'erent de la th\'eorie usuelle des repr\'esentations admissibles 
(construction des repr\'esentations par la m\'ethode de Duflo--Kirillov) nous serons succincts, renvoyant le lecteur \`a
\cite{bouaziz} pour les d\'etails.

Donnons un ensemble explicite de repr\'esentants pour les classes de conjugaison tordue d'\'el\'ements fortement r\'eguliers. Dans $T_{H_c} \cong (T_a(\R))^n$, dont on repr\'esente les \'el\'ements par 
$(w_1,\ldots,w_n)$ : $w_i\in \C$, $|w_i|=1$, un \'el\'ement est fortement r\'egulier si $w_i\not= w_j$, $w_i^2 \not= 1$. 
Dans $\GL(2,\R)$, soit $J_2$ la matrice $\left(
\begin{array}{lll}
    0 &1   \\
  1&0      \\
\end{array}
\right)
$
 et consid\'erons le sous--groupe ${\rm O}(2,\R)=J_2 \SO(2,\R)\cup \SO(2,\R)$. Si $\theta$ est notre automorphisme usuel pour
$\GL(2)$, on a alors :
 $$
g\in \SO(2,\R) \Longrightarrow (g^\theta) g= (J_2 g J_2)g=1
$$
$$
g=J_2 h,\  h\in \SO(2,\R) \Longrightarrow g^\theta g=h^2\,.
$$

On en d\'eduit que la norme est dans ce cas surjective de $J_2 \SO(2,\R)$ vers $T_{H_c} \cong T_a(\R)$ ; si $R_\alpha$ est la rotation d'angle $\alpha$, elle envoie $J_2R_\alpha$ vers $-e^{2i\alpha}$.

Dans $\GL(2n,\R)$, on consid\`ere le tore $T_G \cong T_a(\R)^n$ donn\'e par les rotations d'angle $(\alpha_1,\ldots \alpha_n)$
{\it dans la base} $\mathcal{B}=(e_1,\,e_{2n}\,;\,e_2,\, e_{2n-1}\,; \ldots ;
e_n,\,e_{n+1})$. Si $R_\alpha$ $(\alpha=(\alpha_1,\ldots,\alpha_n))$ est une telle rotation, on a alors
$$
\mathcal{N}(J\, R_\alpha) = (-e^{2i\alpha_1}, -e^{2i\alpha_2}, \ldots, - e^{2i\alpha_n})
$$
vu comme \'el\'ement de $T_{H_c}$.

\`A la repr\'esentation $\pi$, la m\'ethode des orbites associe une forme lin\'eaire $f\in \got g^*$ o\`u $\got g
= \mathrm{Lie}\, (\G/\R) = {\rm M}_{2n}(\R)$. Elle est d\'efinie ainsi. Soit
$\got h$ la sous--alg\`ebre de Cartan de $\got g$ isomorphe \`a $\C^n$, obtenue \`a l'aide d'un plongement d'alg\`ebres 
$\C\subset M_2(\R)$ et de la base $\cal{B}$ ci--dessus. Alors $f: \got h \lgr \R$ est d\'efinie par 
$$(z_1,\ldots,z_n)\mapsto \sum_i p_i (z_i-\overline z_i) / \sqrt{-1}\, =\,
\sum_i
2\, p_i\,
\mathrm{Im}(z_i)
\qquad (z_i\in \C), $$  
c'est donc (au facteur $\sqrt{-1}$ pr\`es) la diff\'erentielle de la repr\'esentation induisant $\pi$. On \'etend $f$, de
fa\c{c}on naturelle, en une forme lin\'eaire sur $\got g$ par dualit\'e de Killing. Le stabilisateur
${\rm G}(f)$ est alors $T_\C\coprod \theta J T_\C$, $T_\C = (\C^\times)^n$ \'etant le tore exponentielle de $\got h$. 

La forme bilin\'eaire altern\'ee $B_f(X,Y) = f([X,Y])$ sur $\got g$ est non--d\'eg\'en\'er\'ee sur $\got g/\got
g(f)$, o\`u $\got g(f)= \mathrm{Lie}\,G(f)$. Soit $M$ le rev\^etement
m\'etaplectique de
$\Sp (\got g/\got
g(f))$. On obtient alors une extension centrale
$$
1 \lgr \{1,\varepsilon\} \lgr \widetilde \G (f) \lgr {\rm G}(f) \lgr 1
$$
(\cite[p. 47]{bouaziz} ; Bouaziz note diff\'eremment $\widetilde {\rm G}(f)$). L'application de chacun 
des facteurs $\C_u^\times = \{w:|w|=1\}\subset \C^\times$ de $T_\C$ dans $\Sp(\got g/\got g(f))$ est quadratique
({\it i.e.} , passe au quotient par $(\C_u^\times)/\pm1$. On en d\'eduit que ce rev\^etement est scind\'e au--dessus de 
$T_\C$, puis au--dessus de $T_\C \coprod \theta\cdot J\cdot T_\C$. Il existe donc un caract\`ere unitaire $\tau$ de $\widetilde
{\rm G}(f)$ tel  que
\[
\begin{array}{lll}
  &\tau(\varepsilon) = -1   &   \\
  &\tau(z) =\chi(z)   &(z\in T_\C )  \\
\end{array}
\]
o\`u $\chi$ est le produit des $\chi_i= z^{p_i}(\overline z)^{-p_i}$.

Soit $s=s(\alpha_1,\ldots, \alpha_n)\in \widetilde \G$, $s = \theta\cdot JR_\alpha$, $R_\alpha$ \'etant la matrice de
rotation ci--dessus. Soit $W$ le groupe de Weyl de $(\G,T_\C)$ : il s'identifie \`a $\got{S}_{n} \times
\{\pm1\}^n:=W_H$. Par ailleurs la donn\'ee de $\tau$ d\'efinit une extension $\pi^+$ de $\pi$ \`a $\G^+$. Alors, d'apr\`es
\cite[Prop. 6.1.2]{bouaziz},

\begin{equation}\label{formulethetabouaziz}
\Theta_{\pi^+}(\theta J R_\alpha) = \Theta_{\pi,\theta}(JR_\alpha) = \sum_{w\in W_H} {\varepsilon(w)\ {}^w
\tau(s)\rho_{\got l}(s)\over D(s)}\,.
\end{equation}

La valeur de $\rho_{\got l}(s)$ et de $D(s)$ d\'epend du choix d'un sous--espace
lagrangien $\got l \subset \got g \otimes \C$ (pour $B_f$) stable par $s$ ; 
soit $\got g =\got g^s \oplus \got q$ la d\'ecomposition de $\got g$ selon les espaces propres de $s$.~Alors
$$
D(s) =\det (1-s) |_{\got l \cap \got q}\,.
$$

Il est clair que $D(s)$ est invariant par conjugaison stable (pour des choix compatibles de $\got l$). 
Par transport de structure on le calcule sur un \'el\'ement de 
$\GL(2n,\C)$ de la forme $\theta x$~o\`u
\[ x\,=\,
\left(
\begin{array}{cccccc}
  x_1 &   &   \\
  & \ddots  &   \\
 &   & x_n  \\
 &&& x_n^{-1}\\
 &&&&\ddots\\ 
 &&&&&  x_1^{-1}\\
\end{array}
\right)\,.
\]

On prend $\got l$ \'egale au radical unipotent de l'alg\`ebre de Borel standard de
${\rm M}_{2n}(\C)$ : on calcule les valeurs propres sur $\got l$ de l'endomorphisme
$$
X \longmapsto -J \ Ad(x)\ {}^t X J\,.
$$

Un calcul simple donne alors
$$
D(x) =\prod_{i=1}^n (1+x_i^2 ) \prod_{i<j} (1-x_i^2 /x_j^2) \prod_{i<j}(1-x_i^2 x_j^2)
$$
qui n'est autre que le d\'enominateur de Weyl $\prod(1-e^\alpha)$ (le
produit portant sur les racines positives de $H$) \'evalu\'e en $(-x_i^2)=\mathcal{N}x$. Le terme
$\rho_{\got l}(s)^{-1}$ le transforme en $\prod(e^{-\alpha/2}- e ^{\alpha/2}):= D_H$. Par ailleurs
$$\tau(s) =\tau(\theta J)\tau(R_\alpha) = \tau(\theta J) \prod e^{2\sqrt{-1} \alpha_i p_i} = \varepsilon \chi_H(\mathcal{N}(JR_\alpha))
$$ o\`u $\chi_H$ est le caract\`ere de $T_H$ param\'etr\'e par $(p_1,\ldots p_n)$ et $\mathcal{N}(JR_\alpha) 
= (-e^{2\sqrt{-1}\alpha_i})$ (noter que $p$ param\`etre le caract\`ere infinit\'esimal de $\pi_H$, donc $p=m +\rho_H$ o\`u $m$ 
est le plus haut poids. L'ind\'etermination due au fait que $p$ n'est pas un caract\`ere dispara\^it puisqu'on \'evalue en des carr\'es). Enfin, 
on v\'erifie que $\varepsilon(w)$ co\"incide, au signe pr\`es, avec le d\'eterminant sur $W_H$ (voir
\cite[(5.5.1) et p.~46]{bouaziz}). L'expression (\ref{formulethetabouaziz}) donne~alors
$$
\Theta_{\pi,\theta}(x) = \varepsilon \sum_w {\varepsilon(w) \mathcal{N}(x)
^{w\chi_H}\over D(\mathcal{N}(x))}\,,\ \ \varepsilon=\pm1
$$
d'o\`u, d'apr\`es la formule de Weyl, le Th\'eor\`eme~2.2.

\subsection{Pseudocoefficients et int\'egrales orbitales
tordues}\label{pseudoorb}

La repr\'esentation $\pi$ de $\G(\R)$ reste fix\'ee (pour simplifier on note simplement $\G$ le groupe $\G(\R)$). 
On v\'erifie ais\'ement que $\pi$ est $\theta$--discr\`ete, {\it i.e.} , isol\'ee parmi les repr\'esentatons temp\'er\'ees, $\theta$--invariantes de $\G$. D'apr\`es le 
th\'eor\`eme de Paley--Wiener de Mezo \cite{Mezo}, il existe une fonction $f_\pi\in \Cc^\infty(\G)$, $K_\infty$--finie (pour un sous--groupe compact maximal
$K_\infty$ de
$\G(\R)$) telle que
\begin{equation}\label{pseudotordu1} \hbox{trace}(\pi(f_\pi)A) =
1\end{equation}
et
\begin{equation}\label{pseudotordu2} \hbox{trace}(\rho(f_\pi)A_\rho)=0.
\end{equation}
\noindent pour toute repr\'esentation temp\'er\'ee, irr\'eductible, $\theta$--invariante de $\G$ ; 
on note $A$ l'op\'erateur d'entrelacement entre $\pi$ et $\pi\circ \theta$, normalis\'e par le
Th\'eor\`eme~\ref{stabilitetheta}, 
et $A_\rho$ un op\'erateur d'entrelacement (non nul) entre $\rho$ et $\rho \circ
\theta$.\ps

Soit $\gamma\in \G$ un \'el\'ement $\theta$--semisimple. Son centralisateur tordu $I$,  \'egal par d\'efinition \`a 
la composante neutre de $I'=\{ g\in \G : g^{-1}\gamma g^\theta =\gamma\}$ est alors r\'eductif. 
On consid\`ere l'int\'egrale orbitale tordue (pour des mesures de Haar $dg$ et $di$ arbitraires)
$$
\TO_\gamma(f) = \int_{I\backslash \G} f(g^{-1} \gamma g^\theta) {dg\over d i}
$$
pour $f\in C_c^\infty(\G)$. Si $\gamma$ est fortement r\'egulier, $I=I'$ est un tore.

Soit $P=MN\subset \G$ un parabolique $\theta$--stable et soit $\gamma\in M$ un \'el\'ement 
fortement $\theta$--r\'egulier. Alors $\gamma$ a la m\^eme propri\'et\'e relativement \`a $M$, 
et son centralisateur tordu est un tore de $M$. Si $f\in C_c^\infty(\G)$ soit 
$\bar f(x) = \int_{K_\infty} f(k^{-1} x k^\theta)dk$ 
(pour la mesure de Haar normalis\'ee). Alors
$$
\TO_\gamma (f) = \int_{I\backslash MN} \bar f(n^{-1} m^{-1} \gamma m^\theta n^\theta) {dmdn\over di}\,.
$$
Si $h \in {\cal C}_c^\infty(P)$ et $m\in M$, 
$$
\int_N h (n^{-1} m \,n^\theta m^{-1})dn = D(m)^{-1} \int_N h(n)dn
$$
o\`u $D(m) = |\det (1- {Ad}(m)\circ \theta)|_{\got n}$ et $\got n =\mathrm{Lie}(N)$ ;  $D(m)$ est non nul si $m$ est $\theta$--r\'egulier. Ainsi
\begin{equation}\label{descenteleviTO}
\TO_\gamma (f) = D(\gamma)^{-1} \int_{I\backslash M}\overline f^{(P)} (m^{-1} \gamma m^\theta){dm\over di} = D(\gamma)^{-1}
\TO_\gamma ^M (\overline f^{(P)}) 
\end{equation}
o\`u l'int\'egrale orbitale tordue est prise dans $M$ et $f^{(P)}$ est d\'efinie selon Harish--Chandra par
$$
f^{(P)}(m) = \int_N f(nm)dn\,.
$$

\begin{lemme}\label{caracellreg}
Soit $\gamma$ un \'el\'ement fortement $\theta$--r\'egulier non--elliptique de $\G$. Alors $\gamma$ est $\theta$--conjugu\'e \`a un \'el\'ement (fortement $\theta$--r\'egulier) de la composante de Levi $M$ d'un parabolique $\theta$--stable propre de~$G$.
\end{lemme}

\begin{pf} Soit en effet $\delta = \gamma^\theta \gamma$. 
C'est un \'el\'ement r\'egulier de $\GL(2n,\R)$ dont l'ensemble des valeurs propres 
(distinctes et complexes) est auto--dual et dont au moins une valeur propre $\lambda \in \C^{\ast}$ n'est pas de module $1$. 
Posons $i=1$ ou $2$ selon que $\lambda$ est r\'eel ou non et consid\'erons le parabolique sup\'erieur standard de type $(i,2n-2i,i)$ de $\G$ ; 
il est $\theta$--stable. \`A conjugaison pr\`es sur $\delta $  
(et donc \`a $\theta$--conjugaison pr\`es sur $\gamma$), on peut supposer que $\delta$ est un \'el\'ement du sous-groupe de Levi standard $M$ 
de ce parabolique :
$$\delta =\left(
\begin{array}{ccc}
\lambda  &   &   \\
  & \ast   &   \\
  &   & \lambda^{-1}
\end{array}
\right)$$
(si $i=2$ on choisit un plongement de $\C$ dans ${\rm M}_2(\R)$). Comme $\delta^{\theta}\gamma=\gamma\delta$ et $\delta$ est r\'egulier,
cela impose que $\gamma \in M$.
\end{pf}

\begin{lemme}\label{annuleregnonell}
Si $\gamma$ est fortement $\theta$--r\'egulier et non--elliptique,
$$
\TO_\gamma (f_\pi ) =0\,.
$$
\end{lemme}

\begin{pf} D'apr\`es (\ref{descenteleviTO}) cette int\'egrale orbitale se calcule dans $M$, 
o\`u $P=MN$ est un parabolique propre, $\theta$--stable et $\gamma \in M$. D'apr\`es \cite{KoR} 
 -- et gr\^ace au th\'eor\`eme de Paley--Wiener de Mezo d\'ej\`a cit\'e -- 
une fonction $h$ sur $M$ a des int\'egrales orbitales tordues nulles si
${\rm trace}\, (\pi_M(h)A)=0$ pour toute repr\'esentation $\theta$--stable temp\'er\'ee 
$\pi_M$ de $M$, $A\not= 0$ \'etant un op\'erateur entrela\c{c}ant $\pi_M$ et $\pi_M\circ \theta$. 
Pour $h=\overline f^{(P)}$, un lemme bien connu d'Harish--Chandra 
(qui s'\'etend formellement au cas tordu) donne
$$
\hbox{trace } (\pi_M(h)A) = \hbox{trace }(\pi_G(f)A_G)
$$
o\`u $\pi_G$ est induite de $\pi_M$ et $A_G$ est l'op\'erateur d'entrelacement induit.
Mais l'expression de droite s'annule d'apr\`es~(\ref{pseudotordu2}).
\end{pf}

Soit $\gamma\in \G$ un \'el\'ement $\theta$--semisimple. On dit que $\gamma$ 
est $\theta$--elliptique si la composante d\'eploy\'ee 
du centre de son centralisateur tordu est r\'eduite \`a l'\'el\'ement neutre. 
L'\'el\'ement $\delta =\gamma^\theta\gamma$ 
de $\GL(2n,\R)$ est conjugu\'e \`a un \'el\'ement diagonal de $\GL(2n,\C)$ de 
la forme $(x_1,\ldots, x_n, x_n^{-1},\ldots,x_1^{-1})$. Il d\'efinit donc 
une classe de conjugaison $\mathcal{N}\gamma$ dans $\SO(2n+1,\C)$, par son spectre
$$
\Lambda (\mathcal{N}\gamma) = \{ -x_1,\ldots,-x_1^{-1} \} \cup \{1\}\,.
$$
\noindent On v\'erifie que $\mathcal{N} \gamma$ (qui est toujours conjugu\'e \`a un \'el\'ement de $\SO^*(2n+1,\R)$) est 
conjugu\'e \`a un \'el\'ement, unique \`a conjugaison pr\`es de $\SO(2n+1,\R)$ si, 
et seulement si, $\gamma$ est $\theta$--elliptique.\footnote{Posons
$\delta=\gamma^\theta\gamma$. La conjugaison
par $\theta\gamma$ induit une involution sur $\G_{\delta}$ 
(un sous-groupe de Levi de $\GL_{2n}$) dont le sous-groupe des points fixes est le centralisateur 
tordu $\G_\gamma$ de $\gamma$. La structure des anti-involutions des
alg\`ebres semisimples complexes montre que le centre $Z$ de $\G_\gamma$ co\"incide
avec le sous-groupe des $\theta\gamma$-invariants du centre $Z'$ de $\G_\delta$. 
Or $Z'=\GG_mD$ o\`u $D$ d\'esigne l'adh\'erence Zariski du sous-groupe
engendr\'e par $\delta$ (central, fix\'e par $\theta\gamma$). Ainsi
$Z={Z'}^{\gamma\theta}=\GG_m^{\theta\gamma}D=\{\pm 1\}D$, ce qui conclut.}

\begin{lemme}\label{TOnonell}
Si $\gamma\in \G$ est $\theta$--semisimple mais non $\theta$--elliptique,
$\TO_\gamma(f_\pi)=0$.
\end{lemme}

\begin{pf}
Soit en effet $I$ le centralisateur tordu de $\gamma$. Une fonction $f\in C_c^\infty(\G)$ \'etant donn\'ee, 
on peut trouver une fonction $h\in {\cal C}_c^\infty(I)$ telle que, pour $x\in I$, voisin de~1,
$$
\TO_{x\gamma}(f) = {\rm O}_x^I(h)
$$
le membre de droite \'etant une int\'egrale orbitale ordinaire dans $I$ (voir la preuve de la Proposition~\ref{connoncon} pour l'argument). Si $x\gamma$ 
est fortement r\'egulier, $x$ est r\'egulier dans $I$ ; pour $x$ assez
proche de $1$, les 
composantes neutres du centralisateur de $x$ et du centralisateur tordu de $x\gamma$ 
co\"incident (\cite[Cor.~3.1.5]{labesse}). Puisque la composante d\'eploy\'ee du centre de $I$ est 
non--triviale, les int\'egrales orbitales ${\rm O}_x^I(h)$ sont donc nulles (au voisinage de
$1$) en les \'el\'ements r\'eguliers et donc $h(1)=\TO_\gamma(f)$ (pour des normalisations 
convenables des mesures) s'annule.
\end{pf}
Si $\gamma\in \G$ est $\theta$--semisimple, l'int\'egrale orbitale tordue stable de 
$f$ associ\'ee \`a $\gamma$ est d\'efinie par Labesse \cite[\S2.7]{labesse}. On la note
$$
\STO_\gamma(f)\,.
$$

\subsection{Int\'egrales orbitales tordues en les \'el\'ements elliptiques}\label{parfinchap2}
Nous pouvons maintenant d\'emontrer le r\'esultat principal de ce chapitre.

\begin{theoreme}\label{staborbreel} Soit $\gamma$ un \'el\'ement $\theta$--semisimple de $\G(\R)$.

\begin{itemize}
\item[(i)] Si $\gamma$ n'est pas $\theta$--elliptique, $\TO_\gamma(f_\pi)= \STO_\gamma(f_\pi)=0$.

\item[(ii)] Soit $\gamma$ un \'el\'ement $\theta$--elliptique, et $I=I_\gamma$. Pour un choix convenable 
de mesures positives sur $\G$ et $I_\gamma$,
\end{itemize}
$$
\TO_\gamma(f_\pi) = e(\gamma) \Theta_{\pi_H}(\mathcal{N}\gamma)\,,
$$
$\pi_H$ \'etant la repr\'esentation (de dimension finie) de $\SO(2n+1,\R)$ associ\'ee \`a $\pi$, et $e(\gamma)=\pm1$ 
un signe ind\'ependant de~$\pi$. En particulier, ces int\'egrales orbitales sont stables.
\end{theoreme}

Nous ne d\'ecrirons pas la normalisation des mesures, qui r\'esulte de \cite[\S A.1]{labesse}. L'important est que, $\gamma$ 
\'etant fix\'e, celle--ci ne d\'epend pas de~$\pi$.\ps

La partie (i) a d\'ej\`a \'et\'e d\'emontr\'ee. Consid\'erons d'abord le cas o\`u $\pi_0$ est l'unique repr\'esentation 
temp\'er\'ee de $\G(\R)$ ayant de la cohomologie \`a coefficients triviaux : c'est celle d\'efinie par 
$p(\pi_0) = \Big({2n-1\over2}\,,\, {2n-3\over2},\ldots , {1-2n\over2}\Big)$ (\S\ref{prelimglnreel}). Dans ce cas Labesse \cite{labesse} 
a donn\'e une construction de $f_{\pi_0}$ par voie cohomologique, qui permet d'en calculer les int\'egrales orbitales tordues 
\cite[Thm.~A.1.1]{labesse}, 
ce qui d\'emontre le th\'eor\`eme dans ce cas, $\pi_{0,H}$ \'etant  la repr\'esentation 
triviale\footnote{Dans cet article notre $I$ est remplac\'e par le \og $\theta$--centralisateur stable\fg\,, cf. \cite[p.~52]{labesse}. 
La description des centralisateurs tordus par Waldspurger \cite[\S I]{Wald} montre qu'ils co\"incident pour $\GL(2n)$ 
avec nos centralisateurs connexes.}. Nul doute que l'on pourrait l'\'etendre au cas g\'en\'eral : comme nous avons d\'emontr\'e des r\'esultats 
plus puissants, nous donnons une d\'emonstration diff\'erente. \ps

Pour $\pi$ donn\'ee, soit $(\rho,V)$ la repr\'esentation alg\'ebrique ($\theta$--stable) telle que la cohomologie de $\pi$ \`a 
coefficients dans $V$ soit non nulle ; on d\'efinit de m\^eme $(\rho_H,V_H)$, $\rho_H$ s'identifiant dans ce cas \`a $\pi_H$. 
On peut r\'ealiser $V$ \`a l'aide du th\'eor\`eme de Borel--Weil, dans la cohomologie de $\G(\C)/B(\C)$ \`a coefficients dans le 
fibr\'e en droites $\mathcal{L}_m$ o\`u $m=m(\pi)$. Noter que $m$ est invariant par $\theta$ ;  on en d\'eduit naturellement une 
repr\'esentation de  $\G^{+}(\C)$. Pour $g\in \G(\C)$, $\theta$--r\'egulier, la trace $(g\times \theta\mid V)$ se calcule 
\`a l'aide du th\'eor\`eme d'Atiyah--Bott ; un calcul simple montre que les points fixes sont param\'etr\'es par le centralisateur 
de $\theta$ dans $W(\G(\C))\cong \got S_{2n}$, isomorphe \`a $W_H$. On obtient alors
\begin{equation}\label{carrhotordu}
\hbox{trace }(g\times \theta \mid V) = \hbox{trace }(\,\mathcal{N} g\mid V_H\,)\,.
\end{equation}
Soit alors $\Theta_{\rho,\theta} $ le caract\`ere tordu de $\rho$ (pour ce choix d'op\'erateur d'entrelacement), 
et soit $$g_\pi=\Theta_{\rho,\theta}f_0.$$ Alors $g_\pi$ a les propri\'et\'es du Th\'eor\`eme~\ref{staborbreel} relativement \`a $\pi$. Pour montrer que 
$g_{\pi}$ a les m\^emes int\'egrales orbitales que $f_{\pi}$, il suffit de montrer que pour toute repr\'esentation temp\'er\'ee $\theta$--stable 
$\tau$ (et op\'erateur d'entrelacement $A_{\theta}$ associ\'e) :
\begin{equation}\label{identitetrace} \trace\,(\,\tau(f_{\pi}) \,A_{\theta}\,)\, = \, \trace\, (\,\tau(g_{\pi})\, A_{\theta}\,)
\end{equation}
(le fait que ceci implique l'\'egalit\'e des int\'egrales orbitales tordues est le th\'eor\`eme de densit\'e de Kottwitz-Rogawski \cite{KoR}).\ps

On dit que $\tau$ est $\theta$--discr\`ete si elle n'est pas induite d'une repr\'esentation $\theta$--stable d'un parabolique $\theta$--stable. 
Dans ce cas le caract\`ere tordu est \`a support dans les \'el\'ements non $\theta$-elliptiques, donc (\ref{identitetrace}) est \'evidente puisque 
les int\'egrales orbitales tordues correspondantes sont nulles. Les repr\'esentations $\theta$--discr\`etes sont, comme on le voit 
ais\'ement, de la forme
\begin{equation}\label{thetadisc1} \tau\,=\,{\rm Ind} (\delta_{1},\ldots,\delta_{n}), \end{equation}
o\`u $\delta_{i}$ est une repr\'esentation de $\GL(2,\R)$, de la s\'erie discr\`ete associ\'ee \`a la repr\'esentation de $\W_{\R}$ induite d'un 
caract\`ere $z \mapsto z^{p_{i}}(\bar z)^{{-p_{i}}}$ de $\W_{\C}=\C^{*}$, avec $p_{i} \in \frac{1}{2}\Z$, les $p_{i}$ \'etant distincts, 
ou bien de la forme
\begin{equation}\label{thetadisc2} \tau\,=\,{\rm Ind}(\delta_{1},\ldots,\delta_{n-1},1,\varepsilon),\end{equation}
les $\delta_{i}$ \'etant comme auparavant, et $\varepsilon$ \'etant le caract\`ere d'ordre $2$ de $\R^{\ast}$.\ps
Si les $p_{i}$ appartiennent \`a $\frac{1}{2}+\Z$ (et donc $\tau$ est cohomologique), on a d'une part 
$$\trace\, (\,\tau(f_{\pi})A_{\theta}^{\tau}\,)\,= \delta\,(\tau,\pi)$$
o\`u $\delta$ est le symbole de Kronecker et $f_{\pi}$ est normalis\'ee par $A_{\theta}^{\pi}$ ; par ailleurs
\begin{equation} \label{identitetrace2} \trace\,(\,\tau(g_{\pi}) A_{\theta}^{\tau}\,)\,=\, \int_{\G} 
\Theta_{\tau,\theta}(g)\,\Theta_{\rho,\theta}(g)\, f_{0}(g)\,dg\,.\end{equation} 
Mais les int\'egrales orbitales tordues de $f_{0}$ s'annulent pour $g$ non $\theta$--elliptique ; si $g$ est $\theta$--elliptique r\'egulier 
(donc de centralisateur tordu ${\rm U}(1)^{n}=T$) et si la mesure sur $T$ est normalis\'ee, on trouve (\cite[Thm. A.1.1]{labesse})
$$\TO_g(f_{0})=1$$
($f_{0}$ \'etant bien s\^{u}r le pseudocoefficient associ\'e \`a une mesure $dg$ qui est celle d\'efinissant l'int\'egrale obitale). Donc 
(\ref{identitetrace2}) s'\'ecrit $$\int_{\G} \Theta_{\tau,\theta}(g)\,\Theta_{\rho,\theta}(g)\,f_{0}(g)\,dg \,\, = \, \, \frac{1}{|W|} \, \int_{T} 
\{ \sum_{\No\delta=\gamma} \Theta_{\tau,\theta}(\delta)\,\Theta_{\rho,\theta}(\delta)\,\Delta(\gamma) \} d\gamma,$$
o\`u $\Delta(\gamma)$ est un d\'enominateur de Weyl (pour la formule d'int\'egration de Weyl relative \`a la conjugaison tordue), qu'on 
v\'erifie \^{e}tre \'egal \`a un facteur $2^{n}$ pr\`es
 (le nombre de $\delta$ de norme $\gamma$) au d\'enominateur de Weyl pour $\SO(2n+1)$. Le Th\'eor\`eme~\ref{stabilitetheta} et la 
 relation (\ref{carrhotordu}) impliquent alors que 
 $\trace\,(\tau(g_{\pi}) \,A_{\theta}^{\tau}))= \delta (\tau, \pi)$, d'apr\`es les relations d'orthogonalit\'e sur $\SO(2n+1)$.\ps
 Consid\'erons enfin les autres repr\'esentations $\tau$ de type (\ref{thetadisc1}) ou (\ref{thetadisc2}). On a 
 $$\begin{array}{ccl} \trace\,(\,\tau(g_{\pi})A_{\theta}^{\tau}\,) & = & \trace\,\left(\,\int_{\G}\tau(x)\,g_{\pi}(x)\,A_{\theta}^{\tau}
 \,\right)\,dx\\
 & = & \trace\,\left(\,\int_{\G}\,\tau(x)\,f_{0}(x)\,\trace\,(\,\rho(x)A_{\theta}^{\rho}\,)\,A_{\theta}^{\tau}\,\right)\,dx\\
 & = & \trace\,\left(\,\int_{\G} f_{0}(x)\,(\,\tau(x)\otimes \rho(x)\,)\,A_{\theta}^{\tau}\otimes A_{\theta}^{\rho}\,\right)\,dx\,.
 \end{array}$$
 Si $\tau$ est du type indiqu\'e et non cohomologique, son caract\`ere infinit\'esimal $\lambda$ (la somme des $p_{i}$ 
 et des $-p_{i}$, avec $p=0$ pour les caract\`eres $1,\varepsilon$) n'appartient pas \`a $\{\frac{1}{2}+\Z\}^{2n}$. On sait que les caract\`eres infinit\'esimaux des sous-quotients de $\tau \otimes \rho$ sont de la forme 
 $\lambda +\mu$ o\`u $\mu$ est un poids (entier) de $\rho$. Il a donc la m\^{e}me propri\'et\'e ; la trace de $f_{0}$ dans $\tau \otimes \rho$ est donc nulle, ce qui termine la d\'emonstration. $\square$

\ps\ps

Si $T \subset \G$ est le centralisateur tordu d'un \'el\'ement $\theta$-semisimple fortement r\'egulier et si ${\cal D}(T,\G,\R)$ classifie 
la conjugaison stable (modulo conjugaison) pour les \'el\'ements de $T$, on a ${\cal D}(T,\G,\R)=H^{1}(\R,T)$ puisque $G$ est 
cohomologiquement trivial. Avec les notations de Labesse \cite[p. 122-123]{labesse}, 
$H^{1}(\R,T)={\cal E}(T,G,\R)$ ; si $\kappa$ est un \'el\'ement du dual de $\cal E$, on en d\'eduit:

\begin{corollaire}\label{corollairekappa} Si $\kappa\neq 1$
$$\TO_{\delta,\kappa}(f_{\pi}) = \sum_{d \in {\cal D}(T,\G,\R)} \langle \kappa, d\rangle \TO_{d\delta}(f_{\pi})=0.$$
\end{corollaire}
Un argument de descente (ibid.) montre que le m\^{e}me r\'esultat d'annulation s'\'etend \`a tous les \'el\'ements $\delta$ de norme elliptique. La fonction $f_{\pi}$ est donc {\it stabilisante} au sens de Labesse.
\bigskip
\bigskip
\section{Asymptotique de la formule des traces tordues de $\GL(2n)$}
\setcounter{equation}{0}
Dans ce chapitre, nous d\'emontrons le Th\'eor\`eme $B$ \'enonc\'e
dans l'introduction en admettant le Th\'eor\`eme $E$.

\subsection{\'Enonc\'e du th\'eor\`eme}\label{ssectenonce} Soient $\ell$ et $p$ deux nombres premiers distincts, $n\geq 1$ un entier. On consid\`ere \`a 
nouveau le groupe alg\'ebrique $G:=\GL_{2n}$ sur $\Q$ muni de son $\Q$-automorphisme d'ordre $2$ 
$$\theta(g)=J {}^t\!g^{-1} J^{-1},$$
o\`u $J=J^{-1} \in \GL_{2n}(\Q)$ est la matrice antidiagonale d\'efinie au
\S\ref{prelimglnreel}. 
On d\'efinit encore $G^+$ comme \'etant le $\Q$-groupe
produit semi-direct de $\Z/2\Z=\langle \theta \rangle$ par $G$ d\'efini par $\theta$.
\par \medskip

\`A la place $p$, fixons une repr\'esentation irr\'eductible 
supercuspidale $\omega$ de $\GL_n(\Q_p)$ dont la contragr\'ediente
$\check{\omega}$ n'est isomorphe \`a aucune tordue non ramifi\'ee de $\omega$. Soient $P$ le
parabolique triangulaire sup\'erieur de $G$ de type $(n,n)$ et $M=\GL_n \times
\GL_n$ son sous-groupe de Levi diagonal. Si $\chi$ est un
caract\`ere non ramifi\'e de $\GL_n(\Q_p)$, on notera $I(\chi)$ l'induite
parabolique normalis\'ee de $P$ \`a $G$ de la repr\'esentation 
$$(\omega\otimes\chi) \times (\check{\omega}\otimes\chi^{-1})$$ de $M(\Q_p)$.
Ces repr\'esentations seront \'etudi\'ees en d\'etail dans un paragraphe
ul\-t\'e\-ri\-eur (\S\ref{pseudocoeffcpomega}). Disons simplement ici que les $I(\chi)$ sont irr\'eductibles
et autoduales, de sorte qu'elles se prolongent (en fait de mani\`ere naturelle
car les donn\'ees ci-dessus sont $\theta$-stables) \`a $G^+(\Q_p)$. \par \medskip

L'objectif principal de ce chapitre est de d\'emontrer le r\'esultat suivant.

\begin{theoreme}\label{mainthm} Il existe une repr\'esentation 
automorphe cuspidale irr\'eductible $\Pi$ de $\GL_{2n}(\AAA)$ ayant les propri\'et\'es suivantes:\begin{itemize}
\item[i)] $\check{\Pi} \simeq \Pi$,
\item[ii)] $\Pi_{\infty}|\cdot|^{(2n-1)/2}$ est alg\'ebrique r\'eguli\`ere,
\item[iii)] $\Pi$ est non ramifi\'ee \`a toutes les places finies diff\'erentes de $\ell$ et $p$,
\item[iv)] $\Pi_p \simeq I(\chi)$ pour un certain $\chi$,
\item[v)] $\Pi_{\ell}$ est la repr\'esentation de Steinberg.
\end{itemize}
\end{theoreme}

Afin d'utiliser la formule des traces d'Arthur, nous aurons besoin de certaines propri\'et\'es des int\'egrales orbitales tordues des pseudocoefficients 
n\'ecessaires en $\infty$, $\ell$ et $p$. Le travail archim\'edien a d\'ej\`a \'et\'e fait dans
le chapitre pr\'ec\'edent. Une premi\`ere sous-section \S\ref{preliminaireLef} sera consacr\'ee 
\`a l'\'etude des pseudocoefficients 
$\theta$-tordus de la repr\'esentation de Steinberg, ce qui fournira les renseignements n\'ecessaires en $\ell$. En ce qui concerne la place $p$, 
nous admettrons au \S\ref{preuvemainthm} le r\'esultat suivant dont la preuve fera l'objet du chapitre $4$.

Consid\'erons l'el\'ement $\theta$-semisimple elliptique $$\gamma_0:= \left( \begin{array}{cc} 1_n & 0 \\ 0 &
-1_n\end{array}\right) \in G(\Q).$$

\begin{theoreme}\label{thmcle} Il existe une fonction $f_p: G(\Q_p)
\longrightarrow \C$ localement constante \`a support compact ayant les propri\'et\'es suivantes:
\begin{itemize}
\item[(i)] Si $\pi$ est une repr\'esentation lisse irr\'eductible autoduale de
$G(\Q_p)$ et si $A: \pi \isomo \pi \circ \theta$ est un $G(\Q_p)$-isomorphisme, alors
$\trace(A\pi(f))=0$ si $\pi$ n'est pas de la forme $I(\chi)$,
\item[(ii)] L'int\'egrale orbitale tordue $$\TO_{\gamma_0}(f_p)=\int_{I_{\gamma_0}(\Q_p)\backslash G(\Q_p)}
f_p(g^{-1}\gamma_0 \theta(g)) \mu$$
est non nulle.\footnote{$\mu$ est ici une mesure
$G(\Q_p)$-invariante non nulle quelconque sur $I_{\gamma_0}(\Q_p)\backslash G(\Q_p)$, $I_{\gamma_0} \simeq \Sp_{2n}$ \'etant le centralisateur de $\gamma_0$.}
\end{itemize}
\end{theoreme}

\subsection{Fonctions d'Euler-Poincar\'e et repr\'esentation de Steinberg dans le cas tordu}

\label{preliminaireLef} Dans cette partie, nous rappelons 
ou \'etablissons certaines propri\'et\'es des fonctions d'Euler-Poincar\'e 
associ\'ees \`a un automorphisme d'un groupe r\'eductif connexe. Les \'enonc\'es \'etant ici tout \`a fait g\'en\'eraux, nous nous pla\c{c}ons dans tout ce \S \ref{preliminaireLef} dans le cadre suivant.

\subsubsection{} \label{notationslef}
On fixe $F$ un corps local non archim\'edien de caract\'eristique nulle, $G$ un groupe r\'eductif connexe sur $F$ et $\theta$ un 
automorphisme $F$-rationnel de $G$ disons d'ordre fini\footnote{En fait, l'hypoth\`ese que $\theta$ est d'ordre fini ne sera utilis\'ee 
que dans le point (ii) de la 
Proposition \ref{nbleftheta}, pour laquelle nous devons supposer que $G^+$ est un groupe
alg\'ebrique lin\'eaire.} $h$. On note $G^+$ le 
produit semi-direct de $\Z/h\Z$ par $G$ d\'efini par $\theta$, $S$ le 
$F$-tore maximal central d\'eploy\'e de $G$ et $X^*(S)=\Hom(S,{\mathbb G}_m)$ le groupe ab\'elien libre des caract\`eres rationnels de $S$. 
Le groupe quotient $G^+/G=\langle \theta \rangle$ agit par conjugaison sur $S$ ; on peut donc 
consid\'erer le plus grand 
sous-tore $S^\theta$ de $S$ fix\'e par $\theta$ ainsi que son polyn\^ome caract\'eristique r\'eciproque 
$$P(z):=\det(1-z\theta_{|X^*(S)}) \in \Z[z].$$ 
On d\'efinit de plus $q(G)\geq 0$ comme 
\'etant le rang d'un $F$-tore d\'eploy\'e maximal du groupe $G/S$.

Nous noterons $1$ la repr\'esentation triviale de $G^+(F)$ et $\St $ la repr\'esentation de Steinberg.
Rappellons que cette derni\`ere est d\'efinie comme suit. Soit $B$ un parabolique minimal de $G$ d\'efini sur $F$. Tous les tels
\'etant $G(F)$ conjugu\'es il existe un \'el\'ement $\theta_B \in G(F)\theta$ tel que
$\theta_B(B)=B$.
Notons $I_B$ l'espace des fonctions complexes lisses sur $B(F)\backslash G(F)$: c'est une repr\'esentation de $G(F)$ 
par translations \`a droite qui se prolonge naturellement \`a $G^+(F)$. La repr\'esentation de Steinberg $St$ de $G^+(F)$ 
est alors l'unique quotient irr\'eductible de $I_B$.
Par exemple, $\St=1$ si $G$ est un tore. 
De plus, $\St$ se factorise par
$G(F)/S(F)=(G/S)(F)$ (via Hilbert 90): c'est la repr\'esentation de Steinberg de $G/S$, et elle est de carr\'e int\'egrable si $S=1$.

Enfin, les $F$-tores d\'eploy\'es maximaux de $B$ \'etant conjugu\'es sous $B(F)$, on peut supposer 
quitte \`a remplacer $\theta_B$ par un \'el\'ement de la forme $b\theta_B$ pour $b \in B(F)$ 
que $\theta_B$ stabilise un tel tore $S'$ de $B$, auquel cas il pr\'eserve aussi le centralisateur de $S'$ dans $B$ qui est un 
sous-groupe de Levi de $B$. L'automorphisme $\theta_B$ agit sur
$\Lambda^{q(G)}(X^*(S'/S))$ par un signe que l'on note $\epsilon(\theta)$ 
(il d\'ecoulera par exemple de la proposition suivante que ce signe ne d\'epend d'aucun choix).

\begin{exemple}\label{exgl2n}{\rm Dans le cas qui nous int\'eresse $F=\Q_l$, 
$G=\GL_{2n}$ et $\theta$ est comme au \S \ref{ssectenonce}. On a $S=\mathbb{G}_m$ et $\theta$ y agit par l'inversion, 
donc $S^\theta=1$ et $P(1)=2$. De plus, $\theta$ pr\'eserve le sous-groupe de Borel triangulaire sup\'erieur, son tore diagonal et 
agit sur ce dernier par $(x_1,x_2,\cdots,x_{2n}) \mapsto (x_{2n}^{-1},\cdots,x_2^{-1},x_1^{-1})$ donc $\varepsilon(\theta)=(-1)^{n-1}$ (et $q(G)=2n-1$).}
\end{exemple}

\subsubsection{Nombre de Lefschetz d'un automorphisme} 
Suivant Borel-Wallach \cite[X.5]{BW}, notons $H^i_e(G(F),-)$ les foncteurs d\'eriv\'es du foncteur des 
$G(F)$-invariants de la cat\'egorie 
des repr\'esentations complexes lisses de $G(F)$ dans celle des espaces vectoriels complexes. 
Il est d\'emontr\'e {\it loc. cit.} que si $V$ est une repr\'esentation lisse admissible de $G(F)$, 
alors les $H^i_e(G(F),V)$ sont de dimension finie, nuls pour $i>q(G)+\dim_F(S)$. Si $V$ est de plus la restriction \`a $G(F)$ d'une repr\'esentation de $G^+(F)$, 
ces espaces sont munis d'une action naturelle de $G^+(F)/G(F)=\langle \theta \rangle$, ce qui d\'efinit 
le \og nombre de Lefschetz\fg\,  de $\theta$:

$$\Lef_G(\theta,V)=\sum_{i\geq 0} (-1)^i \trace(\theta,H^i_e(G(F),V)).$$ 

Ces nombres ont \'et\'e notamment \'etudi\'es par Borel-Labesse-Schwermer \cite[\S 9]{BLS}\footnote{Notons que la Proposition 8.2 (2) de \cite{BLS}
est incorrecte: si $S^{\theta} \neq 1$ alors $\Lef_G(\theta,-)$ est identiquement nul, ce qu'ils d\'emontrent 
d'ailleurs au passage au cours de la preuve de leur Proposition 8.4. 
L'erreur se trouve dans l'appel \`a \cite[X.4.7]{BW}, qui n\'ecessite $G$ simplement connexe.}. La
proposition suivante repose essentiellement sur des r\'esultats de Casselman et Borel-Wallach.

\begin{prop} $\Lef(\theta,\cdot)$ est non identiquement nul si, et seulement si, $S^{\theta}=\{1\}$, ou ce qui est \'equivalent si $P(1) \neq 0$. 
Supposons donc que $P(1)\neq 0$ et fixons $V$ une repr\'esentation irr\'eductible unitaire de $G^+(F)$.
\begin{itemize} 
\item[(i)]  $\Lef_G(\theta,1)=P(1)$ et $\Lef_G(\theta,\St)=P(1)(-1)^{q(G)}\varepsilon(\theta)$.
\item[(ii)] Si $V$ est essentiellement temp\'er\'ee, $\Lef_G(\theta,V) \neq 0$ si, et seulement si, $V=\St$.
\item[(iii)] Si $G/S$ est quasi-simple, $\Lef_G(\theta,V) \neq 0$ si, et seulement si, $V=\St$ ou $V=1$.
\end{itemize}
\end{prop}

\begin{pf} Soit $V$ une repr\'esentation complexe lisse de $G(F)$. D'apr\`es Hitta \cite{Hi}, on dispose pour la 
paire $S(F) \subset G(F)$ d'une suite spectrale de type Hochschild-Serre 
$$E_2^{p,q}=H^p_e(G(F)/S(F),H^q_e(S(F),V)) \Longrightarrow H^{p+q}_e(G(F),V).$$ 
Si $V$ est la restriction \`a $G(F)$ d'une repr\'esentation de $G^+(F)$, cette suite spectrale est par construction 
naturellement munie d'une action de $G^+(F)/G(F)$ compatible \`a l'action rappel\'ee plus haut sur son aboutissement. 
Si de plus $S(F)$ agit trivialement sur $V$, $E_2^{p,q}$ s'\'ecrit alors 
$H^q_e(S(F),1)\otimes_{\C}H^p_e(G(F)/S(F),V)$. Dans le cas o\`u $V$ est de surcro\^{i}t admissible 
on obtient donc l'identit\'e
\begin{equation}\label{identite}\Lef_G(\theta,V)=\Lef_S(\theta,1)\Lef_{G/S}(\theta,V).\end{equation}
\par
Supposons maintenant que $W$ est une repr\'esentation lisse irr\'eductible de $G(F)$ telle que $H^i_e(G(F),W) \neq 0$ pour un certain entier $i\geq 0$.
D'apr\`es un r\'esultat de Borel-Casselman, $W$ est 
un constituant irr\'eductible de $I_B$ (\cite[X.2.4, X.4.3]{BW}). En particulier, $W$ est de caract\`ere central trivial et elle se factorise donc en 
une repr\'esentation de $G(F)/S(F)=(G/S)(F)$. Ainsi, d'apr\`es l'identit\'e (\ref{identite}), il suffit de d\'emontrer la proposition dans les cas $S=\{1\}$ et 
$G=S$ (et dans ce cas pour $V=1$).

Supposons d'abord $G=S$. La valuation $F^* \longrightarrow \Z$ induit un isomorphisme canonique 
$X^*(S) \otimes_{\Z} \C \isomo \Hom_{\Z}(S(F)/S(F)^0,\C)$, o\`u $S(F)^0$ est le sous-groupe compact maximal de $S(F)$. Le calcul standard de la cohomologie des tores nous fournit 
alors des isomorphismes canoniques (\cite[X.2.6]{BW})
\begin{equation}\label{cohotore} H^i_e(S(F),1) \isomo \Lambda^i(X^*(S)\otimes_{\Z} \C), \, \, \, i\geq 0. \end{equation}
d'o\`u l'on tire $\Lef_S(\theta,1)=P(1)$.

Supposons enfin $S=\{1\}$ (i.e. $G(F)$ de centre compact) et fixons $V$ comme dans l'\'enonc\'e. 
D'apr\`es Borel-Wallach \cite[XI.3.8]{BW}, si $V$ est temp\'er\'ee et 
si $H^i_e(G(F),V) \neq 0$, alors $i=q$, $V=\St$. De plus, si $H^i_e(G(F),1)\neq 0$ alors $i=0$. Enfin, si $G$ est quasi-simple et si $H^i_e(G(F),V) \neq 0$, 
alors $V=\St$ ou $V=1$ d'apr\`es \cite[XI.3.9]{BW}. Il ne reste donc qu'\`a verifier que $\theta$ agit sur $H^{d(G)}_e(G(F),\St)$ 
par la multiplication par le signe $\epsilon(\theta)$ d\'efini au \S \ref{notationslef}. Comme les $H^i_e(G(F),\cdot)$ s'annulent si $i> q(G)$, on dispose d'une surjection 
$$H^{q(G)}_e(G(F),I_B) \longrightarrow H^{q(G)}_e(G(F),\St)$$ et il suffit de voir que l'espace de gauche est de dimension $1$ 
et de montrer que l'action de $\theta$, ou ce qui est \'equivalent de $\theta_B$, y agit par $\varepsilon(\theta)$.
Par le lemme de Shapiro (\cite[X.4.2]{BW}) appliqu\'e \`a la repr\'esentation induite $I_B$, on dispose d'une identification $\theta_B$-\'equivariante 
$H^{q(G)}_e(G(F),I_B) \isomo H^{q(G)}_e(M(F),1)$ o\`u $M$ est le Levi de $B$ centralisant $S'$. 
Comme $M(F)/S'(F)=(M/S')(F)$ est compact, 
la suite spectrale de Horschild-Serre rappel\'ee plus haut fournit pour tout $i\geq 0$ un isomorphisme canonique 
$\theta_B$-equivariant $H^i_e(S'(F),1) \isomo H^i_e(M(F),1)$, 
et on conclut alors par l'identit\'e (\ref{cohotore}) appliqu\'e \`a $G=S'(F)$ et $i=q(G)$. 
\end{pf}

\begin{remarque}\label{remcuspidal}{\rm  Soient $P$ un sous-groupe parabolique de $G$ et $M$ un sous-groupe de Levi de $P$ qui sont tous
deux d\'efinis sur $F$ et normalis\'es par un \'el\'ement $\theta_M \in
G(F)\theta$. Par le lemme de Shapiro lisse, pour toute repr\'esentation admissible $W$ de $M(F)
\rtimes \langle \theta_M \rangle$, on a $$\Lef_G(\theta,{\rm Ind}_P^G
W)=\Lef_M(\theta_M,W\delta_P^{1/2})$$
(l'induite de $P(F)$ \`a $G(F)$ \'etant lisse et normalis\'ee). En
particulier, ce nombre est nul si $W$ est unitaire et $P \neq G$, car les constituants
irr\'eductibles de $W\delta_P^{1/2}$ sont de caract\`ere central non
trivial (et m\^eme non unitaire). }
\end{remarque}

\subsubsection{Fonctions d'Euler-Poincar\'e d'un automorphisme}

Fixons une mesure de Haar $\mu$ sur $G^+(F)$. 

\begin{prop}\label{nbleftheta} Supposons $P(1)\neq 0$. Il existe une function $f_{EP}: G^+(F) \longrightarrow \C$ localement constante et \`a support compact dans 
$\theta G(F)$ telle que pour toute representation admissible $(\pi,V)$ de $G^+(F)$ on ait 
$\trace(\pi(f_{EP}))=\Lef_G(\theta,V)$. Elle jouit de plus des propri\'et\'es suivantes:
\begin{itemize}
\item[(i)] La fonction $f'_{EP}:=\frac{\epsilon(\theta)(-1)^q(G)}{P(1)}f_{EP}$ est un pseudo-coefficient de $\St$: 
si $(\pi,V)$ est irr\'eductible et essentiellement temp\'er\'ee, alors $\trace(\pi(f'_{EP}))=1$ si $V=\St$, $0$ sinon.
\item[(ii)] Soient $\gamma \in \theta G(F)$ un \'el\'ement semisimple, $I_\gamma$ la composante neutre
du centralisateur de $\gamma$ dans $G$ ; on munit $I_\gamma(F)\backslash G(F)$
d'une mesure $G(F)$-invariante $\overline{\mu}$. Alors l'int\'egrale orbitale \og tordue\fg\, 
		$$O_{\gamma}(f_{EP}):=\int_{I_\gamma(F)\backslash G(F)}  
f_{EP}(g^{-1}\gamma g) \overline{\mu}$$ est non nulle si et seulement si $I_\gamma(F)$ est de centre compact. 
\end{itemize}
Dans le cas o\`u $G/S$ est quasi-simple, (i) admet la variante plus forte\begin{itemize}
\item[(i)'] Si $G/S$ est quasi-simple et $V\neq \St, 1$ est irr\'eductible et unitaire, alors $\trace(\pi(f'_{EP}))=0$.
\end{itemize}
\end{prop}

\begin{remarque} {\rm Ces fonctions $f_{EP}$ ont \'et\'e 
introduites par Kottwitz dans \cite[\S 2]{Ko} sous le nom de fonctions d'Euler-Poincar\'e dans le cas o\`u $S$ et $\theta$ sont triviaux, il y d\'emontre la proposition dans ce cadre. 
Certaines de leurs propri\'et\'es ont \'et\'e \'etendues par Borel-Labesse-Schwermer dans \cite{BLS} 
dans la g\'en\'eralit\'e adopt\'ee ici (voir cependant les notes de bas de page qui suivent). Les arguments 
de cette section ne sont que des adaptations essentiellement triviales de ces r\'esultats. Notons que l'existence d'une fonction 
$f_{EP}$ satisfaisant (i) pourrait se d\'eduire du th\'eor\`eme de Paley-Wiener, mais la d\'emonstration de Kottwitz 
a l'avantage de fournir une fonction $f_{EP}$ explicite en terme de l'immeuble de Bruhat-Tits de $G$ et de fournir la propri\'et\'e 
(ii) \`a peu de frais.}
\end{remarque}

Le reste de cette partie sera consacr\'e \`a la preuve de la proposition. 

\subsubsection{L'immeuble de Bruhat-Tits de $G$ et d\'efinition de $f_{EP}$} 
\newcommand{\B}{{\cal B}}\label{BruhatTits}
Soit $\B$ l'immeuble de Bruhat-Tits de $G$ (\cite{Ti}). On rappelle que c'est un complexe polysimplicial (\cite[\S 1.1]{BT}) 
muni d'une action cellulaire de $G(F)$, et m\^eme de $G^+(F)$. 
Pour fixer les id\'ees, disons que nous consid\'erons la normalisation canonique de sa partie torale au sens de Rousseau \cite[\S 1.2]{Ti}. 
Rappelons qu'un {\it sommet} de $\B$ est un polysimplexe de dimension $0$ et qu'une {\it chambre} de $\B$ est l'int\'erieur 
d'un polysimplexe de dimension maximale, en l'occurence 
$\dim(S)+q(G)$. L'immeuble $\B$ est contractile et a la propri\'et\'e que le stabilisateur dans $G(F)$ 
de chaque partie compacte de $\B$ est un sous-groupe compact ouvert, de sorte que le complexe de ses chaines 
polysimpliciales orient\'ees permet de calculer la cohomologie des repr\'esentations lisses de $G(F)$ (Casselman-Wigner, Borel-Wallach \cite[X.2]{BW}) de la mani\`ere suivante. 
\renewcommand{\P}{{\cal P}}
Rappelons qu'un polysimplexe est la donn\'ee d'un couple $(p,\epsilon)$ 
o\`u $p$ est un polysimplexe de dimension $>0$ et $\epsilon$ une orientation de $p$. 
Pour $i\geq 0$ notons $\P_{i,\pm}$ l'ensemble des polysimplexes orient\'es de $\B$ de dimension $i$, il est muni 
d'une action naturelle de $G^+(F)$. Si $V$ est une repr\'esentation lisse de $G(F)$ et $i\geq 0$, d\'efinissons $C_i(V)$ comme \'etant 
l'espace des fonctions $G(F)$-\'equivariantes $\varphi: \B_{i,\pm} \longrightarrow V$ telles que 
$\varphi((p,-\varepsilon))=-\varphi(p,\varepsilon)$ pour tout polysimplexe orient\'e $(p,\varepsilon) \in \B_{i,\pm}$ 
de dimension $>0$. Les $C:=(C_i)_{i\geq 0}$ forment alors un complexe de cochaines 
dont la cohomologie est $H^i_e(G(F),V)$. On suppose dor\'enavant que $V$ est une repr\'esentation admissible 
de $G^+(F)$, $G^+(F)/G(F)$ agit alors naturellement sur $C$ par la formule 
$\theta(f)(x)=\theta(f(\theta^{-1}(x))), \forall f\in C_i(V)$. 
\newcommand{\sign}{{\rm sign}}

D'apr\`es Bruhat-Tits, les chambres de $\B$ sont permut\'ees transitivement par $G(F)$, de sorte que 
l'action de $G(F)$ sur les polysimplexes n'a qu'un ensemble fini $\Sigma$ d'orbites. Si $s$ est une 
telle orbite, on note $\dim(s)$ la dimension commune de ses polysimplexes, et on note $\Sigma_i \subset \Sigma$ 
le sous-ensemble des orbites de dimension $i$. Il vient que l'on a une d\'ecomposition $G(F)$-\'equivariante $C_i(V)=\oplus_{s \in \Sigma_i} C_i(V)_s$, o\`u 
$C_i(V)_s$ est le sous-espace des fonctions \`a support dans l'orbite $s$. 

Si $p$ est un polysimplexe, notons $G_p$ (resp. $G^+_p$) le sous-groupe de $G(F)$ (resp. $G^+(F)$) laissant $p$ globalement invariant : 
c'est un sous-groupe compact ouvert. On dispose alors d'un caract\`ere d'orientation\footnote{Si $p$ est un 
simplexe et $g \in G_p^+$, $\sign_p(g)$ est aussi la signature de la permutation des sommets de $p$ induite 
par $g$. Contrairement \`a ce qui semble sous-entendu dans \cite[\S 2]{Ko} et dans la preuve de 
\cite[Prop. 8.4]{BLS}, pr\'ecisons que cela ne vaut plus pour un polysimplexe g\'en\'eral: consid\'erer une r\'eflexion d'un 
carr\'e d'axe passant par le milieu d'un cot\'e (cela se produit par exemple si $G=G^+=PGL_2 \times PGL_2$).} 
continu $\sign_p: G^+_p \rightarrow \{\pm 1\}$ et on note $V_p$ le plus grand sous-espace de $V$ sur lequel 
$G_p$ agit par $\sign_p$. Si $p \in s$, $\phi \mapsto \phi(p)$ induit alors un isomorphisme $C_i(V)_s \isomo V_p$ 
(qui est donc de dimension finie). Il y a deux cas: \begin{itemize}
\item[-] Si $\theta(s)\neq s$, alors $\trace(\theta,\sum_{n \in \Z} C_i(V)_{\theta^n(s)})=0$.
\item[-] Si $\theta(s)=s$, il existe $\theta_s \in G(F)\theta$ tel que $\theta_s(p)=p$. La trace de 
$\theta$ sur $C_i(V)_s$ ne d\'ependant que de l'image de $\theta$ dans $G^+(F)/G(F)$, on peut la calculer pour 
$\theta_s$ par \'evaluation en $p$, et l'on trouve $\sign_p(\theta_s)\trace(\theta_s,V_p)$ 
(noter que $G^+_p=\langle G_p,\theta_s \rangle$ et $\theta_s$ normalise $G_p$). 
\end{itemize} 

Pour chaque $s \in \Sigma$ fix\'e par $\theta$, on choisit alors un $p_s \in s$ et on d\'efinit $f_s$ comme 
\'etant la fonction sur $G^+(F)$ qui vaut $\sign_p$ sur $G_p^+ \cap \theta G(F)$ et qui est nulle ailleurs. 
Si l'on pose\footnote{Noter que strictement, l'expression de $f_{EP}$ donn\'ee par Borel-Labesse-Schwermer 
dans la preuve de \cite[Prop. 8.4]{BLS} n'est valable que quand $G$ est semisimple et simplement connexe, qui est 
l'hypoth\`ese de validit\'e de \cite[X.2.5]{BW} auquel ils renvoient, et auquel cas le complexe $C$ se simplifie. 
C'est une des raisons pour laquelle nous avons redonn\'e l'argument complet.}
\begin{equation}\label{fEP} f_{EP}:=\sum_{s \in \Sigma, \theta(s)=s}
(-1)^{\dim(s)}\mu(G_{p_s})^{-1}f_s,\end{equation}
on a montr\'e que $\trace(\pi(f_{EP}),V)=\Lef_G(\theta,V)$, {\it i.e.} que $f_{EP}$ est une fonction 
d'Euler-Poincar\'e pour $\theta$. Le point (i) de la proposition d\'ecoule alors de la Proposition \ref{nbleftheta}. 
Notons que par construction, la mesure $\mu f_{EP}$ est ind\'ependante du choix de $\mu$, et que les int\'egrales orbitales des 
$f_s$, et donc de $f_{EP}$, ne d\'ependent pas du choix des $p_s$. 
                                                      
\subsubsection{Int\'egrales orbitales de $f_{EP}$ en les \'el\'ements
semisimples
de $G(F)\theta$.} Il ne reste qu'\`a d\'emontrer (ii). La d\'emonstration de Kottwitz \cite[Theorem 2.2]{Ko} s'\'etend
essentiellement verbatim. Fixons agir $G(F)$ sur $G(F)^+$ par la conjugaison \og \`a droite\fg\, 
$g.\gamma:=g^{-1}\gamma g$
et fixons $\gamma \in G(F)\theta$ un \'el\'ement de $G(F)^+$
comme dans l'\'enonce. Comme $\gamma$ est semisimple, l'orbite $\omega(\gamma):=G(F).\gamma$ est alors
ferm\'ee dans $G(F)\theta$ et l'integrale orbitale de l'\'enonc\'e est trivialement convergente.
\newcommand{\FF}{\cal F}

Le sous-espace $\B^\gamma \subset \B$ des points fixes de $\gamma$
est muni d'une structure de complexe polysimplicial
dont les polysimplexes sont les $p \cap \B^{\gamma}$ non vides, {\it i.e.}
tels que $p$ est dans l'ensemble $\FF(\gamma)$ des polysimplexes
de $\B$ stables globalement par $\gamma$. Supposons que $\B^\gamma$ est non vide, ce qui est toujours le cas si le centre de $I_\gamma(F)$ est compact\footnote{En effet, une puissance finie de $\gamma$ est dans le centre de $I_\gamma(F)$, elle admet donc un point fixe 
$x \in B^\gamma$ quand ce centre est compact. L'ensemble fini $\{\gamma^nx, n \in \Z\} \subset \B$ est alors stable par $\gamma$ et son enveloppe convexe contient un \'el\'ement de $\B^\gamma$ d'apr\`es \cite[Prop. 3.2.4]{BT}.}. C'est alors un
ensemble clos au sens de Bruhat-Tits, il est en particulier contractile. Il 
h\'erite de $(\B,G(F))$ les propri\'et\'es suivantes: $I_\gamma(F)$ agit de mani\`ere
cellulaire sur $\B^\gamma$, tout point de $\B^\gamma$ a un stabilisateur
compact ouvert, et tout sous-groupe compact de $I_{\gamma}$ fixe un point de
$\B^\gamma$. Enfin, si $p$ est un polysimplexe de $\B$ et $g \in G(F)$,
\begin{equation} \label{lien} g.\gamma \in G_p^+ \Leftrightarrow g p \in
\FF(\gamma),\end{equation}
ce qui ne d\'epend que de $I_\gamma(F)gG_p$. En
particulier, $I_{\gamma}(F)$ n'a qu'un nombre fini d'orbites sur $G(F)p \cap
\FF(\gamma)$, et donc sur $\FF(\gamma)$, car $\omega(\gamma)\cap G_p^+$ est
compact. D'apr\`es Serre \cite[\S 3.3]{Serre}, toutes les
conditions sont donc satisfaites pour que le complexe cellulaire $B^\gamma$
permette de calculer la mesure d'Euler-Poincar\'e de $I_\gamma$: si
$\mu_I$ est la mesure de Haar sur $I_\gamma$ telle que
$\overline{\mu}=d\mu/d\mu_I$, c'est la
mesure sign\'ee canonique
        $$\mu_{EP,I_\gamma}:=\sum_{p \in I_\gamma(F)\backslash \FF(\gamma)}
(-1)^{\dim(p)}\frac{\mu_I}{\mu_I(I_\gamma(F)\cap G_p)}.$$
D'apr\`es \cite[Prop. 28]{Serre}, cette mesure est non nulle si, et seulement si
$I_\gamma(F)$ est de centre compact, son signe est alors
$(-1)^{q(I_{\gamma})}$.

Il ne reste qu'\`a v\'erifier que  $O_\gamma(f_{EP})\mu_I=\mu_{EP,I_\gamma}$ si $\B^\gamma\neq \emptyset$ et qu'elle est nulle sinon. Notons d'abord que si $\theta(s)\neq s$, $s \cap \FF(\gamma)=\emptyset$.
Soit $s=G(F)p_s \in \Sigma^\theta$, alors par d\'efinition de $f_s$ et (\ref{lien})
on a
        $${\rm O}_\gamma(f_s)=\sum_{p \in I_\gamma(F)\backslash (s\cap
\FF(\gamma))}
\sign_p(\gamma)\mu_I(I_\gamma(F)\cap G_p)^{-1}.$$
On conclut en notant que\footnote{Se ramener \`a v\'erifier cette relation pour le d\'eterminant d'une isom\'etrie
d'un espace affine euclidien ayant un point fixe.} si $p \in \FF(\gamma)$,
$\sign_p(\gamma)=(-1)^{\dim(p)-\dim(p\cap
\B^{\gamma})}$. 
\par \medskip

\subsubsection{Int\'egrales orbitales non semisimples} Un \'el\'ement semisimple $\gamma \in \theta G(F)$ est dit {\it elliptique} si $I_\gamma(F)$ est de centre compact. 

\begin{prop}\label{connoncon} Soit $f\in \Cc^\infty(\theta G(F))$ telle que $\ON_\gamma(f)=0$ pour tout \'el\'ement $\gamma \in \theta G(F)$ semisimple non elliptique, alors $\ON_\gamma(f)=0$ pour tout \'el\'ement $\gamma \in \theta G(F)$ non semisimple.
\end{prop}

\begin{pf} En effet, \'ecrivons
$$\gamma = \nu \sigma = \sigma \nu$$
la d\'ecomposition de Jordan de $\gamma$, o\`u $\nu \in G(F)$ est unipotent et $\sigma \in \theta G(F)$ semisimple. Si $I=I_\gamma$ et $M=I_\sigma$ d\'esignent les centralisateurs connexes respectifs de $\gamma$ et $\sigma$ dans $G$, alors $M$ est r\'eductif, $I$ est unimodulaire et $I \subset M$. De plus, $\nu \in M(F)$. Si $f \in {\cal C}_c^\infty(\theta G(F))$, on peut \'ecrire 

$$\ON_\gamma(f)\,=\, \int_{I(F)\backslash G(F)} f( g^{-1} \gamma g )\frac{dg}{di}\,\, $$ $$\qquad \, \,= \, \, \int_{M(F)\backslash G(F)}\int_{I(F)\backslash M(F)} f(g^{-1}m^{-1} \nu \sigma m g ) \frac{dm}{di} \frac{dg}{dm}$$
$$ \qquad = \, \, \int_{M(F)\backslash G(F)}\int_{I(F)\backslash M(F)} f(g^{-1}m^{-1} \nu m \sigma g ) \frac{dm}{di} \frac{dg}{dm},$$
la seconde int\'egrale \'etant l'int\'egrale orbitale {\it ordinaire} dans le groupe r\'eductif connexe $M$ de la fonction $$f_g^M(m):=f(g^{-1}m\sigma g)$$
en $\nu \in M(F)$ (remarquer que le centralisateur de $\nu$ dans $M$ est exactement $I$).\ps

	D'apr\`es le lemme de compacit\'e usuel, \'etendu par Arthur~\cite[Lemme 2.1]{ArD} au 
cas tordu, il existe un voisinage ouvert $M(F)$-invariant $\cal U$ de $1$ dans $M(F)$ tel que, 
pour tout sous-ensemble compact $\Omega$ de $\theta G(F)$, il existe un sous-ensemble compact 
$\omega$ de $M(F)\backslash G(F)$ ayant la propri\'et\'e suivante : si $g \in G(F)$ et 
$g^{-1} {\cal U} \sigma g \cap \Omega \neq \emptyset$, alors $M(F)g \in \omega$. Soit $M^+$ le centralisateur de $\sigma$ dans $G$. \ps

\begin{lemme}\label{descenteconjugaison} Il existe un voisinage ouvert compact $\cal V$ de $1$ dans $M(F)$ ayant la propri\'et\'e suivante : pour tous $u\in \cal V$ et $g \in G(F)$ tels que $g^{-1}u \sigma g \in \sigma \cal V$, alors $g \in M^+(F)$.
\end{lemme}

\begin{pf} Il s'agit d'une variante de \cite[Lemme 3.1.4]{labesse}. D'apr\`es \cite[Lemme 3.1.1]{labesse}\footnote{Bien que le cadre adopt\'e {\it loc. cit.} soit celui du changement de base, la d\'emonstration de Labesse s'applique verbatim au cas g\'en\'eral.}, il existe une sous-vari\'et\'e analytique $Y=Y^{-1} \subset G(F)$ et un voisinage ouvert $\cal V$ de $1$ dans $M(F)$ telle que l'application $$ Y \times \cal V \longrightarrow \sigma G(F), \, \, \, \, (y,u) \mapsto y^{-1}u\sigma y,\, \, $$
soit un diff\'eomorphisme sur son image $\Omega$, un voisinage ouvert compact de $\sigma$ dans $\sigma G(F)$. De plus, si $(y,u) \in Y \times M(F)$ est tel que $y^{-1}u\sigma y \in \sigma \cal V $, alors $u \in \cal V$ et $y=1$. En particulier, 
\begin{equation}\label{preuvelabesse} \forall\, \, (g,u)\, \, 
\in \, \, (M^+(F)Y) \times \cal V, \, \, \, g^{-1}u\sigma g \in \sigma \cal V \Longrightarrow g \in M^+(F).
\end{equation}

Supposons maintenant par l'absurde que l'\'enonc\'e ne tient pas : il existe une suite d'\'el\'ements $(g_n,u_n,v_n)$ de $(G(F)\backslash M^+(F))\times \cal V \times \cal V$ telle que $(u_n,v_n) \rightarrow (1,1)$ et $g_n^{-1}u_n\sigma g_n=v_n \sigma$. D'apr\`es (\ref{preuvelabesse}), $g_n \notin M^+(F)Y$. D'apr\`es le lemme de compacit\'e, on peut trouver un compact $\omega \subset G(F)$ tel que $\omega \cap M^+(F) =\emptyset$ et tel que pour tout $n$ assez grand, $g_n=m_nw_n$ avec $w_n \in \omega$ et $m_n \in M^+(F)$. On peut donc supposer que $w_n \rightarrow w  \in \omega$, auquel cas $m_n^{-1}u_nm_n$ converge, disons vers $u^* \in M(F)$, qui est n\'ecessairement unipotent car $u_n \rightarrow 1$. De plus, $$w^{-1} u^* \sigma w = \sigma,$$
de sorte que $u^*=1$ par unicit\'e de la d\'ecomposition de Jordan de $\sigma$, puis $w \in M^+(F)$, ce qui est absurde.
\end{pf}

	Si $u \in M(F)$, consid\'erons le centralisateur connexe $I_{u\sigma}$ de $u\sigma$ 
dans $G$. D'apr\`es le Lemme~\ref{descenteconjugaison}, on peut trouver un voisinage ouvert $\cal V$ de $1$ dans $M(F)$ tel que pour tout 
$u \in \cal V$, $I_{u\sigma} \subset M$. Quitte \`a remplacer $\cal V$ par 
$(\cup_{m \in M(F)} m^{-1} \cal V m) \cap \cal U$ on peut supposer que $\cal V$ est 
$M(F)$-invariant et inclus dans $\cal U$. Noter que $\cal V$ (tout comme $\cal U$) contient alors 
tous les \'el\'ements unipotents de $M(F)$.\ps
	
	Soit $\Omega$ le support de $f$, le lemme de compacit\'e lui associe un $\omega$ que 
l'on peut supposer ouvert compact. Soient $\chi \in \Cc^{\infty}(M(F)\backslash G(F))$ la 
fonction caract\'eristique de $\omega$, $\alpha \in \Cc^\infty(G(F))$ telle que 
$\int_{M(F)} \alpha(m\,g) dm = \chi(M(F)\,g)$, et soit $$h(m)=\int_{G(F)} \alpha(g) f_g^M(m) dg.$$
Il est clair que $h$ est localement constante sur $M(F)$ ; son support est compact car si 
$h(m)\neq 0$, $m\sigma \in \Supp(\alpha)\Supp(f)\Supp(\alpha)^{-1}$. Si $u \in \cal V$, 
alors $I_{u\sigma} \subset M$ et on montre comme plus haut que 
$$\ON_{u\sigma}(f)=\int_{M(F)\backslash G(F)}\ON_{u}^M(f_g^M) \frac{dg}{dm} = \int_{M(F)\backslash G(F)} \chi(gM) \ON_{u}^M(f_g^M) \frac{dg}{dm}, $$
cette derni\`ere \'egalit\'e venant de ce que $\cal V \subset \cal U$, puis 
\begin{equation}\label{descentecentrcon}\ON_{u\sigma}(f)\, =\, \int_{G(F)} \alpha(g) \ON_{u}^M(f_g^M) \,dg = \ON_u^M(h).\end{equation}
Par cons\'equent, pour $u\in \cal V$, $\ON_{u\sigma}(f)$ est une int\'egrale orbitale {\it ordinaire} sur le groupe $M$.

Si $M$ contient un tore maximal non compact $T$, alors  $$\ON^M_u(h)=\ON_{u\sigma}(f)=0$$ pour tout \'el\'ement r\'egulier 
$u\sigma$ dans $(T(F)\cap \cal V)\sigma$ d'apr\`es (\ref{descentecentrcon}) et par hypoth\`ese (ces $u\sigma$ sont 
semisimples non elliptiques). D'apr\`es un r\'esultat de Rogawski (voir \cite[p. 636]{Ko}), on en d\'eduit que les int\'egrales orbitales unipotentes de $h$ s'annulent, et en particulier que $$\ON_\nu^M(h)=\ON_{\gamma}(f)=0.$$ Si enfin $M$ est anisotrope, alors $\nu=1$, $\gamma$ est semisimple, et il n'y a rien \`a d\'emontrer. \end{pf}

\begin{remarque}\label{remdescentecentrcon}{\rm  Ce r\'esultat est suppos\'e \og bien connu\fg , mais n'est apparement d\'emontr\'e nulle part.}
\end{remarque}

La proposition~\ref{nbleftheta} (ii) admet le corollaire suivant:

\begin{corollaire}\label{corTOnonss} Soit $\gamma \in \theta G(F)$ non semisimple elliptique, alors ${\rm O}_\gamma(f_{EP})=0$.
\end{corollaire}

\subsection{Preuve du Th\'eor\`eme \ref{mainthm}} \label{preuvemainthm}

Repla\c{c}ons nous sous les hypoth\`eses du \S \ref{ssectenonce}.

\subsubsection{Une version simplifi\'ee de la formule des traces d'Arthur} Soit $A=\R_{+}^{\ast}$ 
la composante neutre topologique du centre de 
$\GL_{2n}(\R)$, \'equipons l'espace homog\`ene\,
$A \, \GL_{2n}(\Q)\backslash \GL_{2n}(\AAA)$ d'une mesure (finie) $G(\AAA)$-invariante \`a droite. La repr\'esentation unitaire 
$R$ de $G(\AAA)$ par translations \`a droite sur l'espace des fonctions cuspidales 
	$$L^2_{\rm cusp}(A\,\GL_{2n}(\Q)\backslash \GL_{2n}(\AAA))$$
s'\'etend en une repr\'esentation unitaire de $G^+(\AAA)$ si l'on fait agir $\theta$ par 
l'op\'erateur $I_{\theta}(\varphi)(x)=\varphi(\theta(x))$. Cette repr\'esentation est discr\`ete. Si $f=f_{\infty}\otimes 
f^{\infty}$ est dans $\Cc^{\infty}(G(\AAA))$, avec $f_{\infty}$ disons 
${\rm SO_{2n}}(\R)$-finie, alors $R(f)I_{\theta}$ est tra\c{c}able et
$$\trace(R(f)I_{\theta})=\sum_{\Pi} \trace(R(f)I_{\theta},\Pi),$$
la somme ci-dessus portant sur les repr\'esentations automorphes cuspidales irr\'eductibles et autoduales de $G(\AAA)$. 
Ces traces d\'ependent toutes d'un choix de mesure de Haar ad\'elique $dg_\AAA$ sur $G(\AAA)$ que nous fixons une fois pour
toutes.

\par \medskip
Pour un choix de fonctions tests $f$ convenables les r\'esultats d'Arthur
donnent une expression g\'eom\'etrique simple de $\trace(R(f)I_{\theta})$. \'Ecrivons
pour cela $f=f_\infty \otimes f_\ell \otimes f_p \otimes f^{\infty,\ell,p}$ o\`u:\begin{itemize}
\item[-] $f_\infty$ est un pseudo-coefficient tordu d'une s\'erie
$\theta$-discr\`ete cohomologique (cf. \S\ref{prelimglnreel},
\S\ref{pseudoorb}). 
\item[-] $f_{\ell}$ est une fonction d'Euler-Poincar\'e fix\'ee pour l'automorphisme $\theta$ de $G_{\Q_l}$ donn\'ee par la Proposition \ref{nbleftheta} (cf. Exemple \ref{exgl2n}).
\item[-] $f_p$ est une fonction donn\'ee par le Th\'eor\`eme \ref{thmcle},
\item[-] $f^{\infty,\ell,p}$ est la fonction caract\'eristique de $\prod_{v\notin \{\infty,\ell,p\}} \GL_{2n}(\Z_v)$.
\end{itemize}
\ps
Rappelons qu'un \'el\'ement $\theta$-semisimple $\gamma \in \G(\Q)$ est
dit elliptique si la composante d\'eploy\'ee du centre de son centralisateur tordu
est triviale. Notons $\{G(\Q)\}_{\rm ell}$ l'ensemble des classes de $\theta$-conjugaison
d'\'el\'ements $\theta$-semisimples elliptiques.
Pour $\gamma \in G(\Q)$ un tel \'el\'ement, on choisit une mesure de Haar ad\'elique $di_\AAA$ sur $I_\gamma(\AAA)$ et on pose $v_\gamma=\mu(I_\gamma(\Q)\backslash I_\gamma(\AAA))>0$ et $$\TO_\gamma(f):=\int_{I_\gamma(\AAA)\backslash G(\AAA)} f(g^{-1}\gamma
\theta(g))\, \, di_\AAA\backslash dg_\AAA.$$ 
On consid\'erera aussi les version locales \'evidentes de ces int\'egrales orbitales tordues.

\begin{prop} \label{formtrace} Pour toute fonction $f$ comme plus haut, 
$$\trace(R(f)I_\theta)=\sum_{\gamma \in \{G(\Q)\}_{\rm ell}} v_\gamma \, \TO_\gamma(f).$$ 
La somme porte sur un sous-ensemble fini de classes qui ne d\'epend que d'un
compact de $G(\AAA)$ contenant le support de $f$.
\end{prop}

\begin{pf} Nous allons appliquer les r\'esultats d'Arthur \cite{ITF} \`a la
composante connexe $\G\theta$. Ainsi que l'explique Arthur (\cite[p.
330]{ITF1},\cite[p. 528]{ITF}), la validit\'e de ces r\'esultats dans ce
cadre d\'epend de la v\'erification d'un argument de cohomologie galoisienne et 
de la validit\'e du th\'eor\`eme de
Paley-Wiener pour $G(\R)\theta$. Le premier a en fait \'et\'e v\'erifi\'e
depuis en toute g\'en\'eralit\'e par Kottwitz et Rogawski \cite{KoR}, ainsi
que le second dans notre cadre par Mezo \cite{Mezo}. 

Les propri\'et\'es de pseudocoefficients des fonctions $f_\infty$ et
$f_\ell$ (\S\ref{pseudoorb}, Rem.~\ref{remcuspidal}) impliquent que la fonction $f$ est cuspidale
au sens d'Arthur \cite[\S 7 p. 538]{ITF} en les deux places $\infty$ et
$\ell$. Mieux, en la place $\ell$, les int\'egrales 
orbitales de $f_\ell$ en les \'el\'ements non semisimples $\Q_l$-elliptiques
s'annulent par le Corollaire~\ref{corTOnonss}, de sorte que \cite[Cor. 7.5]{ITF} s'applique. Ce corollaire identifie le terme de droite de l'\'enonc\'e \`a
la trace de $fI_\theta$ dans la repr\'esentation de $G^+(\AAA)$ sur 
$$L^2_{\rm disc}(A_G\GL_{2n}(\Q)\backslash
\GL_{2n}(\AAA))$$
Notons qu'Arthur consid\`ere {\it loc. cit.} une sommation pr\'ecise pour cette trace
partitionn\'ee par les normes $t$ possibles des caract\`eres
infinit\'esimaux des $\Pi_\infty$. Comme $f_\infty$ ne trace que dans 
des repr\'esentations ayant m\^eme caract\`ere infinit\'esimal, un seul de
ces $t$ intervient. Il ne reste qu'\`a voir que si une $G^+(\AAA)$-repr\'esentation 
automorphe irr\'eductible discr\`ete $\Pi$ n'est pas cuspidale, 
alors $\trace(\Pi(fI_\theta))=0$. Si cette trace est non nulle alors $\Pi_p$ est de la forme $I(\chi)$ par le
Th\'eor\`eme \ref{thmcle} (i). Mais le th\'eor\`eme de Moeglin-Walspurger exclut\footnote{On aurait aussi pu arguer en
$\ell$, ou m\^eme en l'infini en supposant de plus que le caract\`ere
infinit\'esimal de la s\'erie $\theta$-discr\`ete attach\'ee \`a $f_\infty$
est suffisamment r\'eguli\`ere.}
la classe inertielle de $(M,\omega \times \check{\omega})_G$ comme composante locale possible d'une
repr\'esentation r\'esiduelle de $\GL_{2n}$ car $\omega$ n'est isomorphe \`a aucune
torsion non ramifi\'ee de $\check{\omega}$.
\end{pf}

\subsubsection{Preuve du th\'eor\`eme}\label{termeprincal}
Rappelons que la fonction $f_{\infty}$ d\'epend notamment du choix d'une s\'erie 
$\theta$-discr\`ete cohomologique de $G(\R)$, qui sont index\'ees comme on l'a vu au
\S\ref{prelimglnreel} par les repr\'esentations irr\'eductibles du groupe compact $\SO_{2n+1}(\R)$. 
Pour fixer les id\'ees, on choisit $T$ un tore maximal de ce dernier et on note $V_\lambda$ la repr\'esentation 
irr\'eductible de poids extremal $\lambda \in X^*(T)$. Pour chaque tel $\lambda$, on fixe 
un pseudocoefficient $f_{\infty}=f_{\lambda}$ de la s\'erie $\theta$-discr\`ete $\pi_{\lambda}$ associ\'ee 
de sorte que le support de tous ces $f_{\infty}$, $\lambda$ variant, soit contenu dans un m\^eme compact de $G(\R)$, ce qui est loisible, et on applique la 
Proposition \ref{formtrace} aux fonctions $$f^\lambda:=f_\lambda \otimes f_\ell \otimes f_p \otimes f^{\infty,\ell,p}.$$
\ps

Supposons que la trace de $f^{\lambda}$ dans une repr\'esentation cuspidale $\rho$ est non nulle. 
Alors $\rho_\infty$ est g\'en\'erique, 
a le m\^{e}me caract\`ere infinit\'esimal que $\pi_{\lambda}$, et donc lui est isomorphe.
D'apr\`es les propri\'et\'es de pseudocoefficients de $f_\ell$ (Prop. \ref{nbleftheta} (i)', 
noter que la triviale n'est jamais composante locale d'une cuspidale de $G$) et $f_p$ (Th\'eor\^eme \ref{thmcle} (i)), il suffit donc de montrer 
l'on peut choisir $\lambda$ de sorte que $\trace(R(f^{\lambda})I_\theta) \neq 0$,
soit encore que le terme g\'eom\'etrique correspondant de la formule de
trace de la Proposition \ref{formtrace} est non nul. \ps \ps

La somme du c\^ot\'e g\'eom\'etrique est \`a support inclus dans un ensemble fini $$\Sigma \subset \{G(\Q)\}_{\rm ell}$$
ind\'ependant de $\lambda$. On va voir qu'asymptotiquement en
$\lambda$ la classe de $\gamma_0$ porte le terme principal elliptique de la formule des
traces. Les propri\'et\'es essentielles de $\gamma_0$ sont r\'esum\'ees dans
le lemme suivant.

\begin{lemme}\label{generalitegamma0} L'\'el\'ement $\gamma_0$ est 
\`a $\theta$-conjugaison pr\`es l'unique \'el\'ement $\theta$-semisimple elliptique
de $G(\Q)$ tel que $\gamma_0\theta(\gamma_0)=-1$. Son centralisateur tordu est le sous-groupe
symplectique de $G$ de matrice $\gamma_0 J$. Sa classe de
$\theta$-conjugaison stable co\"incide avec sa classe de $\theta$-conjugaison.
\end{lemme}

\begin{pf} On a $\gamma_0\theta(\gamma_0)=-1$ et le centralisateur tordu de $\gamma_0$ est le groupe $\{g \in \GL_{2n},
g\gamma_0J{}^t\!g=\gamma_0J\}$ {\it i.e.} le groupe symplectique usuel ${\rm Sp}_{2n}$ de matrice $\gamma_0 J$ : 
$\gamma_0$ est donc bien $\theta$-semisimple elliptique. 

Les matrices antisym\'etriques inversibles \'etant toutes
congrues \`a $\gamma_0 J$ dans $\GL_{2n}(\Q)$, la classe de $\theta$-conjugaison de $\gamma_0$ co\"incide
exactement avec l'ensemble des \'el\'ements $\gamma$ tels que
$\gamma\theta(\gamma)=-1$. Pour la m\^eme raison, la classe de $\theta$-conjugaison stable de $\gamma_0$ est 
r\'eduite \`a sa classe de $\theta$-conjugaison (ou directement, $H^1(\Q,{\rm Sp}_{2n})=0$).
\end{pf}
\par \medskip

\begin{lemme} \label{nonnullitelemme} $\TO_{\gamma_0}(f^{\lambda,\infty})$ est une constante non nulle. 
\end{lemme}
\begin{pf} En effet, c'est le produit
$$\TO_\gamma(f_\ell)\cdot\TO_\gamma(f_p)\cdot\TO_\gamma(f^{\infty,\ell,p}).$$
Le terme $\TO_\gamma(f_\ell)$ (resp. $  \TO_{\gamma_0}(f_p)$) est non nul d'apr\`es la Proposition
\ref{nbleftheta} (ii) car $\gamma_0$ est $\theta$-semisimple elliptique
(resp. d'apr\`es le Th\'eor\`eme \ref{thmcle} (ii)). Le terme $\TO_\gamma(f^{\infty,\ell,p})$ est un r\'eel 
strictement positif car $\gamma_0 \in \prod_{v\neq \infty,\ell,p}G(\Z_v)$.
\end{pf}

Si $\gamma \in G(\Q)$ est un \'el\'ement $\theta$-semisimple, voyons-le comme un \'el\'ement $\theta$-semisimple 
de $G(\R)$ et consid\'erons sa norme $\No\gamma \in \SO_{2n+1}(\R)$ d\'efinie \`a la fin du
\S\ref{pseudoorb}. Par 
d\'efinition, $\No\gamma$ est un \'el\'ement dont la classe
de conjugaison sous $\SO_{2n+1}(\R)$ ne d\'epend que de la classe de
$\theta$-conjugaison de $\gamma$ sous $G(\R)$. D'apr\`es le Th\'eor\`eme~\ref{staborbreel}, on sait que pour une normalisation 
convenable des mesures on a pour tout $\lambda$, 
\begin{equation}\label{orbinfini} \TO_\gamma(f_\lambda)= \pm \trace(\No\gamma,\lambda).
\end{equation}

Notons que la classe de $\gamma_0$ est l'unique classe de norme centrale ({\it i.e.} triviale) dans $\SO_{2n+1}(\R)$. 
En effet, par d\'efinition $\No\gamma=1$ si et seulement si $\gamma\theta(\gamma)=-1$, et on conclut par le 
Lemme~\ref{generalitegamma0}. De plus,  $$|\TO_{\gamma_0}(f^\lambda)| = c \cdot \dim(V_\lambda) \neq 0$$ pour une 
certaine constante $c > 0$ d'apr\`es (\ref{orbinfini}) et le Lemme~\ref{nonnullitelemme}. 
Le Corollaire~\ref{corlimitecar} montre alors que lorsque $\lambda$ tends vers l'infini dans $X^*(T) \otimes \R$ en s'\'eloignant des murs, 
$$ \frac{\TO_{\gamma}(f^\lambda)}{\dim(V_\lambda)} \longrightarrow 0, \, \, \, \, \, \forall \gamma \in \Sigma\backslash \{\gamma_0\},$$
puis que
$$|\trace(R(f^\lambda)I_\theta)| \sim  v_{\gamma_0} \cdot c \cdot \dim(V_\lambda),$$
ce qui conclut. 

\begin{remarque}{\rm On d\'emontrera au \S \ref{calculdusigne} qu'en normalisant correctement les $f^\lambda$ le signe $\pm$ intervenant dans la formule (\ref{orbinfini}) 
est en fait $+1$ si $\gamma=\gamma_{0}$ (ind\'ependamment de $\lambda$). L'argument ci-dessus prouvera alors la formule
(\ref{geomasymptotique}) annonc\'ee dans l'introduction.}
\end{remarque}

\bigskip
\bigskip
\section{Non--\'evanouissement d'une int\'egrale orbitale
tordue}\label{nonvanorb}
\setcounter{equation}{0}

Le but de ce chapitre est de d\'emontrer les propri\'et\'es de la fonction $f_p$ utilis\'ee dans le 
Ch.~$3$. Bien que le probl\`eme soit local, la d\'emonstration utilise la formule des traces tordue 
pour $\GL(2n)$. Elle s'est r\'ev\'el\'ee assez difficile ; en revanche sa port\'ee est assez grande : 
nous esquissons dans le \S\ref{flocauxcusp} les cons\'equences de la m\'ethode. Celle--ci sugg\`ere des propri\'et\'es 
int\'eressantes de la formule de Plancherel tordue (associ\'ee \`a $\GL(2n)/\Q_p$ et \`a l'automorphisme 
$\theta$), ainsi que des int\'egrales orbitales tordues de coefficients de supercuspidales 
$\theta$--stables selon la parit\'e ({\it i.e.}, la nature symplectique ou orthogonale) de la repr\'esentation 
galoisienne associ\'ee.

La preuve de la non--nullit\'e de $\TO_{\gamma_0}(f_p)$ en l'\'el\'ement \og principal
\fg \,$\gamma_0$ du Ch.~$3$ repose sur le mod\`ele de Whittaker et un argument simple 
mais nouveau de positivit\'e (\S\ref{argumentpositivite}). Il nous a impos\'e de pr\'eciser le
Th\'eor\`eme~\ref{staborbreel}, en \'eliminant le signe implicite dans celui--ci 
gr\^ace \`a la normalisation \og de Whittaker \fg\, pour l'entrelacement $A$ entre $\pi$ et $\pi\circ \theta$.

Ceci nous impose h\'elas de remplacer l'automorphisme $\theta$ \og de
Waldspurger\fg\, (\S\ref{prelimglnreel}) par 
l'automorphisme $\theta_0$ respectant le mod\`ele de Whittaker. Les traductions n\'eces\-saires 
(pour aboutir au r\'esultat utilis\'e dans le Ch.~$3$) sont faites dans
le~\S\ref{waldwitt}.
\vskip2mm

\subsection{Pseudo--coefficients}\label{pseudocoeffcpomega} ${}^{}$ \ps Rappelons qu'on a fix\'e au
\S\ref{ssectenonce} une repr\'esentation supercuspidale $\omega$ de $\GL(n,\Q_p)$. 
Pour simplifier nous supposons $\omega$ {\it non} auto--duale, m\^eme modulo twist non--ramifi\'e :
\begin{equation*}
\omega\not\cong \omega\otimes \chi\, , \qquad \forall \chi \in X_{{\rm nr}}\ (\Q_p^\times)
\end{equation*}
o\`u $X_{{\rm nr}}(\Q_p^\times)\cong \C^\times$ est le tore des caract\`eres non ramifi\'es. 
(Le lecteur se convaincra ais\'ement que cette hypoth\`ese n'est pas fondamentale). 
On \'ecrira simplement $\omega\chi$ pour $\omega\otimes \chi$.

Rappelons que $\theta(g)=J_{2n} {}^t\!g^{-1} J_{2n}$, o\`u nous notons maintenant $J_{2n}$ la 
matrice $J$ du Ch.~$2$. De m\^eme,
\[
J_n\, \, =\, \, \left(
\begin{array}{ccc}
  &   &  1 \\
  &\adots   &   \\
 1 &   &
\end{array}
\right) \qquad \hbox{(taille }n)\,.
\]
\noindent On d\'esignera aussi par $\theta$ l'automorphisme $g\longmapsto
J_n{}^t\!g^{-1}J_n$ de $\GL(n)$.\ps

Nous nous int\'eressons \`a la repr\'esentation induite
\begin{equation}\label{induitecpo}
I(\chi) ={ \rm ind}_P^G
 (\omega\chi \otimes \widetilde\omega\  \chi^{-1})
\end{equation}
de $G=\GL(2n,\Q_p)$. L'induction est unitaire ; $P$ est le parabolique (triangulaire sup\'erieur) 
de type $(n,n)$ ; $\chi \in X_{{\rm nr}}$. Vu notre hypoth\`ese sur $\omega$, 
$I(\chi)$ est irr\'eductible pour tout $\chi$ ; elle est unitaire si, et seulement si, 
$\chi$ est unitaire (en effet si $I(\chi)$ est unitaire, donc hermitienne, on doit avoir
$$
\overline{I(\chi)} ={ \rm ind }(\overline\omega\overline\chi \otimes \overline{\widetilde \omega} \overline\chi^{-1}) \cong {\rm ind} (\widetilde\omega \chi^{-1}\otimes
\omega\chi)=\widetilde{I(\chi)}\,.
$$
Or $\overline\omega\cong \widetilde\omega$~; vu notre {\hypo} sur $\omega$, 
ceci implique $\omega\overline\chi^{-1}\cong \omega\chi$ soit $\omega\cong \omega(\chi\overline\chi)$~; 
alors $\chi\overline\chi$, donc $\chi$, est unitaire).\ps

Il est clair que $I(\chi)$ est isomorphe \`a sa duale et donc \`a $I(\chi)^\theta$. Pour obtenir un 
entrelacement explicite, remarquons que $\widetilde\omega
\chi^{-1}=\widetilde{\omega\chi}$ est isomorphe \`a $(\omega\chi)\circ \theta =(\omega\chi)^\theta$. 
Nous r\'ealisons donc $I(\chi)$  comme l'induite, {\it isomorphe} \`a
(\ref{induitecpo}) et d\'esign\'ee par la m\^eme notation :
\begin{equation*}
I(\chi)= {\rm ind}_P^G (\omega\chi \otimes (\omega\chi)^\theta)\,.
\end{equation*}

Si $V$ est l'espace de $\omega$, l'espace $\mathcal{L}(\chi)$ de $I(\chi)$ est donc form\'e des 
fonctions $f : G\lra V\otimes V$ v\'erifiant, pour
$p=\mathrm{diag}(m_1,m_2)n\, \, \in\, \,  P = (\GL(n)\times \GL(n))N$ :
$$
f(pg) = \delta_P(p)^{1/2}((\omega\chi)(m_1)\otimes(\omega\chi)^\theta m_2)f(g)\,.
$$
Soit $A_M : V\otimes V \lra V\otimes V$ l'op\'erateur 
\begin{equation*}
v\otimes w \longmapsto w \otimes v\,,
\end{equation*}
entrela\c{c}ant $\omega\chi\otimes (\omega\chi)^\theta$ et $(\omega\chi)^\theta \otimes \omega\chi$. On d\'efinit
\begin{equation*}
A_\theta : \mathcal{L}(\chi) \lra \mathcal{L} (\chi)
\end{equation*}
par
\begin{equation*}
A_\theta\ f(g) = A_M f(\theta g)\,.
\end{equation*}
\noindent On v\'erifie aussit\^ot que $A_\theta$ pr\'eserve $\mathcal{L}(\chi)$ ; il entrelace \'evidemment
les repr\'esentations $I(\chi)$ et $I(\chi)\circ \theta$ ; il est involutif. 
Noter que, dans la r\'ealisation compacte de $I(\chi)$, $A_\theta$ est ind\'ependant  de $\chi$.

\begin{prop}\label{enoncebisfp}
Il existe une fonction $f=f_p$ sur $\G=\G(\Q_p)$ ayant les propri\'et\'es suivantes :
\begin{enumerate}
\item[(i)] $f\in {\cal C}_c^\infty (\G)$,
\item[(ii)] $\trace\,(A_\theta I(\chi)(f))$ est une fonction alg\'ebrique sur $X_{nr}$, $>0$ sur 
les caract\`eres unitaires,
\item[(iii)] si $\pi$ est une repr\'esentation $\theta$--stable de $\G$ et 
$A:\pi\cong \pi\circ \theta$ est un op\'erateur d'entrelacement,
\end{enumerate}
\begin{equation*}
\trace\,(A\pi(f)))\,=\,0
\end{equation*}
si $\pi$ n'est pas une induite $I(\chi)$.
\end{prop}

On appellera parfois une telle fonction $f$ un {\it pseudo--coefficient positif}. 
Nous donnons une d\'emonstration simple de ce r\'esultat ; une autre d\'emonstration utiliserait le 
th\'eor\`eme de Paley--Wiener tordu de Rogawski~\cite{paleyweinerrog} : cf. \S\ref{argumentpositivite}.

\begin{pf} Le groupe $X_{nr}$ des caract\`eres non--ramifi\'es de $\Q_p^\times$ est un tore complexe ; 
soit $\Gamma \subset X_{nr}$ le sous--groupe (fini) des caract\`eres $\chi$ tels que 
$\omega\chi\cong \chi$ et soit $T=X_{nr}/\Gamma \cong \C^\times$. L'orbite de
Bernstein~\cite{bernstein} contenant les repr\'esentations $I(\chi)$ s'identifie \`a $T\times T$, 
par 
\begin{equation*}
(\chi_1,\, \chi_2) \longmapsto \mathrm{ind}_P^G (\omega\chi_1 \otimes \widetilde\omega\widetilde\chi_2)  = I
(\chi_1,\,\chi_2)\,,
\end{equation*}
et l'ensemble des $I(\chi)$ au sous-tore diagonal. Parmi les $I(\chi_1,\chi_2)$ seules les $I(\chi)$ sont 
autoduales. D'apr\`es des arguments bien connus, il existe $f\in C_c^\infty (G)$ telle que 
$\mathrm{trace}(A_\theta I (1)(f))\not= 0$. D'apr\`es Bernstein, on peut alors remplacer $f$ par une 
fonction v\'erifiant, de surcro\^{\i}t, (iii). La fonction $F(\chi)=\mathrm{trace}(A_\theta I(\chi)(f))$ 
est alg\'ebrique, non nulle, sur $T$ ; vu les propri\'et\'es de $\omega$, l'image du centre de Bernstein 
$\mathcal{Z}$ dans $\C[T\times T]$ est form\'ee de toutes les fonctions
alg\'ebriques.\ps
Soit $z$ le param\`etre sur $T\cong \C^\times$, et 
\begin{equation*}
F(\chi) =F(z) =\sum  a_n\,z^n\,
\end{equation*}
une s\'erie de Laurent finie. Si $h \in \mathcal{Z}$ est d'image
\begin{equation*}
H(\chi) =H(z) =\sum a_n z^{-n}
\end{equation*}
et si $g=h\star  f$, on a alors :
\begin{equation*}
\mathrm{trace}(A_\theta I(\chi)(g))= G(z)F(z)\,,
\end{equation*}
fonction positive sur les caract\`eres unitaires, d'o\`u (ii), la positivit\'e stricte r\'esultant aussit\^ot d'un argument de
compacit\'e (prendre une somme de telles fonctions). \end{pf}

Rappelons (Lemme~\ref{generalitegamma0}) que l'\'el\'ement
\[
\gamma_0 = 
\left(
\begin{array}{ccc}
  1_n&      \\
  &   -1_n   \\
\end{array}
\right) \in G
\]
\noindent donne le \og terme principal\fg \,  de la formule des traces tordue
(\S\ref{termeprincal}). On a $\mathcal{N}\gamma_0 = 1 \in \SO(2n+1)$. 
Si $f\in {\cal C}_c^\infty(G(\Q_p))$, l'int\'egrale orbitale tordue de $f$ en $\gamma_0$ est
\begin{equation}\label{niemeto}
\TO_{\gamma_0} (f) =\int_{\G(\Q_p)/I(\Q_p)} f(g\, \gamma_0 \, g^{-\theta}) \frac{dg}{di}. 
\end{equation}
Elle est stable (Lemme~\ref{generalitegamma0}).
\begin{theoreme}\label{ancien42}
Si $f$ est un pseudo--coefficient positif (Prop.~\ref{enoncebisfp}), $\TO_{\gamma_0}(f)\not= 0$.
\end{theoreme}

\vskip2mm

\subsection{Stabilit\'e d'un caract\`ere tordu}\label{staibilitecaractcpo} ${}^{}$\ps

Consid\'erons la distribution sur $G=G(\Q_p)$~:
\begin{equation}\label{tracetorcp}
f \longmapsto \mathrm{trace}(A_\theta\ I(x) (f))\,.
\end{equation}
Comme dans le chapitre $2$, \S\ref{parnorme}, on dispose sur $G(\Q_p)$ de la notion d'\'el\'ements 
$\theta$--semisimples, $\theta$--r\'eguliers : c'est le cadre original de Waldspurger~\cite{Wald}. 
D'apr\`es des r\'esultats g\'en\'eraux \cite[Thm. 1]{ClAB} on sait que la distribution (\ref{tracetorcp}), 
que l'on notera $\Theta_{\chi,\theta}$, est une fonction localement int\'egrable sur $G$, ${\cal C}^\infty $ 
sur les \'el\'ement $\theta$--r\'eguliers, invariante par $\theta$--conjugaison.

\begin{prop}\label{stabtchit}
Le caract\`ere tordu $\Theta_{\chi,\theta}$ est invariant par $\theta$--conjugaison {\it stable} 
(sur les \'el\'ements fortement $\theta$--r\'eguliers).
\end{prop}

\begin{pf} En effet, toute formule pour le caract\`ere tordu de la repr\'esentation induite --- 
par exemple, le th\'eor\`eme d'Atiyah--Bott \cite[Prop.~6]{ClAB} --- montre que le support de 
$\Theta_{\chi,\theta}$ est contenu dans l'ensemble des \'el\'ements de $G$ qui sont 
$\theta$--conjugu\'es \`a un \'el\'ement de $M$. Si $\gamma\in G$ est un \'el\'ement 
fortement $\theta$--r\'egulier, son centralisateur tordu $I$ est un tore (\S\ref{parnorme}), de dimension $n$. L'ensemble des classes de $\theta$--conjugaison dans la classe de conjugaison stable de $\gamma$ s'identifie \`a 
\begin{equation*}
H^1 (\Q_p,I) = \ker [H^1 (\Q_p,I) \lra H^1(\Q_p,G)]\,.
\end{equation*}
Or le centralisateur tordu dans $M$ d'un \'el\'ement fortement $\theta$--r\'egulier $g$ de $M$ s'identifie 
\`a un tore maximal de $\GL(n)$. En effet, si $\theta$ d\'esigne l'automorphisme $m\longmapsto J_n {}^t m^{-1}J_n$ 
de $\GL(n)$, l'automorphisme $\theta$ de $\GL(2n)$ restreint \`a $M$ est $(m_1,m_2)\longmapsto (\theta m_2,\theta m_1)$. 
Par l'isomorphisme $(m_1,m_2)\longmapsto (m_1,{}^\theta m_2)$ il est conjugu\'e \`a $$(m_1,m_2) \longmapsto (m_2,m_1).$$ 
On est alors ramen\'e au cas, facile, de la $\sigma$--conjugaison pour le changement de base en une place d\'ecompos\'ee 
\cite[Ch.~I.5]{AC}. Au vu des dimensions, on en d\'eduit que les centralisateurs tordus $I$
 de $g$ dans $M$ et $G$ co\"{\i}ncident. Donc $I$ est cohomologiquement trivial, d'o\`u la proposition.
\end{pf}

\begin{corollaire}\label{coriost}
Si $f\in {\cal C}_c^\infty (G)$ est un pseudo--coefficient positif, les int\'egrales orbitales tordues stables de 
$f$, en les \'el\'ements fortement $\theta$--r\'eguliers, ne sont pas identiquement nulles.
\end{corollaire}

\begin{pf} Il suffit d'appliquer la formule d'int\'egration de Weyl sur $\theta \G$ \`a la 
fonction $f\cdot\Theta_{\chi,\theta}$.
\end{pf}

\vskip2mm

\subsection{De Waldspurger \`a Whittaker}\label{waldwitt}

Nous avons jusqu'ici utilis\'e l'automorphisme $\theta$ de $G$ en suivant Waldspurger ; mais les arguments qui 
suivent vont reposer sur le mod\`ele de Whittaker, auquel il n'est  pas adapt\'e. Soit~donc

\begin{equation*}
D=\mathrm{diag}(1,-1,1,\ldots,1,-1)\,,
\end{equation*}

\[J_0 = D J=
\left(
\begin{array}{ccccc}
  & & & & 1 \\
  & & & \adots & \\
  & & -1 & & \\
  & 1 & & & \\
-1 & & & &   
\end{array}
\right)\,,\qquad J_0^2 = -1
\]
et $\theta_0:g\longmapsto J_0{}^t g^{-1}J_0^{-1}$ $(g\in \GL(2n))$.\ps

Adaptons rapidement les r\'esultats pr\'ec\'edents, en rempla\c{c}ant $\theta$ par $\theta_0$.\footnote{Le lecteur choqu\'e par 
notre inconstance m\'ethodologique refera la 
construction de la norme (Ch. $2$) en suivant Kottwitz et Shelstad \cite{KS}.} 
Rappelons tout d'abord que $\theta_0$ fixe un \'epinglage pour le couple 
$(B,T)$ form\'e du groupe de Borel triangulaire sup\'erieur et du tore diagonal. Si $\alpha$ est un caract\`ere additif non trivial 
(d'un corps local $F=\R$ ou $\Q_p$, ou bien de $\A)$ et si $\psi$ est le caract\`ere du sous--groupe unipotent $N$ de $B$ donn\'e par
\[
\psi
\left(
\begin{array}{cccccc}
 1 &x_1 &  &  & \ast \\
  &1   & &  \\
  &  & \ddots&  &\\
  && & 1 &x_{2n-1}\\
  && & &1   
\end{array}
\right) = \alpha (x_1+\cdots + x_{2n-1})\,,
\]
on a $\psi(\theta_0 n) =\psi(n)$ $(n\in N)$. Soit $\pi$ une repr\'esentation irr\'eductible g\'en\'erique 
autoduale (de $G(F)$ ou $G(\A))$. Il existe alors un unique op\'erateur d'entrelacement involutif 
$A_{\theta_0}^{\Wh}$ entre $\pi$ et $\pi\circ \theta$ tel que $A_{\theta_0}^{\Wh}(\lambda)=\lambda$, 
$\lambda$ \'etant l'unique fonctionnelle de Whittaker associ\'ee \`a $\psi$ (unique \`a un scalaire pr\`es) 
et $A_{\theta_0}$ op\'erant par l'action duale. On dira que $A_{\theta_0}^\Wh$ est la normalisation de 
Whittaker pour l'entrelacement involutif (\`a priori d\'efini modulo $\pm 1$).

Si on n\'eglige le signe, on peut construire, en $p$, un op\'erateur d'entrelacement $A_{\theta_0}$ sur $\pi=I(\chi)$ en posant :
\begin{equation*}
A_{\theta_0} = \pi (D) A_\theta \,.\hbox to 1cm{}
\end{equation*}

On a alors
\[
\begin{array}{ccc}
 A_{\theta_0}^2 &= \pi(DJ {}^t D^{-1} J) A_{\theta}^2     \\
  &=\pi(J_0^2)A_\theta^2   \hfill  \\
  & =\pi(-1)=1 \hfill   
\end{array}
\]
puisque $\pi$ est de caract\`ere central trivial. Par ailleurs, notons, pour
$\delta\in \G$, 
${\rm T}_0{\rm O}_\delta(f)$ l'int\'egrale orbitale tordue d\'eduite de $\theta_0$. 
Si $h(g) = f(gD)$, on v\'erifie aussit\^ot que
\begin{equation*}
{\rm T}_0{\rm O}_\delta (h)=\TO_\gamma(f)
\end{equation*}
pour
\begin{equation*}
\delta =\gamma D\,.
\end{equation*}
Par ailleurs $\trace\,(\,A_{\theta_0}\pi(h)\,) = \trace (\,A_\theta \pi(f)\,)$. La fonction $h$ v\'erifie donc les propri\'et\'es de la 
Proposition~\ref{enoncebisfp} relativement \`a $\theta_0$. Le 
Th\'eor\`eme~\ref{ancien42} est \'equivalent~\`a
\begin{equation}\label{nunnulw}
{\rm T}_0{\rm O}_{\delta _0}(h) \not= 0\,,
\end{equation}
o\`u
\begin{equation*}
\delta _0 = \gamma_0 D\,.
\end{equation*}

Des consid\'erations analogues s'appliquent \`a la place r\'eelle, 
en rempla\c{c}ant la norme du Ch.~$2$  par 
$\mathcal{N}_0(\delta) = \mathcal{N}(\delta D)$. En particulier $\delta _0$ est de norme~$1$, l'analogue du Th\'eor\`eme~\ref{staborbreel} est v\'erifi\'e, et $\delta _0$ est, 
localement ou globalement, l'unique \'el\'ement dans sa classe de conjugaison tordue stable. Si $\pi$ est une  repr\'esentation cohomologique de $G(\R)$ comme dans le  Ch.~$2$, on d\'efinit $A_{\theta_0}$ comme dans le cas $p$--adique~; alors, en posant encore $h_\pi(g) =f_\pi(gD)$ :
\begin{equation}\label{toinfiniw}
{\rm T}_0{\rm O}_{\delta_0}(h_\pi)= \varepsilon(\pi) \dim (\pi_H)\, ,
\end{equation}
o\`u $\varepsilon$ est un signe. La d\'emonstration va nous amener \`a pr\'eciser (\ref{toinfiniw}) 
ainsi que le Th\'eor\`eme~\ref{ancien42}. 
Pour simplifier les notations, $\TO$, $\mathcal{N}$\,\dots\,\, d\'esignent dor\'enavant les variantes relatives 
\`a~$\theta_0$. 

\begin{prop}\label{ancienprop45} Si l'op\'erateur $A_{\theta_0}=A_{\theta_0}^\Wh$ est normalis\'e pour le mod\`ele de Whittaker,
\begin{enumerate}
\item [(i)] $\TO_{\delta_0} (h_p) > 0$,
\item[(ii)] $\TO_{\delta_0}(h_\pi) = \dim (\pi_H)$ pour toute repr\'esentation cohomologique $\pi$ de 
$G(\R)$ (pour une normalisation fixe des mesures sur $G(\R)$ et le centralisateur tordu de~$\delta_0$).
\end{enumerate}
\end{prop}

Nous reportons au \S\ref{argumentpositivite} et \S\ref{calculdusigne} la d\'emonstration de cette proposition. Nous aurons enfin besoin du 

\begin{lemme}\label{lemmenormnonram}
Soit $\pi$ une repr\'esentation auto--duale de la s\'erie principale : 
$\pi={\rm ind} (\chi_1,\ldots \chi_n,\, \chi_n^{-1},\ldots,\chi_1^{-1})$, les $\chi_i$ \'etant non ramifi\'es. 
Alors l'op\'erateur $A_{\theta_0}^\Wh$ donn\'e par $f(g) \longmapsto f(\theta_0(g))$ entrelace $\pi$ et $\pi\circ \theta_0$~; il est involutif et op\`ere par $+1$ sur l'espace de Whittaker~; il op\`ere trivialement sur le vecteur non ramifi\'e.
\end{lemme}

C'est \'evident (rappelons que la fonctionnelle de Whittaker est donn\'ee par
\begin{equation*}
f\longmapsto \int_N \ \ f(w_0 n)\psi(n) dn
\end{equation*}
si $f$ est \`a support dans $B w_0 N$, $w_0$ \'etant l'\'el\'ement  de plus grande longueur du groupe de Weyl).

\vskip2mm

\subsection{Formule des traces (stable)}\label{preliftst}${}^{}$\ps

Soit $f$ une fonction ${\cal C}^\infty$ \`a support compact sur $G(\A)$, d\'ecompos\'ee, 
donc $f=\otimes'_v f_v$. 
\`A la place r\'eelle, on choisit une repr\'esentation cohomologique $\pi_\infty$ autoduale de $G(\R)$~; 
on notera parfois $\mu$ le para\-m\`etre
\begin{equation*}
(m_1,\ldots,m_n) \in \N^n,\ m_1 \geq \cdots \geq m_n \geq 0
\end{equation*}
de $\pi=:\pi_{\mu}$, cf. \S\ref{prelimglnreel} ; c'est donc le plus haut poids de la repr\'esentation de $\SO(2n+1)$ associ\'ee \`a 
$\pi$. Soit $\theta_{\mu}$ le caract\`ere de celle-ci. Alors $f_\pi = f_\mu$ est le pseudo--coefficient de $\pi$ pour la normalisation de Whittaker de l'entrelacement $A_{\theta_0}$ (fonction not\'ee $h_\pi$ dans le paragraphe pr\'ec\'edent).

En la place $\ell$, $f$ est un pseudo--coefficient tordu pour la repr\'esentation de Steinberg (\S\ref{preliminaireLef}). En $p$, c'est un pseudo--coefficient positif pour $I(\chi)$ (\S\ref{pseudocoeffcpomega}). On utilise de nouveau $A_{\theta_0}$, avec la normalisation de Whittaker. En $v\not= \infty,\, \ell,\, p$, 
$f_v$ est pour l'instant arbitraire.\ps

Soit $I_{\theta_0} : \varphi (x) \longmapsto \varphi(\theta_0 x)$ l'op\'erateur d'entrelacement de $L^2(A\, \G(\Q)\ba \G(\A))$ donn\'e par $\theta_0$. Gr\^ace aux propri\'et\'es particuli\`eres des fonctions $f_\infty$ et $f_\ell$, on a tout d'abord~:
\begin{equation} \label{tform}
\sum_{\delta \in \{\G(\Q) \}_{\rm ell}} {\rm vol} (I_\delta ) \TO_\delta (f) = \sum_\pi \mathrm{trace}
 (I_{\theta_0}\pi(f))\,. 
\end{equation}

On a fix\'e une mesure de Haar $dg_\A$ sur $\G(\A)$, qui d\'efinit la trace dans le membre de droite. Dans le membre de gauche, $\{\G(\Q)\}_{\rm ell}$ est l'ensemble des classes de conjugaison tordue d'\'el\'ements $\theta$--elliptiques \cite[p.~508]{ITF} de $\G(\Q)$ ; $\TO_\delta $ est l'int\'egrale orbitale tordue (ad\'elique) d\'efinie par $dg_\A$ et une mesure de Haar $di_\A$ sur le centra\-lisateur tordu connexe $I_\delta $ de $\delta $~; ${\rm vol}(I_\delta )$ est la mesure de $I_\delta (\Q)\ba I_\delta (\A)$. La somme est finie. Dans le membre de droite, $\pi$ parcourt les repr\'esentations {\it cuspidales} $\theta_0$--stables ($\equiv$ autoduales) de~$\G(\A)$.

L'\'egalit\'e (\ref{tform}) entre le c\^ot\'e g\'eom\'etrique et le c\^ot\'e spectral est d\'emontr\'ee par Arthur \cite[\S7]{ITF}. L'expression du c\^ot\'e g\'eom\'etrique est \cite[Cor.~7.4]{ITF}, \'etant donn\'ees les propri\'et\'es de $f_\infty$ et $f_\ell$. De m\^eme, le c\^ot\'e spectral est \cite[Cor.~7.2]{ITF}. Noter que le membre de droite de \cite[Cor.~7.2]{ITF} est \`a priori plus compliqu\'e : il contient les repr\'esentations du spectre r\'esiduel, ainsi que certaines induites $\theta$--discr\`etes : ces repr\'esentations sont \'elimin\'ees par $f_p$ (cf. Prop.~\ref{formtrace}). Il contient aussi une sommation sur $t$ (norme du caract\`ere infinit\'esimal de $\pi_\infty$) mais $t$ est fix\'e par le choix de $f_\infty$. En particulier, $f_{{\rm fin}}=\otimes'_{v \neq \infty} f_v$ \'etant fix\'ee, et donc aussi la ramification, le nombre de termes du membre de droite est fini.\ps

Notons ${\rm T}(f)$ l'expression (\ref{tform}). Si $\delta \in G(\Q)$ est un \'el\'ement $\theta$--semisimple elliptique, $\TO_\delta (f)$ est un produit
\begin{equation*}
\prod_v \TO_\delta (f_v)
\end{equation*}
d'int\'egrales orbitales locales. Si $\delta$ est $\theta$--semisimple et fortement r\'egulier, posons
\begin{equation}\label{stofreg}
 \STO_\delta (f) =\displaystyle \prod_v \STO_\delta (f_v)
\end{equation}
o\`u
\begin{equation*}
\STO_\delta (f_v) =\displaystyle \sum_{\delta'} \TO_\delta (f_v) 
\end{equation*}
la somme portant sur les \'el\'ements $\delta '$ stablement $\theta$ conjugu\'es \`a $\delta $. Il r\'esulte alors de la stabilisation du terme r\'egulier de la formule des traces par Kottwitz--Shelshad \cite{KS} que :

\begin{lemme}\label{ancienlemme47}
Soit $q$ un nombre premier diff\'erent de $p$, $\ell$ et soit $f=\otimes'_v f_v$ une fonction v\'erifiant les conditions pr\'ec\'edentes, $f_q$ \'etant de plus de support $\theta$--fortement r\'egulier. Alors
\begin{equation*}
{\rm T}(f) = \alpha(\G) \sum_\delta  \STO_\delta (f)\,,
\end{equation*}
la somme portant sur les m\^emes \'el\'ements qu'en $(\ref{tform})$ mais modulo $\theta$--conjugaison stable dans~$\G$.
\end{lemme}

Ceci r\'esulte de \cite[Ch. 7]{KS} (cf. en particulier (7.4.3)~; $\alpha(\G)$ est une constante $>0$), et des propri\'et\'es de $f_\infty$ : c'est une fonction stabilisante (\cite[\S3.8]{labesse} et notre Corollaire \ref{corollairekappa}). On en d\'eduit :

\begin{prop}\label{ancienprop48} Il existe une repr\'esentation cuspidale, $\theta$--stable $\pi$ de $\G(\A)$ telle que :
\begin{enumerate}
\item [(i)] $\pi_\infty \simeq \pi_{\mu}$,
\item[(ii)] $\pi_\ell \simeq \St_\ell$,
\item[(iii)] $\pi_p \simeq I(\chi)$ pour un caract\`ere non ramifi\'e $\chi$,
\end{enumerate}
\end{prop}

\begin{pf} Soit en effet $f_q$ la fonction caract\'eristique d'un ouvert compact $\omega_q$ de $\G(\Q_q)$ 
constitu\'e d'\'el\'ements $\theta$-semisimples fortement r\'eguliers. Alors
\begin{equation}\label{expressionT}
{\rm T}(f) =\alpha(G) \sum_\delta  \STO_\delta (f), 
\end{equation}
la somme \'etant uniform\'ement finie quand le support de $f$ est fix\'e \cite[p.~106]{KS}. L'application $\mathcal{N}_0$ donne une bijection entre classes de $\theta_0$--conjugaison semisimples $\delta$ stable et classes de conjugaison semisimples stable dans $\SO(2n+1)^*$ (quasi--d\'eploy\'e) \cite{KS}~; celles--ci sont donn\'ees par les polyn\^omes r\'eciproques
\begin{equation*}
P(X)\in \Q[X],\qquad P(X) = X^{2n} P(X^{-1})
\end{equation*}
de degr\'e $2n$. On \'ecrira $P=P_\delta $. En $v=\infty,\, \ell, \, p$,
choisissons un \'el\'ement $\theta$-semisimple fortement r\'egulier $\delta_v \in
\GL_{2n}(\Q_v)$ de sorte que 
$$\STO_{\delta _v}(f_v) \neq 0,$$ 
ce qui est loisible d'apr\`es le Th\'eor\`eme~\ref{staborbreel}, la
Proposition~\ref{nbleftheta} et le Corollaire~\ref{coriost}.
Cette non-annulation persiste alors dans un voisinage ouvert compact assez
petit $\omega_{\infty}\times\omega_\ell \times \omega_p$ de $(\delta_{\infty},\delta_{\ell},\delta_{p})$ dans
$\GL(2n,\AAA_{{\infty\,\ell\,p}})$. \ps

Par approximation faible on peut trouver $\delta' \in \GL(2n,\Q)$ tel que (pour le plongement diagonal) $\delta'$ appartient \`a $\prod_{v \in S}\omega_{v}$, 
o\`u $S=\{\infty,\ell, p, q\}$. Soient $\omega^S=\prod_{v\notin S}\omega_v$ un voisinage ouvert compact
d\'ecompos\'e fixe de $\delta'$ dans $\GL_{2n}(\AAA^S)$, et $f^{S}$ la fonction
caract\'eristique de $\omega^S$. Si 
$$\STO_{\delta} ( f_{S} \otimes f^{S} ) \neq 0 \qquad (\delta \in \G(\Q)),$$
 le polyn\^ome $P_\delta(X)$ a les propri\'et\'es suivantes. 
Ses coefficients sont des
$S'$-entiers o\`u $S'=S \cup \{v, \omega_v \neq \GL_{2n}(\Z_v)\}$. En les places $v \in
S'$ ils sont born\'es par la donn\'ee de $\prod_{v\in S'}\omega_v$ : il n'y a
qu'un nombre fini de tels polyn\^omes.

Choisissons $\omega_q$ contenant $\delta'_q$ assez petit. Si $\STO_\delta(f) \neq 0$ alors les coefficients 
de $P_\delta$ sont uniquement 
d\'etermin\'es donc $P_\delta=P_{\delta'}$ : $\delta$ est stablement conjugu\'e \`a $\delta'$. 
La somme (\ref{expressionT}) ne porte que sur la classe de conjugaison stable de $\delta'$.
Comme $\delta'$ est fortement r\'egulier $\STO_{\delta'}(f_q\otimes f^S) > 0$, d'o\`u la proposition.
\end{pf}

\vskip2mm

\subsection{D\'emonstration du Th\'eor\`eme~\ref{ancien42}}\label{argumentpositivite}

Dans ce paragraphe nous d\'emontrons le th\'eor\`eme, en supposant pour l'instant la partie archim\'edienne (ii) de la Proposition~\ref{ancienprop45}. Noter que d'apr\`es le Ch.~$2$, celle--ci est \'equivalente~\`a 

\begin{prop}\label{propequivancien45} Pour $\delta \in G(\R)$ de norme elliptique
$$
\TO_\delta (h_\pi) = \varepsilon(\delta )\ \trace\, \pi_H(\mathcal{N}_0 \gamma)
$$
o\`u le signe $\varepsilon(\delta )$ ne d\'epend pas de $\mu$~; de plus $\varepsilon(\delta _0)=1$.
\end{prop}

Comme dans le Ch.~$2$, cette assertion n'est vraie que pour des choix convenables des mesures (positives) sur les centralisateurs tordus, choix que nous ne pr\'eciserons pas.

Nous utiliserons le r\'esultat \'el\'ementaire suivant :

\begin{prop}\label{ancien49} Soit $G$ un groupe de Lie compact, $\gamma_0=1,\gamma_1,\ldots, \gamma_N$ des classes de conjugaison distinctes dans $G$, et $\lambda_0,\ldots, \lambda_N$ des nombres complexes. Supposons que, pour tout caract\`ere irr\'eductible $\rho$ de~$G$,
$$\sum_{i=0}^N \lambda_i \ \hbox{\rm trace}\ \rho(\gamma_i)$$
est un nombre r\'eel $\geq 0$, et que cette somme soit strictement positive pour un caract\`ere $\rho_0$. Alors $\lambda_0$ est un r\'eel $>0$.
\end{prop}

\begin{pf} Soit en effet $\OO_i$ la classe de conjugaison de $\gamma_i$ et $\mu_i$ la mesure invariante (normalis\'ee) sur $\OO_i$, vue comme une distribution sur $G$. Si $\mu=\Sigma\ \lambda_i\ \mu_i$,
\begin{equation*}
\hbox{\rm trace }(\rho(\mu)) \geq 0\qquad 
(\rho \hbox{ \rm  irr\'eductible)}.
\end{equation*}

D'apr\`es l'extension \`a $\G$ de la transformation de Fourier des distributions, la distribution $\mu$ s'\'ecrit
\begin{equation}\label{transffourieriorb}
\mu=\sum_\rho a_\rho\ \Theta_\rho
\end{equation}
o\`u $\Theta_\rho$ est le caract\`ere de $\rho$ et $a_\rho$ est une fonction \`a croissance lente sur $\widehat G$. Il r\'esulte des relations d'orthogonalit\'e et de (\ref{transffourieriorb}) que
\begin{equation*}
a_\rho \geq 0,\qquad a_{\rho_0}>0\,.
\end{equation*}

Soit $f$ une fonction ${\cal C}^\infty$ sur $G$ telles que le support de $f \ast f^*$ ne rencontre pas $\OO_i$ pour $i>0$, $f^*(g) $ \'etant $\overline{f(g^{-1})}$. Alors
\begin{equation*}
\mu(f\ast f^*) = \lambda_0 \| f\|^2 = \sum_\rho a_\rho \hbox{ \rm trace } \rho(f\ast f^*)\,,
\end{equation*}
d'o\`u le r\'esultat.
\end{pf}

Soit $\pi$ la repr\'esentation exhib\'ee dans la
Proposition~\ref{ancienprop48}. Soit $q$ un nombre premier tel que $\pi_{q}$ est ramifi\'ee. La repr\'esentation $\pi_q$ est g\'en\'erique, $\theta_0$--stable et m\^eme temp\'er\'ee d'apr\`es Harris et Taylor \cite{HT}. 
Elle s'\'ecrit donc :
\begin{equation*}
\pi_q = {\rm ind}_P^G \ \ (\delta _1\otimes\cdots \otimes \delta _r) 
\end{equation*}
o\`u $P$ est un parabolique de type $(n_1,\ldots n_r)$, $\delta _i$ est unitaire et de carr\'e int\'egrable pour $\GL(n_i)$, et o\`u l'on peut supposer : \begin{itemize}\ps
\item[(i)] pour $i\leq t$, $n_{r+1-i}=n_i$ et\, $\delta _i=\delta _{r+1-i}\circ \theta_0$,\ps
\item[(ii)] pour $t<i<r+1-t$, $\delta _i \cong \delta _i \circ \theta_0$.\ps
\end{itemize}

(On a d\'esign\'e par $\theta_0$ l'automorphisme du \S\ref{waldwitt}, pour $\GL(n_i)$, peut-\^etre en rempla\c{c}ant $D$ par $-D$).

S'il n'y a que des blocs de type (i), on construit comme dans le \S\ref{pseudocoeffcpomega} un entrelacement $A_{\theta_0}$, 
Whittaker--normalis\'e, pour $\pi=\pi_q$, ainsi que pour toutes les induites des $\delta _i$ tordues par des caract\`eres non--ramifi\'es. 
En g\'en\'eral, un tel $A_{\theta_0}$\break existe~: composer avec un op\'erateur d'entrelacement normalis\'e pour les blocs de type (ii)~; 
on obtient un op\'erateur op\'erant \`a priori par $\{\pm 1\}$ sur la fonctionnelle de Whittaker. L'op\'erateur d'entrelacement normalis\'e \'etant holomorphe, le signe est constant.

Il r\'esulte alors du th\'eor\`eme de Paley--Wiener de Rogawski \cite{paleyweinerrog} que :

\begin{lemme}\label{consPWR}
Il existe une fonction $f_q\in C_c^\infty(G(\Q_q))$ telle que

\begin{enumerate}
\item [(i)] parmi les repr\'esentations $\theta$--stables g\'en\'eriques $\pi$ de $G(\Q_q)$, les $I(\delta _i\otimes \chi_i)$ (pour des $\chi_i$ non ramifi\'es tels que l'induite est $\theta$--stable) sont les seules telles que
\begin{equation*}
{\rm trace}(A_{\theta_0} \ \pi\  (f_q)) \not= 0,
\end{equation*}

\item[(ii)] ${\rm trace}\ (A_{\theta_0}^\Wh I(\delta _i\otimes \chi_i)(f_q)) >0$ 

\noindent ($\chi_i$ unitaires non ramifi\'es, $I$ suppos\'ee $\theta$--stable).
 
\end{enumerate}

\end{lemme}

Nous fixons ainsi $f_q$ pour tout nombre premier $q \neq \ell, p$ en lequel $\pi_{q}$ est ramifi\'ee ; $f_\ell$ et $f_p$ l'ont \'et\'e, et nous faisons varier $f_\infty=f_\mu$ avec le poids $\mu$. 
Nous allons appliquer la formule des traces (\ref{tform}), $f_{\ell'}$ (pour les autres nombres premiers) \'etant l'unit\'e de 
l'alg\`ebre de Hecke sph\'erique. Notons $S$ la r\'eunion de $\{\infty\}$ et de l'ensemble des 
nombre premiers $q$ tels que $\pi_{q}$ est ramifi\'ee.

Sur l'espace des formes paraboliques sur $\G(\A)$, on dispose d'une fonctionnelle de Whittaker globale, $\theta_0$--invariante :
\begin{equation*}
\varphi \longmapsto \int_{N(\Q)\ba N(\A)} \psi^{-1}(n) \varphi (n) dn = \lambda(\varphi)\,.
\end{equation*}

Soit $\pi\subset L^2(A\,G(\Q)\ba G(\A))$ une repr\'esentation cuspidale $\theta$--stable, o\`u $A=\R_+^\times \subset \G(\R)$, et soit  $A_{\theta_0}= I_{\theta_0} |_\pi$. Le Lemme~\ref{lemmenormnonram} implique que $A_{\theta_0}$ s'\'ecrit 
\begin{equation*}
\left(\bigotimes_{v\in S} A_{\theta_{0,v}}\right) \otimes A_{\theta_0}^S,
\end{equation*} 
les $A_{\theta_{0,v}}$ sont normalis\'es et  $A_{\theta_{0}}^S$ est le produit tensoriel (bien d\'efini) des op\'erateurs non ramifi\'es.

Nous utilisons la formule des traces (\ref{tform}). Le membre de gauche s'\'ecrit 
\begin{equation*} \sum_{\pi} \prod_{v\in S} \trace\,(\,A_{\theta_0}\pi_v(f_v)\,). 
\end{equation*}
Pour tout $\pi$ (cuspidale, $\theta$-stable), chaque terme du produit est $\geq 0$ ; la somme est en fait finie car $\pi_\infty$ est cohomologique et le niveau de $\pi$ est fix\'e. Le membre de droite comporte un nombre fini, fixe de termes, m\^eme quand $\mu$ (et $f_\mu$) varie. Il s'\'ecrit

\begin{equation}\label{mmbredte}
a(\delta _0) \deg(\theta_\mu) + \sum_{\delta \not= \delta _0} a(\delta ) \theta_{\mu}(\mathcal{N}\delta ),
\end{equation}

o\`u 
\begin{equation*}
a(\delta ) = \varepsilon(\delta ) {\rm vol} (I_\delta ) \prod_{v\not=\infty} \TO_\delta (f_v)
\end{equation*}
pour tout $\delta $, et $\varepsilon(\delta _0)=1$. Le seul \'el\'ement de norme $1$ dans $\SO(2n+1)$ 
est $\delta _0$~; on peut r\'e\'ecrire (\ref{mmbredte}) en regroupant les $\delta \not= \delta _0$ selon 
la conjugaison (tordue) stable, de sorte que (\ref{mmbredte}) est de la forme
\begin{equation*}
a(\delta _0) \deg (\theta_\mu) + \sum_{\gamma'} a(\gamma') \theta_{\mu}(\gamma')
\end{equation*}
et les $\gamma'\in \SO(2n+1)$ sont des classes de conjugaison distinctes.

On obtient alors une somme de valeurs de caract\`eres de $\SO(2n+1)(\R)$ v\'erifiant les hypoth\`eses 
de la Proposition~\ref{ancien49}, la positivit\'e stricte r\'esultant de la Proposition~\ref{ancienprop48}. Donc $$\prod_{v\not= \infty} \TO_{\delta _0}(f_v) >0,$$ d'o\`u le Th\'eor\`eme~\ref{ancien42}.

\vskip2mm

\subsection{D\'emonstration de la Proposition~\ref{ancienprop45}}\label{calculdusigne}

Consid\'erons d'abord une extension quadratique r\'eelle $F$ de $\Q$. 
Soit $\infty$, $\infty'$ les deux places archim\'ediennes, et soit 
$f= \otimes'_v f_v \in {\cal C}_c^\infty(\G(\A_F))$ une fonction \'egale \`a la fonction unit\'e sph\'erique aux places finies, alors que
\begin{equation*}
f_{\infty}=f_{\mu}\,,\qquad f_{\infty'}= f_{\mu'}\,,
\end{equation*}
$\mu$, $\mu'$ \'etant deux poids dominants pour $\SO(2n+1)$. 
Les consid\'erations du \S\ref{preliftst} s'appliquent, $f_{\infty}$ et $f_{\infty'}$ 
\'etant \og cuspidales\fg\,  au sens d'Arthur. L'argument du \S\ref{argumentpositivite} 
donne, pour tous $\mu$, $\mu'$~:
\begin{equation}\label{encoreuneformule}
a(\delta_0) \deg(\theta_\mu) \deg(\theta_{\mu'}) + \sum_{\gamma_1 \not= 1} a(\gamma_1,\mu,\mu') 
\theta_\mu(\gamma) \theta_{\mu'}(\gamma) \geq 0
\end{equation}
o\`u
\begin{equation*}
a(\delta _0) = \varepsilon(\mu) \varepsilon(\mu') {\rm vol} (I_{\delta _0}) \prod_{v \nmid \infty} \TO_{\delta_0}(f_v)\,,
\end{equation*}
et o\`u $\gamma_1=(\gamma,\gamma')$ parcourt un ensemble fini de classes de conjugaison du groupe $\SO(2n+1)(\R)^2$, 
et $a(\gamma_1,\mu,\mu')=a(\gamma_1)\varepsilon(\mu)\varepsilon(\mu')$. 
Noter que si $\delta \not=\delta _0$ et si $\gamma_1=(\gamma,\gamma')$ est associ\'e \`a 
$\delta$, $\gamma$ et $\gamma'$ sont $\not= 1$, la conjugaison stable (globale) \'etant indiqu\'ee par la conjugaison stable en {\it une} place. Les signes $\varepsilon(\mu)$, pour l'instant inconnus, sont donn\'es par l'\'egalit\'e
\begin{equation*}
A_{\theta_0}^{\Wh}= \varepsilon(\mu) A_{\theta_0}
\end{equation*}
o\`u les deux op\'erateurs sont des entrelacements involutifs pour $\pi(\mu)$, celui de gauche \'etant Whittaker--normalis\'e, et celui de droite donnant, comme dans le Ch.~2, l'identit\'e correcte (sans signe) pour l'int\'egrale orbitale en~$\delta _0$.

Fixons $\mu'$ et faisons tendre $\mu$ vers l'infini \og loin des murs\fg\, . Le terme dominant de (\ref{encoreuneformule}) est celui relatif \`a 
$\delta_0$. Les int\'egrales orbitables tordues, aux places finies, \'etant $>0$ en $\delta _0$, on voit que $\varepsilon(\mu)\varepsilon(\mu')$ est $>0$ pour $\mu$ \og assez loin des murs\fg\,. Variant maintenant $\mu'$, on en d\'eduit que $\varepsilon(\mu')$ est constant.

Appliquant enfin cet argument \`a un corps $F$ cubique, on voit que $\varepsilon(\mu)=1$. Un argument analogue d\'emontre le r\'esultat en $p$, ainsi d'ailleurs qu'en $\ell$ (le pseudo--coefficient de $\St_\ell$ \'etant Whittaker--normalis\'e).

\vskip2mm

\subsection{Facteurs locaux des repr\'esentations cuspidales autoduales}\label{flocauxcusp}

Les arguments de cet article r\'ev\`elent  des propri\'et\'es remarquables des repr\'esentations autoduales, de $\GL(2n)$, 
tant globales (repr\'esentations automorphes) que locales, qui peuvent \^etre comprises du point de vue de la fonctorialit\'e 
entre motifs et repr\'esentations automorphes, ainsi qu'\`a l'aide des r\'esultats conjecturaux d'Arthur.

Soit $\pi$ une repr\'esentation cuspidale autoduale de $\GL(2n,\A)$ o\`u $\A=\A_\Q$. Conjecturalement, 
Langlands lui associe une repr\'esentation complexe, de degr\'e~$2n$~:
\begin{equation*}
r : \mathcal{L}_\Q \lra \GL(2n,\C)\,,
\end{equation*}
(cf. \cite{Lgl}) qui devrait \^etre irr\'eductible. Si $\pi$, donc $r$, est autoduale, il existe sur l'espace de 
$r$ une forme bilin\'eaire non--d\'eg\'en\'er\'ee unique, donc orthogonale ou symplectique, invariante par $r$. 
On dira que $r$ est symplectique ou orthogonale.

Si $\pi_\infty$ est cohomologique, la repr\'esentation $r_\infty$ de $\W_\R$ obtenue 
par restriction est symplectique et non orthogonale, donc $r$ doit \^etre symplectique. Noter que ceci n'est pas vrai si on suppose que 
$\widetilde\pi \cong \pi\otimes \varepsilon$ o\`u $\varepsilon$ est un caract\`ere, m\^eme d'Artin~: les formes de poids 
$k$ impair sur $\GL(2)$ donnent des repr\'esentations $\pi$ (normalis\'ees, comme chez Langlands, de fa\c{c}on \`a \^etre temp\'er\'ees) 
telles que le param\`etre de $\pi_\infty$ est sur~$\C^\times$~:
\[
z\longmapsto
\left(
\begin{array}{ccc}
 (z/\overline z)^{\frac{k-1}{2}}    &   \\
  &(\overline z/z)^{\frac{k-1}{2}}      
\end{array}
\right)
\]
et la repr\'esentation associ\'ee de $\W_\R$ est orthogonale.

S'il existe un nombre premier $p$ tel que $\pi_p$ appartient \`a la s\'erie (essentiellement) discr\`ete, le param\`etre de Langlands
\begin{equation*}
r_p : \WD_{\Q_p} \lra \GL(2n,\C)\,,
\end{equation*}
bien d\'efini d'apr\`es Harris et Taylor, est ind\'ecomposable, et en fait irr\'eductible si on le consid\`ere comme une repr\'esentation
\begin{equation*}
r_p' : \W_{\Q_p} \times {\rm SU}(2) \lra \GL(2n,\C)\,.
\end{equation*}

De nouveau, si $r_p'$ est symplectique (on dira que $\pi_p$ est symplectique), $r$ doit \^etre symplectique. 

\begin{remarque}\label{remarqueGL2} {\rm En dimension $2$, une repr\'esentation irr\'eductible autoduale est symplectique si, 
et seulement si, son d\'eterminant est trivial, ce qui permet de donner un sens non conjectural aux pr\'edictions ci-dessus et de les prouver. 
Soit $\pi$ une cuspidale autoduale de $\GL(2,\AAA_{F})$ qui est discr\`ete en une place $v$, et soit 
$\omega$ le caract\`ere central de $\pi$: alors $\omega=1$ si, et seulement si, $\omega_{v}=1$. En effet, 
$\omega^{2}=1$ donc si $\omega \neq 1$ c'est le signe $\omega_{{E/F}}$ d'une extension quadratique $E$ de $F$, et on a alors
$$\check{\pi} \simeq \pi \otimes \omega_{{E/F}}\circ \det$$ par multiplicit\'e $1$. 
La th\'eorie du changement de base quadratique montre que $\pi$ est l'induite automorphe d'un caract\`ere de Hecke 
$\chi$ de $\AAA_E^{\ast}$  (non isomorphe \`a son
conjugu\'e). On peut alors considerer de mani\`ere non conjecturale le $L$-param\`etre de $\pi$, \`a savoir $L(\pi)={\rm Ind}_{\W_E}^{\W_F} \chi$,
qui est compatible \`a toutes les places avec $\pi$. Comme $\pi_{v}$ est discr\`ete, $L(\pi)_{|\W_{F_{v}}}$ est irr\'eductible, en particulier 
$v$ n'est pas d\'ecompos\'e dans $E$, {\it i.e.} $\omega_{v}=\omega_{E_v/F_v} \neq 1$.}
\end{remarque}

On peut aussi comprendre ces ph\'enom\`enes du point de vue des r\'esultats annonc\'es par Arthur \cite[Ch.~30]{arthurlivre}. 
D'apr\`es ceux-ci, une repr\'esentation $\theta$--stable (disons, cuspidale) de 
$\GL(2n,\A)$ proviendra par fonctorialit\'e 
d'une repr\'esentation automorphe d'un groupe endoscopique $H$ o\`u $H$~est~:
\smallskip
\begin{tabbing}
\noindent (a)\hskip1cm \=$\SO^*(2n+1)$\,, \hskip1cm\=$\widehat H = \Sp(2n,\C)$\\

\noindent (b) \>$\SO^*(2n,\chi)$\,, \>$\widehat H = \SO(2n,\C)$
\end{tabbing}

\noindent o\`u les groupes $\SO^*$ sont quasi--d\'eploy\'es 
(donc d\'eploy\'es en dimension impaire) et o\`u le groupe $\SO^*(2n,\chi)$ est sp\'ecifi\'e par un caract\`ere d'Artin $\chi$ 
d'ordre $2$ d\'ecrit par Arthur \cite{arthurlivre} (et d\'etermin\'e par $\pi$). Des arguments analogues montrent alors que sous les hypoth\`eses 
le type (symplectique ou orthogonal) de $\pi$ sera d\'etermin\'e par celui de $\pi_p$~; 
les r\'esultats d'Arthur devraient rendre ceci inconditionnel.

Les m\'ethodes de cet article, combin\'ees aux r\'esultats d'Harris et Taylor associant 
des repr\'esentations galoisiennes \`a $\pi$ (\cite{HT}, compl\'et\'es par \cite{TY} si 
$\pi_v$ est une repr\'esentation de Steinberg g\'en\'eralis\'ee), permettent toutefois 
d'obtenir le r\'esultat suivant.

\begin{theoreme}\label{thmsymplectique} Soient $F$ un corps totalement r\'eel et $\pi$ une
repr\'esentation automorphe cuspidale de $\GL_{2n}(\AAA_F)$. On suppose que
$\pi$ est autoduale, essentiellement de carr\'e int\'egrable en au moins une
place finie $w$, et cohomologique \`a toutes les places archim\'ediennes. Alors : 
\begin{itemize}
\item[(i)] Pour toute place $v$ de $F$, $\pi_v$ est symplectique. 
\item[(ii)] Si $V_\ell$ est une repr\'esentation $\ell$-adique de
$\Gal(\overline{F}/F)$ associ\'ee \`a $\pi|\cdot|^{(2n-1)/2}$, 
avec $w$ premier \`a $\ell$, alors il existe un accouplement symplectique non d\'eg\'en\'er\'e et Galois-equivariant 
		$$V_\ell \otimes V_\ell \longrightarrow\,
{\overline{\Q}}_\ell\,\,(2n-1).$$
\end{itemize}
\end{theoreme}

\begin{remarque}\label{remarqueunpeugrosse} {\rm 
L'hypoth\`ese "$w$ premier \`a $\ell$" dans (ii) ne serait pas
n\'ecessaire s'il on disposait, pour les repr\'esentations galoisiennes
d'Harris et Taylor, d'un analogue du th\'eor\`eme de Saito identifiant en
$w \,|\,\ell$ la repr\'esentation de Weil-Deligne donn\'ee par la th\'eorie de
Fontaine.}
\end{remarque}

Si $v$ est finie, \og $\pi_v$ symplectique\fg\, signifie que son param\`etre de Langlands sous sa forme ${\rm SU}(2)$ pr\'eserve une forme bilin\'eaire altern\'ee non d\'eg\'en\'eree. Du point de vue de la repr\'esentation de Weil-Deligne $(r,N)$ associ\'ee, il est \'equivalent de demander que $r$ (la repr\'esentation de $\W_{F_v}$) pr\'eserve une forme bilin\'eaire altern\'ee non d\'eg\'en\'er\'ee, et que $N$ est dans l'alg\`ebre de Lie du groupe symplectique associ\'e.

Dans le contexte de l'\'enonc\'e, comme $\pi_v$ est temp\'er\'ee (\cite{HT}), elle s'\'ecrira alors 
\begin{equation}\label{induitegenerique}
\pi_v ={ \rm ind}_P^G(\delta _1, \delta _1^\theta,\ldots \delta _r, \delta _r^\theta,\ 
\delta _{2r+1},\ldots,\delta _s) 
\end{equation}
\noindent o\`u les $\delta _i$ sont des repr\'esentations de la s\'erie (essentiellement) discr\`ete de $\GL(n_i,F_v)$, 
et les repr\'esenta\-tions $\delta _{2r+1},\ldots,\delta _s$ sont autoduales, distinctes et symplectiques. 

D\'emontrons maintenant le th\'eor\`eme~\ref{thmsymplectique}. D'apr\`es le \S\ref{prelimglnreel}, le (i) vaut par hypoth\`ese si $v$ est archim\'edienne. 

Fixons $v$ finie, choisissons un $\ell$ premier \`a $w$ et $v$, puis consid\'erons une repr\'e\-sen\-ta\-tion galoisienne $\ell$-adique $V_\ell$ attach\'ee\footnote{L'existence et les propri\'et\'es de compatibilit\'e de $V_\ell$ se d\'eduisent de mani\`ere standard de \cite{HT} par changement de base \`a tous les $EF$, $E$ \'etant un corps quadratique imaginaire tel que $EF$ est d\'ecompos\'e au dessus de $w$ (voir \cite[Thm 3.6]{Tay} pour le cas $F=\Q$).}  \`a $\pi|.|^{(2n-1)/2}$. D'apr\`es \cite{HT} et
\cite{TY}, $V_\ell$ est irr\'eductible, compatible avec $\pi$ (en toutes les places finies premi\`eres \`a $\ell$) selon la correspondance de Langlands Frobenius semi-simplifi\'ee ; en particulier, cela nous fournit un accouplement non d\'eg\'en\'er\'e Galois-\'equivariant comme dans l'\'enonc\'e (unique \`a un scalaire pr\`es, et dont il faut montrer qu'il est symplectique). Ainsi\footnote{On pourra remarquer que si une repr\'esentation de Weil-Deligne pr\'eserve un accouplement symplectique, alors sa Frobenius simplification (qui associe \`a chaque \'el\'ement du groupe de Weil sa partie semisimple dans sa d\'ecomposition de Jordan) pr\'eserve le m\^eme accouplement.}, il suffit de prouver (i) quand $v=w$ pour l'avoir pour tout $v$, et (i) implique (ii). De plus, si $\pi_w$ est
la repr\'esentation de Steinberg, alors (i) et (ii) sont satisfaits. Il suffit donc de d\'emontrer le th\'eor\`eme suivant (en fait, (i) $\Rightarrow$ (ii) $\Rightarrow$ (v) suffit).

\begin{theoreme}\label{propsymplstein} Soient $\pi_w$ une repr\'esentation autoduale (essentiellement) discr\`ete de $\GL_{2n}(F_w)$ et $f_w$ un pseudocoefficient tordu de $\pi_w$. Les conditions suivantes sont \'equivalentes : \begin{itemize}
\item[(i)] $\pi_w$ est composante locale d'une repr\'esentation automorphe cuspidale $\pi$ de $\GL_{2n}(\AAA_F)$ qui est autoduale, et cohomologique \`a toutes les places ar\-chi\-m\'e\-dien\-nes, 
\item[(ii)] $\TO_{\gamma_0}(f_w) \neq 0$,
\item[(iii)] les int\'egrales orbitales tordues stables de $f_w$ ne sont pas identiquement nulles, 
\item[(iv)] les int\'egrales orbitales tordues stables fortement r\'eguli\`eres de $f_w$ ne sont pas identiquement nulles, 
\item[(v)] $\pi_w$ est composante locale d'une repr\'esentation automorphe cuspidale $\pi$ comme en (i) et qui de plus est la repr\'esentation de Steinberg \`a une autre place finie. \par \smallskip
\end{itemize} 

Si l'une de ces propri\'et\'es est satisfaite, alors $\pi_w$ est symplectique. 

\end{theoreme}

\begin{pf} L'implication (i) $\Rightarrow$ (ii)  d\'ecoule du m\^eme raisonnement de positivit\'e que dans le \S\ref{argumentpositivite}. 
Quitte \`a faire un changement de base quadratique r\'eel d\'ecompos\'e en $w$, on peut supposer que $F \neq \Q$. Pour toutes les places 
finies $v \neq w$ de $F$, choisissons un pseudocoefficient positif $f_v$ de $\pi_v$ comme dans le Lemme~\ref{consPWR}. Enfin, pour chaque 
place archim\'edienne $v$ de $F$ on prend pour $f_v$ un pseudocoefficient de s\'eries th\'eta-discr\`etes cohomologiques quelconque comme 
au \S\ref{pseudoorb} ; il d\'epend d'un $[F:\Q]$-uple $\lambda$ de poids dominants de $\SO_{2n+1}(\R)$. En appliquant la formule des traces 
tordue d'Arthur \`a la fonction $\prod_v' f_v$ (qui est simplifi\'ee car cette fonction est cuspidale en deux places archim\'ediennes 
au moins) et en faisant varier $\lambda$, l'argument du \S\ref{argumentpositivite} montre que 
\begin{equation}\label{mainpropff} \forall v, \, \, \, \, \, 
\TO_{\delta_0}(f_v) \neq 0,
\end{equation}
donc (ii). 

L'implication (ii) $\Rightarrow$ (v) est par exemple cons\'equence de la m\'ethode du \S\ref{termeprincal}, si l'on choisit la fonction test \'egale \`a $f_w$ en la place $w$ : l'appel au Th\'eor\`eme~\ref{thmcle} est remplac\'e par (\ref{mainpropff}).

(ii) $\Rightarrow$ (iii) est \'evident. Rappelons que les int\'egrales orbitales stables 
sont d\'efinies en tous les \'el\'ements $\theta-$semisimples, et sont limites faibles de telles int\'egrales pour des \'el\'ements fortement $\theta$--r\'eguliers, donc (iii) $\Leftrightarrow$ (iv).

\`A partir de (iv), les arguments de la Proposition~\ref{ancienprop48} d\'emontrent (v). 

	Enfin, d'apr\`es les consid\'erations pr\'ec\'edentes, le
Th\'eor\`eme \ref{thmsymplectique} s'ensuit, ainsi donc que la derni\`ere
propri\'et\'e. \end{pf}

\begin{remarque}{\rm  En fait, ces propri\'et\'es sont sans doute {\it \'equivalentes} \`a \og $\pi_w$ est symplectique \fg. 
C'est connu si $\pi_w$ est supercuspidale par les r\'esultats de Henniart
\cite{He} et Shahidi \cite[Prop. 5.1]{Shahidi2} rappel\'es dans l'introduction.
Nous reviendrons sur ce point dans un travail ult\'erieur.}\end{remarque}

 Notons enfin le ph\'enom\`ene purement local qui devrait r\'esulter de nos m\'ethodes. Consid\'erons, sur $G(\Q_p)$, 
 l'int\'egrale orbitale tordue $\TO_{\delta _0}$ (on utilise l'automorphisme $\theta_0$ Whittaker--normalis\'e). 
 D'apr\`es des principes g\'en\'eraux, elle devrait s'\'ecrire
 \begin{equation}\label{planchereltheta}
\TO_{\delta _0}(f) =\displaystyle \int_{{\widehat \G}_{\rm temp}^{\theta_0}} {\rm trace}(A_{\theta_0} \pi(f))\ d\mu(\pi) 
\end{equation}
o\`u ${\widehat G}_{\rm temp}^{\theta_0}$ est la vari\'et\'e des repr\'esentations temp\'er\'ees, $\theta_0$--stables de $G$, 
et $d\mu(\pi)$ une distribution (une mesure ?). La distribution $\TO_{\delta _0}(f)$ est stable~; 
d'apr\`es une extension simple du Th\'eor\`eme~\ref{propsymplstein}, les int\'egrales orbitales tordues stables de pseudo--coefficients 
(au sens du Lemme~\ref{consPWR}) d'une famille d'induites non symplectiques --- i.e., ne v\'erifiant pas la condition 
suivant (\ref{induitegenerique}) --- s'annulent. On dispose de plus de l'argument de positivit\'e du \S\ref{argumentpositivite}. 
Par cons\'equent, on s'attend au r\'esultat suivant, sans doute accessible~:

\begin{conj}\label{conjectureplancherel}

\noindent{\rm (i)} Dans la formule $(\ref{planchereltheta})$, l'int\'egrale ne porte que sur les re\-pr\'e\-sen\-ta\-tions symplectiques.

\noindent{\rm (ii)} La distribution $d\mu(\pi)$ est une mesure positive.

\end{conj}

Terminons en remarquant que de telles propri\'et\'es feront sans doute partie des d\'emonstrations d'Arthur~!
 
\bigskip
\bigskip
\section{Application au groupe de Galois absolu de $\Q$}
\setcounter{equation}{0}

Repla\c{c}ons nous dans le contexte de l'introduction. Soient $S$ un ensemble fini de nombres premiers, 
$\Q_S$ une extension alg\'ebrique maximale de $\Q$ non ramifi\'ee hors de $S$ (et de l'infini) et $\G_S=\Gal(\Q_S/\Q)$ son groupe de Galois. 
Supposons $S\neq \emptyset$ et fixons $p \in S$ un nombre premier.

\begin{theoreme}\label{applicationgalois} Si $|S|\geq 2$, les applications naturelles $\Gal(\Qpb/\Q_p) \longrightarrow G_S$ sont 
injectives.
\end{theoreme}

\begin{pf} On peut supposer que $S=\{\ell,p\}$ o\`u $\ell$ est un nombre
premier diff\'erent de $p$. Fixons un plongement $\Q_S \longrightarrow \Qpb$, et notons $$E:=\Q_p.\Q_S$$ 
le compositum de $\Q_p$ et $\Q_S$ dans $\Qpb$ : 
il s'agit de d\'emontrer que $E=\Qpb$. Comme $\Gal(\Qpb/\Q_p)$ n'admet pas de sous-groupe ferm\'e normal non trivial rencontrant 
trivialement le sous-groupe d'inertie par \cite[Lemme 4 (i)]{Ch}, il suffit de d\'emontrer que $\Q_p^{\rm nr}.E=\Qpb$, soit encore que 
$\Q_p^{\rm ab}.E=\Qpb$ puisqu'il est \'evident que $\Q_p(\mu_{p^{\infty}}) \subset E$.

D'apr\`es \cite[Lemme 4 (iii)]{Ch}, la cl\^oture alg\'ebrique s\'eparable d'un corps $k$ est le compositum de ses sous-extensions galoisiennes finies $K/k$ dont le groupe de Galois
admet une representation lin\'eaire complexe \`a la fois injective et irr\'eductible. Fixons donc une telle extension $K$ de 
$k:=\Q_p$, il s'agit de montrer que $K \subset E.\Q_p^{\rm ab}$.

Vu le choix de $K$, il existe une repr\'esentation injective irr\'eductible 
$\rho: \Gal(K/\Q_p) \rightarrow \GL_n(\C)$ pour un certain entier $n\geq 1$, que l'on voit par inflation comme une repr\'esentation continue irr\'eductible
$$\psi: \W_{\Qp} \longrightarrow \GL_n(\C)$$
du groupe de Weil $\W_{\Qp}$ de $\Qp$. Quitte \`a tordre $\psi$ par un caract\`ere lisse $\eta$ bien choisi de $\W_{\Q_p}$, nous pouvons supposer que la repr\'esentation duale $\psi^*$ n'est 
isomorphe \`a aucune torsion de $\psi$ par un caract\`ere non ramifi\'e. En effet, si $I \subset \W_{\Q_p}$ est le sous-groupe d'inertie, 
$I \cap \Ker \psi$ est un sous-groupe ouvert de $I$. On peut donc trouver un \'el\'ement $g \in I \cap \Ker \psi$ 
agissant sur $\Q_p(\mu_{p^{\infty}})$ par un automorphisme d'ordre infini. Il suffit alors de choisir un caract\`ere 
$\eta: \W_{\Qp} \rightarrow \Gal(\Q_p(\mu_{p^{\infty}})/\Q_p) \rightarrow \C^*$ tel que $\eta(g)^2 \neq 1$. Par construction et injectivit\'e de $\rho$, 
notons que l'extension de $\Qp^{\rm ab}$ fix\'ee par le noyau la restriction de $\psi$ \`a $\W_{\Q_p}\cap \Gal(\Qpb/\Qp^{\rm ab})$
est $K.\Qp^{\rm ab}$.

Soit $\omega$ la repr\'esentation supercuspidale de $\GL_n(\Q_p)$ de $L$-param\`etre $\psi$ donn\'ee par la correspondance de Langlands locale (\cite{HT}). 
Par construction, sa contragr\'ediente $\check{\omega}$ n'est isomorphe \`a aucune tordue de
$\omega$ par un caract\`ere non ramifi\'e. 
D'apr\`es le Th\'eor\`eme~\ref{mainthm}, il existe une repr\'esentation automorphe cuspidale irr\'eductible $\Pi$ de $\GL_{2n}(\AAA_{\Q})$ ayant les propri\'et\'es suivantes:\begin{itemize}
\item[i)] $\check{\Pi} \simeq \Pi$,
\item[ii)] $\Pi$ est non ramifi\'ee aux places finies diff\'erentes de $p$ et $\ell$, 
\item[iii)] $\Pi_p \simeq I(\chi)$ pour un certain caract\`ere non ramifi\'e
$\chi: \Q_p^* \longrightarrow \C^*$ (cf. \S\ref{ssectenonce},\S\ref{pseudocoeffcpomega})
\item[iv)] $\Pi_{\infty}|\cdot|^{(2n-1)/2}$ est alg\'ebrique r\'eguli\`ere et
$\Pi_\ell$ est la repr\'esentation de Steinberg.
\end{itemize}

Sous les conditions i) et iv), les travaux de Kottwitz, Clozel et Harris-Taylor montrent qu'un choix quelconque de 
plongements $\ell$-adiques et complexes de $\overline{\Q}$ \'etant fait, on peut associer \`a $\Pi|\cdot|^{(2n-1)/2}$ une repr\'esentation continue
$$\rho_{\Pi,\ell}: \G_S \longrightarrow \GL_{2n}(\overline{\Q}_\ell)$$
compatible \`a toutes les places finies $\neq \ell$ \`a la correspondance de Langlands locale \og Frobenius semi-simplifi\'ee\fg\,  
(\cite{HT}, \cite[Thm 3.6]{Tay}).

Que cette repr\'esentation se factorise par $\G_S$ d\'ecoule alors bien s\^ur de la condition ii). Le $L$-param\`etre de
$\Pi_p \simeq I(\chi)$ est la somme directe $\psi \otimes \chi \oplus  \psi^*
\otimes \chi^{-1}$
(cf. \S\ref{pseudocoeffcpomega}). La compatibilit\'e de ${\rho_{\Pi,\ell}}_{|W_{\Q_p}}$ \`a $\Pi_p|\cdot|^{(2n-1)/2}$ assure en particulier que si $\Q(\Pi) \subset \Q_S$ est 
le sous-corps de $\Qb$ fix\'e par $\Ker(\rho_{\Pi,\ell})$, alors $$\Q(\Pi).\Q_p^{\rm ab}=K.\Qp^{\rm ab} \subset E.\Qp^{\rm ab},$$ ce qui conclut.
\end{pf}

Puisqu'il existe une surjection continue de $\Gal(\Qpb/\Q_p)$ vers $\widehat{\Z}$, le th\'eor\`eme ci-dessus admet le corollaire suivant 
en direction d'une question de J. Milne.

\begin{corollaire}\label{corapplicationgalois} Si $|S| \geq 2$, alors $|G_S|=|\widehat{\Z}|$. \end{corollaire}

Autrement dit: si $|S|\geq 2$, pour tout entier $n\geq 1$ il existe un corps 
de nombres de degr\'e divisible par $n$ qui est non ramifi\'e hors de $S$.
Notons que ce corollaire admet par exemple le corollaire imm\'ediat suivant: si $F$ est un corps
local et $d\geq 1$, il n'existe pas de repr\'esentation continue injective $G_S \rightarrow
\GL_d(F)$.

\begin{remarque}{\rm \begin{itemize}
\item[(i)] Les r\'esultats de cette section admettent des
g\'en\'eralisations imm\'ediates, avec les m\^emes preuves, dans le cas o\`u le corps de base $\Q$ est remplac\'e par un
corps totalement r\'eel $F$ quelconque. L'\'enonc\'e est alors le suivant:
soit $S$ un ensemble fini non vide de places finies de $F$, $v \in S$, et
supposons que $S \backslash \{v\}$ contient toutes les places divisant un
nombre premier $\ell$ donn\'e, alors les applications naturelles
$$\Gal(\overline{F_v}/F_v) \longrightarrow \Gal(F_S/F)$$
sont injectives.
\item[(ii)] Une cons\'equence du th\'eor\`eme est que l'image de $\Q_p.\Q_S$ dans $\Qpb$
contient l'extension maximale non ramifi\'ee de $\Q_p$. Il est remarquable
que par la m\'ethode employ\'ee ici ce simple fait ne semble pas pouvoir s'obtenir \`a
beaucoup moins de frais que le r\'esultat tout entier. Le point est que nous ne pouvons pas prescrire
totalement $\Pi_p$ pour les $\Pi$ que nous construisons mais seulement sa classe inertielle, et c'est bien \'evident : 
nous ne pouvons pas
deviner \`a priori les nombres de Weil attach\'es au $L$-param\`etre de
$\Pi_p$ (\`a savoir le $\chi$ tel que $\Pi_p=I(\chi)$). Du coup,
nous ne contr\^olons jamais vraiment la partie non ramifi\'ee des corps de
nombres que nous contruisons, bien que nous maitrisions parfaitement leurs
groupes d'inerties. Celle-l\`a n'est r\'ecup\'er\'ee qu'\`a la
fin \`a cause de la structure du groupe de Galois local.
\item[(iii)] La question de savoir si le th\'eor\`eme vaut si $S=\{p\}$
semble nettement plus d\'elicate (voir \cite[\S 4.2]{Ch}). Il est cependant tentant de
conjecturer que le th\'eor\`eme est aussi vrai dans ce cas.
\end{itemize}}
\end{remarque}


\begin{thebibliography}{9999}


\bibitem[A1]{ITF1} J. Arthur, {\it The invariant trace formula I.
Local theory}, J. Amer. Math. Soc. {\bf 1} (1988), 323--383.

\bibitem[A2]{ITF} J. Arthur, {\it The invariant trace formula II. Global
theory}, J. Amer. Math. Soc. {\bf 1} (1988), 501--554.

\bibitem[A3]{ArD} J. Arthur, {\it The local behaviour of weighted orbital integrals}, Duke Math. J. {\bf 56} (1988), 223--293.

\bibitem[A4]{arthurlivre} J. Arthur, {\it An introduction to the trace formula}, 
in {\it Harmonic Analysis, the trace formula and Shimura varieties}, A.M.S.\,, Clay Math. Institute (2005).

\bibitem[AC]{AC} J. Arthur \& L. Clozel, {\it Simple algebras, base change, and the advanced 
theory of the trace formula}, Ann. of Math. Studies {\bf 120} (1989). 

\bibitem[Be]{bernstein} J.-N. Bernstein (r\'edig\'e par P. Deligne), {\it Le centre de Bernstein}, in 
Bernstein, Deligne, Kazhdan, Vign\'eras, {\it Repr\'esentations des groupes r\'eductifs sur un corps local}, Hermann (1984).

\bibitem[BLS]{BLS} A.Borel, J.-P. Labesse \& J. Schwermer, 
{\it On the cuspidal cohomology of $S$-arithmetic subgroups of reductive groups over number fields}, 
Compositio Math. {\bf 102} (1996), 1--40. 

\bibitem[BW]{BW} A. Borel \& N. Wallach, {\it Continuous cohomology, discrete subgroups, and representations of reductive groups},  
Annals of Mathematics Studies {\bf 94} (1980).

\bibitem[Bou]{bouaziz} A. Bouaziz, {\it Sur les caract\`eres des groupes de Lie r\'eductifs non connexes}, J. Funct. Analysis {\bf 70} (1987), 1--79.

\bibitem[BT]{BT} F. Bruhat \& J. Tits, {\it Groupes r\'eductifs sur un corps local: I. Donn\'ees radicielles valu\'ees}, Pub. Math. IHES 
{\bf 41} (1972).

\bibitem[Cas]{Cas} W. Casselman, {\it A new nonunitarity argument for $p$-adic representations}, 
  J. Fac. Sci. Univ. Tokyo {\bf 28} (1981), 907--928.

\bibitem[Ch]{Ch} G. Chenevier, {\it On number fields with given ramification}. \`A para\^itre \`a Compositio Math.\, . 

\bibitem[Clo1]{ClAB} L. Clozel, {\it Th\'eor\`eme d'Atiyah-Bott pour les vari\'et\'es $p$-adiques et 
caract\`eres des groupes r\'eductifs}, Harmonic analysis on Lie groups and symmetric spaces (Kleebach, 1983), 
M\'em. Soc. Math. France {\bf 15} (1984), 39--64.

\bibitem[Clo2]{AA} L. Clozel, {\it Motifs et formes automorphes : applications du principe de fonctorialit\'e}, in 
{\it Automorphic forms, Shimura varieties, and $L$-function}, Clozel, Milne eds. vol. I, Perspectives in Math. {\bf 10}, Academic Press (1990), 77--159.

\bibitem[CloD]{CloDel} L. Clozel \& P. Delorme, {\it Le th\'eor\`eme de Paley-Wiener invariant 
pour les groupes de Lie r\'eductifs II},  Ann. Sci. \'Ecole Norm. Sup. {\bf 23}  (1990), 
193--228. 

\bibitem[HT]{HT} M. Harris \& R. Taylor, {\it The geometry and cohomology of some simple Shimura varieties},
Annals of Math. Studies {\bf 151} (2001).

\bibitem[He1]{henniartftr} G. Henniart, {\it La conjecture de Langlands locale pour $\GL(3)$}, M\'em. S.M.F.\,, no. 11--12 (1984).

\bibitem[He2]{He} G. Henniart, {\it Correspondance de Langlands et fonctions
$L$ des carr\'es ext\'erieur et sym\'etrique}. Pr\'epublications math\'ematiques de l'IHES (2003).

\bibitem[Hi]{Hi} A. Hitta, {\it On the continuous (co)homology of locally profinite groups and the
K\"unneth theorem}, J. Algebra  {\bf 163}  (1994),  481--494. 

\bibitem[K1]{Ko} R. Kottwitz, {\it Tamagawa numbers}, Ann. of Math. {\bf 127} (1988), 629--646. 

\bibitem[K2]{kott} R. Kottwitz, {\it Rational conjugacy classes in reductive groups}, Duke Math. J. {\bf 49} (1982), 785--806.

\bibitem[KR]{KoR} R. Kottwitz \& J. Rogawski, {\it The distributions in the invariant trace formula are supported on
characters}, Canad. J. Math. {\bf 52} (2000), 804--814. 

\bibitem[KS]{KS} R. Kottwitz \& D. Shelstad, {\it Foundations of twisted endoscopy}, Ast\'erisque  {\bf 255}  (1999). 

\bibitem[La1]{labesse} J.-P. Labesse, {\it Cohomologie, stabilisation et
changement de base}, Ast\'erisque {\bf 257} (1999).

\bibitem[La2]{Labch1} J.-P. Labesse, {\it Pseudo-coefficients tr\`es cuspidaux et $K$-th\'eorie}, Math. Ann. {\bf 291} (1991), 607--616.

\bibitem[Lg]{Lgl} R. P. Langlands, {\it Ein M\"archen}, in {\it Automorphic representations, Shimura varieties, and motives}, in Proc. Sympos. Pure Math. 33 vol. II, A.M.S. , 
Providence (1979), 205--246. 

\bibitem[M]{Mezo} P. Mezo, {\it Twisted trace Paley-Wiener theorems for special and 
general linear groups}, Compositio Math. {\bf 140} (2004),  205--227.

\bibitem[PR]{PR} D. Prasad \& D. Ramakrishnan, {\it On self-dual
repr\'esentations of division algebras over local fields}. Pr\'epublication.

\bibitem[Ro]{paleyweinerrog} J. Rogawski, {\it Trace Paley-Wiener theorem in the twisted case}, Trans. A.M.S. {\bf 309} (1988), 215--229.

\bibitem[S1]{Serre}  J.-P. Serre, {\it Cohomologie des groupes discrets}, 
Prospects in mathematics, Ann. of Math. Studies {\bf 70} (1971), 77--169.

\bibitem[S2]{SerreJams} J.-P. Serre, {\it R\'epartition asymptotique des
valeurs propres de l'op\'erateur de Hecke $T\sb p$}, J. Amer. Math. Soc. {\bf 10} (1997),  75--102. 

\bibitem[Sh1]{Shahidi1} F. Shahidi, {\it A proof of Langlands' conjecture on Plancherel measures;
complementary series for $p$-adic groups},  Ann. of Math. {\bf 132} (1990), 273--330.

\bibitem[Sh2]{Shahidi2} F. Shahidi, {\it Twisted endoscopy and reducibility
of induced representations for $p$-adic groups}, Duke Math. J.  {\bf 66}  (1992), 1--41. 

\bibitem[She]{Shelstadlemme} D. Shelstad, {\it Characters and inner forms of a quasi-split group over $\R$}, Comp. Math. {\bf 39} (1979), 11--45.

\bibitem[T]{Tay} R. Taylor, {\it Galois representations}, Annales de la Facult\'e des Sciences de Toulouse {\bf 13} 
(2004), 73--119.

\bibitem[TY]{TY} R. Taylor \& T. Yoshida, {\it Compatibility of local and global Langlands correspondences}, J. Amer. Math. Soc. {\bf 20} (2007), 467--493. 

\bibitem[Ti]{Ti} J. Tits, {\it Reductive groups over local fields}, in 
{\it Automorphic forms, representations and $L$-functions}, Proc. Sympos. Pure Math. Part 1, Corvallis (1977), 29--69.

\bibitem[V]{V} D. Vogan, {\it Unitarizability of certain series of representations}, Ann. of Math. {\bf 120} (1984), 141--187.

\bibitem[W1]{Wald} J.-L. Waldspurger, {\it Le groupe $\GL(n)$ tordu sur un
corps $p$-adique, partie I}. Pr\'epublication.

\bibitem[W2]{Wald2} J.-L. Waldspurger, {\it Le groupe $\GL(n)$ tordu sur un
corps $p$-adique, partie II}. Pr\'epublication.
 


\end{thebibliography}
\end{document}